\documentclass{amsproc}
\usepackage{amsxtra}
\usepackage{amsmath}
\usepackage{graphicx}
\usepackage{latexsym}
\usepackage{amsfonts}
\usepackage{amssymb}
\theoremstyle{definition}
\newtheorem{definition}{Definition}[section]
\newtheorem{remark}[definition]{Remark}
\newtheorem{example}[definition]{Example}

\theoremstyle{plain}
\newtheorem{theorem}[definition]{Theorem}
\newtheorem{lemma}[definition]{Lemma}
\newtheorem{corollary}[definition]{Corollary}
\newtheorem{proposition}[definition]{Proposition}
\theoremstyle{remark}
\newtheorem*{acknowledgements}{Acknowledgements}
\newbox\ipbox
\newcommand{\ip}[2]{\left\langle #1\mathrel{\mathchoice
{\setbox\ipbox=\hbox{$\displaystyle \left\langle\mathstrut #1#2\right\rangle$}
\vrule height\ht\ipbox width0.25pt depth\dp\ipbox}
{\setbox\ipbox=\hbox{$\textstyle \left\langle\mathstrut #1#2\right\rangle$}
\vrule height\ht\ipbox width0.25pt depth\dp\ipbox}
{\setbox\ipbox=\hbox{$\scriptstyle \left\langle\mathstrut #1#2\right\rangle$}
\vrule height\ht\ipbox width0.25pt depth\dp\ipbox}
{\setbox\ipbox=\hbox{$\scriptscriptstyle \left\langle\mathstrut #1#2\right\rangle$}
\vrule height\ht\ipbox width0.25pt depth\dp\ipbox}
} #2\right\rangle}

\def\openone
{\mathchoice
{\hbox{\upshape \small1\kern-3.3pt\normalsize1}}
{\hbox{\upshape \small1\kern-3.3pt\normalsize1}}
{\hbox{\upshape \tiny1\kern-2.3pt\SMALL1}}
{\hbox{\upshape \Tiny1\kern-2pt\tiny1}}}

\newlength{\customskipamount}
\newlength{\customleftmargin}
\newlength{\qedskip}

\input cyracc.def
\font\eightcyr=wncyr8

\renewcommand{\emptyset}{\varnothing}

\def\hooklongrightarrow{\lhook\joinrel\longrightarrow}

\begin{document}
\title[Unitary Representations and O-S Duality]{Unitary Representations and Osterwalder-Schrader Duality}
\author{Palle E. T. Jorgensen}
\address{Department of Mathematics\\
University of Iowa\\
Iowa City\\
IA 52242}
\email{jorgen@math.uiowa.edu}
\author{Gestur \'Olafsson}
\address{Department of Mathematics\\
Louisiana State University\\
Baton Rouge\\
LA 70803}
\email{olafsson@math.lsu.edu}
\thanks{Both authors are supported in part by the U.S. National Science Foundation.}
\thanks{The second author was supported by LEQSF grant
(1996-99)-RD-A-12.}
\dedicatory{Dedicated to the memory of Irving E. Segal}

\begin{abstract}
The notions of reflection, symmetry, and positivity from quantum field theory
are shown to induce a duality operation for a general class of unitary
representations of Lie groups. The semisimple Lie groups which have this
$c$-duality are identified and they are placed in the context of
Harish-Chandra's legacy for the unitary representations program. Our paper
begins with a discussion of path space measures, which is the setting where
reflection positivity (Osterwalder-Schrader) was first identified as a useful
tool of analysis.
\end{abstract}\maketitle

\begin{center}
\textit{Le plus court chemin entre deux v\'{e}rit\'{e}s dans}\linebreak
\textit{le domaine r\'{e}el passe par le domaine complexe.\bigskip}%
\linebreak \settowidth{\qedskip}{\textit{le domaine r{\'{e}}%
el passe par le domaine complexe.}}\makebox[\qedskip]{\hfill\textsc
{---Jacques Hadamard}}\vspace{24pt}
\end{center}

\section*{\label{Introduction}Introduction}

In this paper, we present an idea which serves to unify the following six
different developments:

\begin{enumerate}
\item \label{IntCat(1)}reflection positivity of quantum field theory,

\item \label{IntCat(2)}the role of reflection positivity in functional integration,

\item \label{IntCat(3)}the spectral theory of unitary one-parameter groups in
Hilbert space,

\item \label{IntCat(4)}an extension principle for operators in Hilbert space,

\item \label{IntCat(5)}the Bargmann transform, and

\item \label{IntCat(6)}reflection positivity and unitary highest weight
modules for semisimple Lie groups.
\end{enumerate}

The emphasis is on (\ref{IntCat(6)}), but the common thread in our paper is
the unity of the six areas, which, on the face of it, may perhaps appear to be
unrelated. We also stress the interconnection between the six projects, and
especially the impact on (\ref{IntCat(6)}) from (\ref{IntCat(1)}%
)--(\ref{IntCat(5)}).

Readers who may not be familiar with all six of the subjects (\ref{IntCat(1)}%
)--(\ref{IntCat(6)}) may wish to consult the bibliography; for example
\cite{GlJa87} is an excellent background reference on (\ref{IntCat(1)}%
)--(\ref{IntCat(2)}), \cite{FOS83} is especially useful on (\ref{IntCat(2)}),
and \cite{Arv84} covers the theory underlying (\ref{IntCat(3)}). Area
(\ref{IntCat(4)}) is covered in \cite{Jor80} and \cite{JoMu80}, while
\cite{BH94} and \cite{OO96} cover (\ref{IntCat(5)}). Background references on
(\ref{IntCat(6)}) include \cite{JO97}, \cite{Jor86}, and \cite{Jor87}.

Our main point is to show how the concepts of reflection, symmetry and
positivity, which are central notions in quantum field theory, are related to
a duality operation for certain unitary representations of semisimple Lie
groups. On the group level this duality corresponds to the $c$-duality of
causal symmetric spaces, a duality relating the ``compactly causal'' spaces to
the ``non-compactly causal'' spaces. On the level of representations one
starts with a unitary representation of a group $G$ with an involution $\tau$
corresponding to a non-compactly causal space, then uses an involution on the
representation space satisfying a certain positivity condition on a subspace
to produce a contraction representation of a semigroup $H\exp\left(  C\right)
$, where $H=G^{\tau}$ and $C$ is an $H$-invariant convex cone lying in the
space of $\tau$-fixed points in the Lie algebra. Now a general result of
L\"uscher and Mack and the co-authors can be used to produce a unitary
representation of the $c$-dual group $G^{c}$ on the same space. (See
\cite{HO95}, \cite{HiNe93},
\cite{LM75}, and \cite{O98} for these terms.)

We aim to address several target audiences: workers in representation theory,
mathematical physicists, and specialists in transform theory. This diversity
has necessitated the inclusion of a bit more background material than would
perhaps otherwise be called for: certain ideas are typically explained
slightly differently in the context of mathematical physics from what is
customary among specialists in one or more of the other areas.

The symmetry group for classical mechanics is the \textit{Euclidean motion
group} $E_{n}=SO(n)\times_{\mathrm{sp}}\mathbb{R}^{n}$, where the subscript
${}_{\mathrm{sp}}$ stands for semidirect product. Here the action of
$(A,\mathbf{x})\in E_{n}$ on $\mathbb{R}^{n}$ is given by $(A,\mathbf{x}%
)\cdot\mathbf{v}=A(\mathbf{v})+\mathbf{x}$, that is $SO(n)$ acts by rotations
and $\mathbb{R}^{n}$ acts by translations. The connected symmetry group for
the space-time of relativity is the \textit{Poincar\'{e} group} $P_{n}%
=SO_{o}(n-1,1)\otimes_{sp}\mathbb{R}^{n}$. Let $x_{n}$ stand for the time
coordinate, that is $t=x_{n}$. Those two symmetry groups of physics are
related by transition to imaginary time, that is multiplying $x_{n}$ by $i$.
This in particular changes the usual Euclidean form $(\mathbf{x}\mid
\mathbf{y})=x_{1}y_{n}+\cdots+x_{n}y_{n}$ into the Lorentz form $q_{n-1,1}%
(\mathbf{x},\mathbf{y})=x_{1}y_{1}+\cdots+x_{n-1}y_{n-1}-x_{n}y_{n}$, which is
invariant under the group $SO_{o}(n-1,1)$. Those groups and the corresponding
geometry can be related by the $c$-duality. Define an involution $\tau\colon
E_{n}\rightarrow E_{n}$ by
\[
\tau(A,\mathbf{x})=(I_{n-1,1}AI_{n-1,1},I_{n-1,1}\mathbf{x})\,,\qquad
I_{n-1,1}=\left(
\begin{array}
[c]{cc}%
I_{n-1} & 0\\
0 & -1
\end{array}
\right)  \,.
\]
The differential $\tau\colon\frak{e}_{n}\rightarrow\frak{e}_{n}$ is given by
the same formula, and $\tau$ is an involution on $\frak{{e}_{n}}$. Let
\begin{align*}
\frak{h}:=  &  \{X\in\frak{e}_{n}\mid\tau(X)=X\}\simeq\frak{s}\frak{o}%
(n-1)\times_{sp}\mathbb{R}^{n-1}\\
\frak{q}:=  &  \{X\in\frak{e}_{n}\mid\tau(X)=-X\}
\end{align*}
The $c$-dual Lie algebra $\frak{e}^{c}$ is defined by:
\[
\frak{e}_{n}^{c}:=\frak{h}\oplus i\frak{q}\,.
\]
A simple calculation shows that $\frak{e}_{n}\simeq\frak{p}_{n}$. Let
$G=E_{n}$ and let $G^{c}$ denote the simply connected Lie group with Lie
algebra $\frak{{e}^{c}}$. Then $G^{c}=\tilde{P}_{n}$, the universal covering
group of $P_{n}$.

Given that a physical system is determined by a \textit{unitary
representation} $(\pi,\mathbf{H}(\pi))$ of the symmetry group, in our case
$E_{n}$ and $P_{n}$, the problem is reduced to use ``analytic continuation''
to move unitary representations of $E_{n}$ to unitary representations of
$P_{n}$ by passing over to imaginary time. This idea was used in the paper by
J. Fr\"{o}hlich, K. Osterwalder, and E. Seiler in \cite{FOS83}, see also
\cite{KlLa83}, to \textit{construct quantum field theoretical systems using
Euclidean field theory}. In this paper we will give some general constructions
and ideas related to this problem in the context of the applications
(\ref{IntCat(1)})--(\ref{IntCat(5)}) mentioned above, and work out some simple examples.

In \cite{Sch86} R. Schrader used this idea to construct, from a
\textit{complementary series representation} of $SL(2n,\mathbb{C})$, a unitary
representation of the group $SU(n,n)\times SU(n,n)$. In that paper the
similarities to the \textit{Yang-Baxter} relation were also discussed, a theme
that we will leave out in this exposition. What was missing in R.~Schrader's
paper was the \textit{identification} of the resulting representations and a
general procedure \textit{how} to construct those representations. We will see
that we can do this for all simple Lie groups where the associated Riemannian
symmetric space $G/K$ is a tube domain, and the that the duality works between
complementary series representations and highest weight representations.

In general this problem can be formulated in terms of $c$\textit{-duality} of
Lie groups and the analytic continuation of unitary representations from one
real form to another. (See \cite{HO95} for these terms.) The representations
that show up in the case of semisimple groups are generalized principal series
representations on the one side and highest weight representations on the
other. The symmetric spaces are the causal symmetric spaces, and the duality
is between non-compactly causal symmetric spaces and compactly causal
symmetric spaces. The latter correspond bijectively to real forms of bounded
symmetric domains. Therefore both the geometry and the representations are
closely related to the work of Harish-Chandra on bounded symmetric domains and
the holomorphic discrete series \cite{HCIV,HCV,HCVI}. But the ideas are also
related to the work of I. Segal and S. Paneitz through the notion of causality
and invariant cones, \cite{Segal,Paneitz81,Paneitz84}. A more complete
exposure can be found in the joint paper by the coauthors: \textit{Unitary
Representations of Lie Groups with Reflection Symmetry}, \cite{JO97}.

There are other interesting and related questions, problems, and directions.
In particular we would like to mention the analytic continuation of the
$H$-invariant character of the highest weight representations, and the
reproducing kernel of the Hardy space realization of this representation to a
spherical distribution character (spherical function) of the generalized
principal series representation that we started with. This connects the
representations that show up in the duality on the level of distribution
characters. For this we refer to \cite{O98, BK98}.

The paper is organized as follows. The list also includes some sources for
additional background references:\bigskip\newline
\ref{Considerations}.
Some Spectral Considerations Related to Reflection Positivity
\cite{FOS83,OsSc73,OsSc75}
\hrulefill\ \pageref{Considerations}\newline
\ref{uopg}.
Unitary One-Parameter Groups and Path Space Integrals
\cite{Arv84,KlLa81}
\hrulefill\ \pageref{uopg}\newline
\ref{axbg}.
The $(ax+b)$-Group
\cite{JO97,Ped90}
\hrulefill\ \pageref{axbg}\newline
\ref{Hilbert}.
The Hilbert Transform
\cite{GF89,Seg98}
\hrulefill\ \pageref{Hilbert}\newline
\ref{RemConsiderationsMar.1}.
One-Parameter Groups
\hrulefill\ \pageref{RemConsiderationsMar.1}\newline
\ref{Setting}.
The General Setting
\cite{JO97}
\hrulefill\ \pageref{Setting}\newline
\ref{Preliminaries}.
Preliminaries
\cite{GlJa87,Phil,Sto51,Mac}
\hrulefill\ \pageref{Preliminaries}\newline
\ref{Basic}.
A Basic Lemma
\cite{JO97}
\hrulefill\ \pageref{Basic}\newline
\ref{S-Hrep}.
Holomorphic Representations
\cite{HC65,HC70,Nee94,'OO88a,'OO88b}
\hrulefill\ \pageref{S-Hrep}\newline
\ref{LM}.
The L\"uscher-Mack Theorem
\cite{LM75}
\hrulefill\ \pageref{LM}\newline
\ref{SSS}.
Bounded Symmetric Domains
\cite{HiNe93,HO95,Hel62,He78}
\hrulefill\ \pageref{SSS}\newline
\ref{S-hwm}.
Highest Weight Modules
\cite{EHW83,Stanton,Ne94,Ne99,FK94,HO,Ja83,KNO98}
\hrulefill\ \pageref{S-hwm}\newline
\ref{S:Ex}.
An Example: $SU(1,1)$
\cite{JO97}
\hrulefill\ \pageref{S:Ex}\newline
\ref{S:Sssp}.
Reflection Symmetry for Semisimple Symmetric Spaces
\cite{Jor86,Jor87,JO97,Sch86}
\hrulefill\ \pageref{S:Sssp}\newline
\ref{S:Bargmann}.
The Segal-Bargmann Transform
\cite{OO96}
\hrulefill\ \pageref{S:Bargmann}\newline
\ref{Diagonal}.
The Heisenberg Group
\cite{BH94,Jor88,JO97}
\hrulefill\ \pageref{Diagonal}\newline
\ref{axb}.
The $(ax+b)$-Group Revisited
\cite{JO97,LaPh,Hel64}
\hrulefill\ \pageref{axb}\newline
References
\hrulefill\ \pageref{refpage}

\section{\label{Considerations}Some Spectral Considerations Related to
Reflection Positivity}

\subsection{\label{uopg}Unitary One-Parameter Groups and Path Space Integrals}

The term \emph{reflection positivity} is from quantum field theory (QFT) where
it refers to a certain reflection in the time-variable. As we explained in the
introduction this reflection is also the one which makes the analytic
continuation between the (Newtonian) group of rigid motions and the
Poincar\'{e} group of relativity. This is the approach to QFT of Osterwalder
and Schrader \cite{OsSc73,OsSc75}. The approach implies a change of the inner
product, and the new Hilbert space which carries a \emph{unitary}
representation $\tilde{\pi}$ of the Poincar\'{e} group $P_{4}$ results from
the corresponding ``old'' one by passing to a subspace where the positivity
(the so-called Osterwalder-Schrader positivity) is satisfied. An energy
operator may then be associated with this ``new'' unitary representation
$\tilde{\pi}$ of $P_{4}$. This representation $\tilde{\pi}$ is ``physical'' in
that the corresponding energy is positive. The basic connection between the
two groups may further be understood from the corresponding quadratic forms on
space-time $\left(  x,t\right)  $, $x=\left(  x_{1},x_{2},x_{3}\right)  $,
$\left(  x,t\right)  \mapsto\left\|  x\right\|  ^{2}+t^{2}$ with $\left\|
x\right\|  ^{2}=x_{1}^{2}+x_{2}^{2}+x_{3}^{2}$. The analytic continuation
$t\mapsto it$, $i=\sqrt{-1}$, turns this into the form $\left\|  x\right\|
^{2}-t^{2}$ of relativity. This same philosophy may also be used in an
analytic continuation argument connecting Feynman measure with the Wiener
measure on path space. This is important since the Wiener measure seems to be
the most efficient way of making precise the Feynman measure, which involves
infinite ``renormalizations'' if given a literal interpretation. We refer to
\cite{GlJa87} and \cite{Nel64} for more details on this point.

For the convenience of the reader we include here a simple instance of
reflection positivity for a path space measure which will be needed later: Let
$\mathcal{D}=\mathcal{D}\left(  \mathbb{R}\right)  =C_{c}^{\infty}\left(
\mathbb{R}\right)  $ denote the usual test functions on $\mathbb{R}$, and the
dual space $\mathcal{D}^{\prime}=\mathcal{D}^{\prime}\left(  \mathbb{R}%
\right)  $ of distributions. (We use the notation $q\left(  f\right)  $,
$q\in\mathcal{D}^{\prime}$, $f\in\mathcal{D}$.) Let $H_{0}=-\frac{1}%
{2}\bigtriangleup+\frac{1}{2}q^{2}-\frac{1}{2}$ be the harmonic oscillator
Hamiltonian, and form $\hat{H}=H_{0}-E_{0}$ picking $E_{0}$ such that $\hat
{H}\geq0$ and $\hat{H}\Omega=0$ for a ground state vector $\Omega$. Then by
\cite{GlJa87} there is a unique path-space measure $\phi_{0}$ on
$\mathcal{D}^{\prime}$ such that
\begin{align*}
\int q\left(  t\right)  \,d\phi_{0}\left(  q\right)   &  =0,\\
\int q\left(  t_{1}\right)  q\left(  t_{2}\right)  \,d\phi_{0}\left(
q\right)   &  =\frac{1}{2}e^{-\left|  t_{1}-t_{2}\right|  }.
\end{align*}
One rigorous interpretation (see \cite{Nel64,Nel73}) is to view $\left(
q\left(  t\right)  \right)  $ here as a stochastic process, i.e., a family of
random variables indexed by $t$. Further, for each $t>0$, and for each ``even
+ linear'' real potential $V\left(  q\right)  $, there is a unique measure
$\mu_{t}$ on $\mathcal{D}^{\prime}$ such that
\[
d\mu_{t}=Z_{t}^{-1}\exp\left(  -\int_{-t/2}^{t/2}V\left(  q\left(  s\right)
\right)  \,ds\right)  \,d\phi_{0}\left(  q\right)
\]
with
\[
Z_{t}=\int\exp\left(  -\int_{-t/2}^{t/2}V\left(  q\left(  s\right)  \right)
\,ds\right)  \,d\phi_{0}\left(  q\right)  .
\]
The consideration leading to the relation between the measures $d\phi_{0}$ and
$d\mu_{t}$ is the Trotter approximation for the (analytically continued)
semigroup,%
\[
e^{-it\left(  -\frac{1}{2}\bigtriangleup+V\right)  }=\lim_{n\rightarrow\infty
}\left(  e^{\left(  it/2n\right)  \bigtriangleup}e^{-\left(  it/n\right)
V}\right)  ^{n}.
\]
See \cite{Nel64} for more details on this point. When the operator $\left(
\cdots\right)  ^{n}$ on the right-hand side is computed, we find the integral
kernel
\[
K^{\left(  n\right)  }\left(  x_{0},x_{n},t\right)  =\frac{1}{N_{n}}\int
e^{i\mathcal{S}\left(  x_{0},\dots,x_{n},t\right)  }\,dx_{1}\,\cdots\,dx_{n-1}%
\]
with%
\[
N_{n}=\left(  \frac{2\pi it}{n}\right)  ^{3n/2}%
\]
and%
\[
\mathcal{S}\left(  x_{0},\dots,x_{n},t\right)  =\frac{1}{2}\sum_{j=1}%
^{n}\left|  x_{j-1}-x_{j}\right|  ^{2}\left(  \frac{t}{n}\right)  ^{-1}%
-\sum_{j=1}^{n}V\left(  x_{j}\right)  \left(  \frac{t}{n}\right)  .
\]
The heuristic motivation for $\phi_{0}$ and $\mu_{t}$ is then the ``action''
$\mathcal{S}$ approximating%
\[
\mathcal{S}\left(  q\right)  =\frac{1}{2}\int_{0}^{t}\left(  \dot{q}\left(
s\right)  \right)  ^{2}\,ds-\int_{0}^{t}V\left(  q\left(  s\right)  \right)
\,ds
\]
via%
\[
q\left(  t_{j}\right)  =x_{j},\qquad\bigtriangleup t_{j}=\frac{t}{n},
\]
and%
\[
\dot{q}\left(  t_{j}\right)  \sim\frac{x_{j}-x_{j-1}}{t/n}%
\]
If $-t/2<t_{1}\leq t_{2}\leq\dots\leq t_{n}<t/2$, then
\begin{multline}%
\ip{\Omega}{A_{1}e^{-\left( t_{2}-t_{1}\right) \hat
{H}}A_{2}e^{-\left( t_{3}-t_{2}\right) \hat{H}}A_{3}\cdots A_{n}\Omega}%
\label{eqConsiderationsNew.1}\\
=\lim_{t\rightarrow\infty}\int\prod_{k=1}^{n}A_{k}\left(  q\left(
t_{k}\right)  \right)  \,d\mu_{t}\left(  q\left(  \,\cdot\,\right)  \right)  .
\end{multline}
(The reader can find more details in \cite{GlJa87}.) Using further Minlos'
theorem
(see, e.g., \cite{GlJa87} or \cite{ReSi75}),
it can be shown that there is a measure $\mu$ on $\mathcal{D}%
^{\prime}$ such that
\begin{equation}
\lim_{t\rightarrow\infty}\int e^{iq\left(  f\right)  }\,d\mu_{t}\left(
q\right)  =\int e^{iq\left(  f\right)  }\,d\mu\left(  q\right)  =:S\left(
f\right)  \label{eqConsiderationsNew.2}%
\end{equation}
for all $f\in\mathcal{D}$. Since $\mu$ is a positive measure, we have
\begin{equation}
\sum_{k}\sum_{l}\bar{c}_{k}c_{l}S\left(  f_{k}-\bar{f}_{l}\right)  \geq0
\label{eqConsiderationsNew.3}%
\end{equation}
for all $c_{1},\dots,c_{n}\in\mathbb{C}$, and all $f_{1},\dots,f_{n}%
\in\mathcal{D}$. When combining (\ref{eqConsiderationsNew.1}) and
(\ref{eqConsiderationsNew.2}), we note that this limit-measure $\mu$ then
accounts for the time-ordered $n$-point functions which occur on the left-hand
side in formula (\ref{eqConsiderationsNew.1}). This observation will be
further used in the analysis of a corresponding stationary stochastic process
$\mathbf{X}_{t}$, $\mathbf{X}_{t}\left(  q\right)  =q\left(  t\right)  $,
which we proceed to analyze and review. But, more importantly, it can be
checked from the construction that we also have the following reflection
positivity: Let $\left(  \theta f\right)  \left(  s\right)  :=f\left(
-s\right)  $, $f\in\mathcal{D}$, $s\in\mathbb{R}$, and set
\[
\mathcal{D}_{+}=\left\{  f\in\mathcal{D}\mid f\text{ real valued, }f\left(
s\right)  =0\text{ for }s<0\right\}  \,.
\]
Then
\begin{equation}
\sum_{k}\sum_{l}\bar{c}_{k}c_{l}S\left(  \theta\left(  f_{k}\right)
-f_{l}\right)  \geq0 \label{eqConsiderationsNew.4}%
\end{equation}
for all $c_{1},\dots,c_{n}\in\mathbb{C}$, and all $f_{1},\dots,f_{n}%
\in\mathcal{D}_{+}$. (A small technical point: In the use of Minlos' theorem
it is often more convenient to use the duality $\mathcal{S}\leftrightarrow
\mathcal{S}^{\prime}$ of tempered distributions, as opposed to the
$\mathcal{D}\leftrightarrow\mathcal{D}^{\prime}$ duality mentioned above. The
reason for this is that the Hermite functions are in $\mathcal{S}$ but not in
$\mathcal{D}$.) These concepts are covered in detail in \cite{OsSc73,OsSc75}.

While our work here was centered on the role of reflection positivity for the
structure of representations of non-compact Lie groups $G$, the basic concepts
are important even in the case when $G=\mathbb{R}$ (the real line). In that
case, the data is as follows:

\begin{itemize}
\item $\mathbf{H}$: a complex Hilbert space (compare with
(\ref{eqConsiderationsNew.3}))

\item $\mathbf{K}_{0}$: a closed subspace (compare with $\mathcal{D}_{+}$)

\item $\left\{  U\left(  t\right)  \mid t\in\mathbb{R}\right\}  $: a unitary
one-parameter group acting on $\mathbf{H}$

\item $J\colon\mathbf{H}\rightarrow\mathbf{H}$: a period-2 unitary operator
(compare with $f\mapsto\theta\left(  f\right)  $) satisfying
\begin{align}
JU\left(  t\right)  &=U\left(  -t\right)  J, & &  \label{eqConsiderations.1}\\
P_{0}JP_{0}&\geq0  &  &\text{(compare with (\ref{eqConsiderationsNew.4}%
)),}\label{eqConsiderations.2}\\
P_{0}U\left(  t\right)  P_{0} &=U\left(  t\right)  P_{0}  &  &\text{for all
}t\geq0,\label{eqConsiderations.3}
\end{align}%
where $P_{0}$ is the orthogonal projection of $\mathbf{H}$ onto $\mathbf{K}%
_{0}$.

\item  Define $\mathbf{N}:=\left\{  k_{0}\in\mathbf{K}_{0}\mid%
\ip{k_{0}}{Jk_{0}}%
=0\right\}  $.
\end{itemize}

To illustrate that the axiom system (\ref{eqConsiderations.1}%
)--(\ref{eqConsiderations.3}) is very restrictive, we note that it is
\emph{not} satisfied for the usual translation group on $\mathbb{R}$.
Specifically, let $\mathbf{H}=L^{2}\left(  \mathbb{R}\right)  $, and $U\left(
t\right)  f\left(  x\right)  =f\left(  x-t\right)  $, $f\in L^{2}\left(
\mathbb{R}\right)  $, $x,t\in\mathbb{R}$. Since the subspaces $\mathbf{K}_{0}$
for which (\ref{eqConsiderations.3}) hold are known by Beurling's theorem, see
\cite{Hel64} and \cite{LaPh}, we see that we cannot have all three
(\ref{eqConsiderations.1})--(\ref{eqConsiderations.3}) in this example, unless
$\mathbf{K}_{0}=\mathbf{N}$. If (\ref{eqConsiderations.3}) holds, then either
(case 1) $\mathbf{K}_{0}$ is invariant under all $U\left(  t\right)  $,
$t\in\mathbb{R}$, and then it consists of all functions whose Fourier
transform is supported on some measurable subset of $\mathbb{R}$ (depending on
$\mathbf{K}_{0}$), or else (case 2) it consists of the transforms of functions
in $qH^{2}\left(  \mathbb{R}\right)  $ where $q$ is a unitary function and
$H^{2}\left(  \mathbb{R}\right)  $ is the usual Hardy space. Hence we may
assume that $J$ is given by%
\[
\left(  Jf\right)  \sphat\left(  \xi\right)  =\hat{f}\left(  -\xi\right)
,\qquad\xi\in\mathbb{R},
\]
where $\hat{f}$ is the Fourier transform. It is then easy to check in case 1
that we cannot have (\ref{eqConsiderations.2}), and that in case 2,
$\mathbf{N}=\mathbf{K}_{0}$.

The simplest instance of (\ref{eqConsiderations.1})--(\ref{eqConsiderations.3}%
) arises in quantum field theory and for certain stochastic processes; see,
e.g., \cite[Lecture 4]{Arv84}, \cite{Kle78}. Let $\left\{  \mathbf{X}_{t}\mid
t\in\mathbb{R}\right\}  $ be a stationary stochastic process on a probability
space $\left(  \Omega,P\right)  $ which is \emph{symmetric} in the sense that
$\mathbf{X}_{-t}$ has the same distribution as does $\mathbf{X}_{t}$. In that
case $\mathbf{H}=L^{2}\left(  \Omega,P\right)  $ and $\mathbf{K}_{0}$ may be
taken to be the subspace in $L^{2}\left(  \Omega,P\right)  $ generated by the
functions which are measurable with respect to the $\sigma$-field generated by
$\left\{  \mathbf{X}_{t}\mid t\geq0\right\}  $, and $P_{0}$ may be taken to be
the corresponding conditional expectation. Since the process is stationary, it
generates a unitary one-parameter group in $L^{2}\left(  \Omega,P\right)  $,
and the choice of $P_{0}$ makes it clear that (\ref{eqConsiderations.1}) and
(\ref{eqConsiderations.3}) will be satisfied. Condition
(\ref{eqConsiderations.2}) is an extra condition, which is called
Osterwalder-Schrader positivity (or O-S positivity, for short), although the
O-S positivity concept was first formulated in a different context; see, e.g.,
\cite[Axioms 1--2]{OsSc73,OsSc75}.

Both in the case of unitary one-parameter groups and in the general context of
representations of Lie groups, there is an operator-theoretic step used in
passing from $\mathbf{H}$ to the new Hilbert space. It underlies two technical
points involved in the construction: \emph{positivity} and
\emph{norm-estimates}. It is given in Section \ref{Basic} and referred to as
the \emph{Basic Lemma}. More details are in \cite{JO97}.

An application of our Basic Lemma, Section \ref{Basic}, produces a contraction
$W$,
\begin{equation}%
\begin{array}
[c]{ccc}%
& \makebox[0pt]{\hss$\displaystyle\mathbf{K}_{0}/\mathbf{N}$\hss} & \\
\llap{\raisebox{2pt}{\small quotient}}\rlap{$\displaystyle\nearrow$} &  &
\llap{$\displaystyle\searrow$}\rlap{\raisebox{2pt}{\small completion}}\\
\mathbf{K}_{0} & \underset{W}{\longrightarrow} & \rlap{$\displaystyle\left(
\mathbf{K}_{0}/\mathbf{N}\right) \sptilde$}\hphantom{( \mathbf{K}_{0}}%
\end{array}
\label{eqConsiderations.4}%
\end{equation}
such that $\left\{  U\left(  t\right)  \mid t\geq0\right\}  $ is realized in
$\mathbf{\hat{H}}_{+}:=\left(  \mathbf{K}_{0}/\mathbf{N}\right)  \sptilde$ as
a contractive and selfadjoint semigroup $\vphantom{\hat{U}}\left\{
\smash{\hat{U}}\left(  t\right)  \mid t\geq0\right\}  $, and
\begin{equation}
WU\left(  t\right)  |_{\mathbf{K}_{0}}=\hat{U}\left(  t\right)  W,\qquad
t\geq0. \label{eqConsiderations.5}%
\end{equation}
Since this induced semigroup is contractive and selfadjoint, it can be shown
to have the form
\begin{equation}
\hat{U}\left(  t\right)  =e^{-tH},\qquad t\geq0 \label{eqConsiderations.6}%
\end{equation}
for a selfadjoint operator $H$, $H\geq0$, $H$ acting in $\mathbf{\hat{H}}_{+}%
$. The semigroup $\hat{U}\left(  t\right)  $ in (\ref{eqConsiderations.6}) is
constructed here by a procedure which is also applicable to a large class of
non-compact Lie groups. In the present case, the idea is simple: {}From the
axioms (\ref{eqConsiderations.1})--(\ref{eqConsiderations.3}), we get
\begin{equation}%
\ip{k_{0}}{JU\left( t\right) k_{0}}%
\leq%
\ip{k_{0}}{Jk_{0}}%
,\text{\qquad for }k_{0}\in\mathbf{K}_{0}\text{ and for all }t\geq0.
\label{eqConsiderations.7}%
\end{equation}
Hence $\mathbf{N}:=\left\{  k_{0}\in\mathbf{K}_{0}\mid%
\ip{k_{0}}{Jk_{0}}%
=0\right\}  $ is invariant under $\left\{  U\left(  t\right)  \mid
t\geq0\right\}  $ which then passes to the quotient $\mathbf{K}_{0}%
/\mathbf{N}$ as a contraction semigroup on the completed Hilbert space. The
selfadjointness of the induced semigroup follows from the identity
\begin{equation}%
\ip{k_{1}}{JU\left( t\right) k_{2}}%
=%
\ip{U\left( t\right) k_{1}}{Jk_{2}}%
, \label{eqConsiderations.8}%
\end{equation}
valid for all $k_{1},k_{2}\in\mathbf{K}_{0}$, and $t\geq0$. The proof of
(\ref{eqConsiderations.8}) in turn is immediate from axiom
(\ref{eqConsiderations.1}).

When applied to the stochastic process example, we get a concrete realization
of this semigroup $\hat{U}\left(  t\right)  =e^{-tH}$ on $\mathbf{\hat{H}}%
_{+}$. In pure operator-theoretic terms, what results is the following data:

\begin{enumerate}
\item \label{uopgitem(1)}$\mathbf{K}_{0}$ and $\mathbf{\hat{H}}_{+}$, two
Hilbert spaces;

\item \label{uopgitem(2)}$\left\{  V\left(  t\right)  \right\}  _{t\geq0}$, a
semigroup of isometries in $\mathbf{K}_{0}$;

\item \label{uopgitem(3)}$H\geq0$, a selfadjoint, generally unbounded
operator, in $\mathbf{\hat{H}}_{+}$;

\item \label{uopgitem(4)}$W\colon\mathbf{K}_{0}\rightarrow\mathbf{\hat{H}}%
_{+}$, a contractive linear operator which has the following properties:
$\ker\left(  W\right)  =0$; and $\ker\left(  W^{\ast}\right)  =0$, i.e., the
range of $W$ is dense, which is to say, $W\left(  \mathbf{K}_{0}\right)  $ is
dense in $\mathbf{\hat{H}}_{+}$ relative to the norm of $\mathbf{\hat{H}}_{+}$;

\item \label{uopgitem(5)}$W$ intertwines the semigroups in the respective
spaces, i.e.,%
\[
e^{-tH}W=WV\left(  t\right)  \text{,\qquad for all }t\geq0.
\]
\end{enumerate}

Note that the properties listed in (\ref{uopgitem(4)}) for $W$ imply that the
polar decomposition $W=W_{0}\left(  W^{\ast}W\right)  ^{1/2}$ has the partial
isometric factor $W_{0}\colon\mathbf{K}_{0}\rightarrow\mathbf{\hat{H}}_{+}$ a
unitary isomorphism, but, of course, $W_{0}$ will \emph{not} intertwine
$e^{-tH}$ and $\left\{  V\left(  t\right)  \right\}  _{t\geq0}$, i.e., the
relation in (\ref{uopgitem(5)}) does not pass to the polar decomposition.

It is this realization which we call the O-S construction. Since the constant
function $\openone\in L^{2}\left(  \Omega,P\right)  $ is cyclic for the
$L^{\infty}\left(  \Omega\right)  $-multiplication algebra acting on
$L^{2}\left(  \Omega,P\right)  $, we get $\Omega:=W\left(  \openone\right)
\in\mathbf{\hat{H}}_{+}$ satisfying a similar cyclicity property in
$\mathbf{\hat{H}}_{+}$ and also $\hat{U}\left(  t\right)  \Omega=\Omega$ for
all $t\geq0$, or equivalently, $\Omega$ is in the domain of $H$, and
$H\Omega=0$.

If $L^{2}\left(  \Omega,P\right)  $, $\left\{  U\left(  t\right)  \right\}  $
are constructed from a stationary stochastic process $\mathbf{X}_{t}$ as
described, then we may get (\ref{eqConsiderations.2}) satisfied if
\begin{multline}
\int_{\Omega}f_{1}\circ\mathbf{X}_{t_{1}}\,f_{2}\circ\mathbf{X}_{t_{2}%
}\,\cdots\,f_{n}\circ\mathbf{X}_{t_{n}}\,dP\geq0\text{\qquad for all }%
f_{1},\dots,f_{n}\in C_{c}\left(  \mathbb{R}\right)
,\label{eqConsiderations.9}\\
\text{and all }t_{1},\dots,t_{n}\in\mathbb{R}\text{ such that }-\infty
<t_{1}\leq t_{2}\leq\dots\leq t_{n}<\infty.
\end{multline}
Moreover, condition (\ref{eqConsiderations.9}) is alternately denoted the
Osterwalder-Schrader positivity condition. It can be obtained if instead we
start with a semigroup $\vphantom{\hat{U}\hat{\mathbf{H}}}\left(  \smash
{\hat{U}}\left(  t\right)  ,\smash
{\hat{\mathbf{H}}}\right)  $ and a representation $f\mapsto\pi\left(
f\right)  $ of $C_{c}\left(  \mathbb{R}\right)  $ on $\mathbf{\hat{H}}$ such
that, for a vector $v\in\mathbf{\hat{H}}$, we have%
\begin{multline}%
\ip{v}{\pi\left( f_{1}\right) \smash{\hat{U}}\left( t_{2}-t_{1}\right
) \pi\left( f_{2}\right) \smash{\hat{U}}\left
( t_{3}-t_{2}\right) \cdots\smash{\hat{U}}\left( t_{n}-t_{n-1}\right) \pi
\left( f_{n}\right) v}%
\label{eqConsiderations.10}\\
\text{for all }f_{1},\dots,f_{n}\in C_{c}\left(  \mathbb{R}\right)  \text{,
and }-\infty<t_{1}\leq t_{2}\leq\dots\leq t_{n}<\infty.
\end{multline}
Starting with one of the two, (\ref{eqConsiderations.9}) or
(\ref{eqConsiderations.10}), the other can be constructed such that the
expressions in (\ref{eqConsiderations.9}) and (\ref{eqConsiderations.10}) are
equal. This is the content of the O-S construction in the formulation of
E.~Nelson and others \cite{Nel73,KlLa75}.

To understand (\ref{eqConsiderations.5}) better, it is useful to compare it to
the Nagy dilation of a semigroup; see \cite{SzN74}. The Nagy theorem states
that every contraction semigroup $\vphantom{\hat{U}\hat{\mathbf{H}}}\left(
\smash{\hat{U}}\left(  t\right)  ,\smash
{\hat{\mathbf{H}}},t\geq0\right)  $ admits a representation $\left(  U\left(
t\right)  ,\mathbf{H},t\in\mathbb{R}\right)  $ where $\mathbf{\hat{H}}%
\subset\mathbf{H}$, and $U\left(  t\right)  $ is a one-parameter unitary group
in $\mathbf{H}$ which satisfies%
\begin{equation}
P_{\mathbf{\hat{H}}}U\left(  t\right)  |_{\mathbf{\hat{H}}}=\hat{U}\left(
t\right)  ,\qquad t\geq0. \label{eqConsiderations.11}%
\end{equation}
While there are some generalizations of (\ref{eqConsiderations.11}) to
(non-abelian) Lie groups, e.g., \cite{JoMu80}, the relation
(\ref{eqConsiderations.5}) is the focus of the present paper, wherein we show
that it characterizes a class of ``physical'' representations of semisimple
non-compact Lie groups. In (\ref{eqConsiderations.11}), $P_{\mathbf{\hat{H}}}$
denotes the orthogonal projection of $\mathbf{H}$ onto $\mathbf{\hat{H}}$.
Hence (\ref{eqConsiderations.11}) may be rephrased in terms of an isometry
$V\colon\mathbf{\hat{H}}\rightarrow\mathbf{H}$, and we get%
\begin{equation}
V^{\ast}U\left(  t\right)  V=\hat{U}\left(  t\right)  ,\qquad t\geq0.
\label{eqConsiderations.12}%
\end{equation}

In case the stochastic process $\left(  \mathbf{X}_{t},t\in\mathbb{R}\right)
$ is given on path space $\mathbf{X}_{t}\left(  \omega\right)  =\omega\left(
t\right)  $, Arveson found (in \cite[Proposition 4.6]{Arv84}) a reformulation
of (\ref{eqConsiderations.5}) much in the spirit of (\ref{eqConsiderations.12}%
), but Arveson produces a simultaneous ``dilation'' of a given semigroup
\emph{and} a representation of $C_{c}\left(  \mathbb{R}\right)  $. First
consider the following two representations of $C_{c}\left(  \mathbb{R}\right)
$:

\begin{itemize}
\item $\sigma$ representing $C_{c}\left(  \mathbb{R}\right)  $ on
$L^{2}\left(  \Omega,P\right)  $ given by
\begin{equation}
\sigma\left(  f\right)  F\left(  \omega\right)  =f\left(  \omega\left(
0\right)  \right)  F\left(  \omega\right)  , \label{eqConsiderations.13}%
\end{equation}
for $f\in C_{c}\left(  \mathbb{R}\right)  $, $F\in L^{2}\left(  \Omega
,P\right)  $, and

\item $M$ representing $C_{c}\left(  \mathbb{R}\right)  $ on $\mathbf{\hat{H}%
}$ by%
\begin{equation}
\left(  M_{f}h\right)  \left(  t\right)  =f\left(  t\right)  h\left(
t\right)  , \label{eqConsiderations.14}%
\end{equation}
for $f\in C_{c}\left(  \mathbb{R}\right)  $, $h\in\mathbf{\hat{H}}$.
\end{itemize}

Let $v$ be the vector in (\ref{eqConsiderations.10}). Arveson then shows that
the linear mapping%
\begin{equation}
M_{f}v\longmapsto\sigma\left(  f\right)  \openone\label{eqConsiderations.15}%
\end{equation}
extends uniquely to an isometry%
\[
V\colon L^{2}\left(  \mathbb{R}\right)  \longrightarrow L^{2}\left(
\Omega,P\right)
\]
which satisfies%
\begin{multline}
V^{\ast}U\left(  t_{1}\right)  \sigma\left(  f_{1}\right)  U\left(
t_{2}\right)  \sigma\left(  f_{2}\right)  \cdots U\left(  t_{n}\right)
\sigma\left(  f_{n}\right)  V\label{eqConsiderations.16}\\
=\hat{U}\left(  t_{1}\right)  M_{f_{1}}\hat{U}\left(  t_{2}\right)  M_{f_{2}%
}\cdots\hat{U}\left(  t_{n}\right)  M_{f_{n}}%
\end{multline}
for all $n=1,2,\dots$, $t_{i}\in\mathbb{R}$, $t_{i}\geq0$, and all $f_{i}\in
C_{c}\left(  \mathbb{R}\right)  $.

The space $L^{2}\left(  \mathbb{R}\right)  $ is thereby identified as a closed
subspace in $L^{2}\left(  \Omega,P\right)  $. It is the time-zero subspace,
and is ``much smaller'' than the subspace generated by the $\sigma$-algebra of
$\left(  \mathbf{X}_{t},t\geq0\right)  $. There is a selfadjoint semigroup
induced from both of the subspaces. In general, it is the $t=0$ semigroup
which has the Markoff property in the sense of Nelson \cite[Theorem 1]{Nel73}.

\subsection{\label{axbg}The $(ax+b)$-Group}

The general case breaks down into the analysis of the solvable case, and the
semisimple case.
Therefore
it is appropriate to start with the $2$-dimensional solvable
case, that is the $(ax+b)$-group. We may realize this group $G$ as the affine
transformations $ax+b$ of $\mathbb{R}$, $a\in\mathbb{R}_{+}$, $b\in\mathbb{R}%
$. Hence up to a trivial scale, the Lie algebra is determined by the relation
$\left[  A,B\right]  =B$. Let $\tau(a,b)=(a,-b)$. A unitary representation of
$G$ on a Hilbert space $\mathbf{H}$ is therefore specified by two unitary
one-parameter groups $U_{A}\left(  s\right)  $, $U_{B}\left(  t\right)  $,
$s,t\in\mathbb{R}$, satisfying%
\begin{equation}
U_{A}\left(  s\right)  U_{B}\left(  t\right)  U_{A}\left(  -s\right)
=U_{B}\left(  e^{s}\cdot t\right)  . \label{eqConsiderations.17}%
\end{equation}
The spectra of the two groups form subsets of $\mathbb{R}$, and
(\ref{eqConsiderations.17}) shows that the spectrum $\Lambda\left(  B\right)
=\operatorname*{spec}\left(  U_{B}\right)  $ is invariant under
positive dilations,
that is
\begin{equation}
\mathbb{R}_{+}\cdot\Lambda\left(  B\right)  =\Lambda\left(  B\right)  .
\label{eqConsiderations.18}%
\end{equation}
But (\ref{eqConsiderations.18}) implies that either $\Lambda\left(  B\right)
=\mathbb{R}_{+}$, $\Lambda\left(  B\right)  =\mathbb{R}_{-}$, or
$\Lambda\left(  B\right)  =\mathbb{R}$.

Based on these considerations, we have the following:

\begin{theorem}
\label{ThmNoRefl}There are no nontrivial reflection symmetries for infinite
dimensional unitary representations of the $(ax+b)$ group.
\end{theorem}

\begin{proof}
[Proof \textup{(}sketch\/\textup{)}]It is enough to exclude the reflection
which sends $A\mapsto A$, and $B\mapsto-B$. We must show that in every
$\mathbf{K}_{0}\subset\mathbf{H}$ which is invariant under $\left\{
U_{A}\left(  s\right)  \mid s\in\mathbb{R}\right\}  $ and under $\left\{
U_{B}\left(  t\right)  \mid t\geq0\right\}  $, and for every $\left(
J,\mathbf{K}_{0}\right)  $ such that (\ref{eqConsiderations.2}) holds, and
$JU_{A}\left(  s\right)  =U_{A}\left(  s\right)  J$ and $JU_{B}\left(
t\right)  =U_{B}\left(  -t\right)  J$, we must have $\mathbf{\hat{H}}$
($=\left(  \mathbf{K}_{0}/\mathbf{N}\right)  \sptilde
$) the trivial zero-dimensional space.

Recalling the new norm in $\mathbf{\hat{H}}$, $f\mapsto%
\ip{f}{Jf}%
=\left\|  f\right\|  _{\mathbf{\hat{H}}}^{2}$, and
\[
\mathbf{N}=\left\{  f\in\mathbf{K}_{0}\mid%
\ip{f}{Jf}%
=0\right\}  ,
\]
we conclude that the induced transformations $\hat{U}_{A}$ and $\hat{U}_{B}$
on $\mathbf{\hat{H}}$ satisfy:

\begin{enumerate}
\item \label{ThmNoRefl.proof(1)}$\hat{U}_{A}$ is a unitary one-parameter group
on $\mathbf{\hat{H}}$ and its spectrum is a subset of $\Lambda\left(
A\right)  $;

\item \label{ThmNoRefl.proof(2)}$\vphantom{\hat{U}}\left\{  \smash{\hat{U}%
}_{B}\left(  t\right)  \mid t\geq0\right\}  $ is a contraction semigroup on
$\mathbf{\hat{H}}$ satisfying $\hat{U}_{B}\left(  t\right)  ^{\ast}=\hat
{U}_{B}\left(  t\right)  $; and

\item \label{ThmNoRefl.proof(3)}$\hat{U}_{A}\left(  s\right)  \hat{U}%
_{B}\left(  t\right)  \hat{U}_{A}\left(  -s\right)  =\hat{U}_{B}\left(
e^{s}\cdot t\right)  $ for all $s\in\mathbb{R}$, and $t\geq0$.
\end{enumerate}

\noindent Combining (\ref{ThmNoRefl.proof(2)})--(\ref{ThmNoRefl.proof(3)}),
and using a theorem of \cite{Ped90}, we may assume that the selfadjoint
generator $H_{B}$ of $U_{B}\left(  t\right)  $ ($=e^{itH_{B}}$) and $J$ have
the representation%
\begin{equation}
H_{B}=%
\begin{pmatrix}
H & 0\\
0 & -H
\end{pmatrix}
,\qquad J=%
\begin{pmatrix}
0 & I\\
I & 0
\end{pmatrix}
\label{eqConsiderations.pound}%
\end{equation}
relative to $\mathbf{H}=%
\begin{pmatrix}
\mathbf{H}_{+}\\
\mathbf{H}_{-}%
\end{pmatrix}
$ for closed subspaces $\mathbf{H}_{\pm}$. In this representation it is
possible to check all the candidates for $\mathbf{K}_{0}\subset\mathbf{H}$
with the stated properties; see also Sections \ref{Preliminaries} and
\ref{axb} below. The presence of a nontrivial $\mathbf{K}_{0}$ leads to
(\ref{eqConsiderations.pound}) and the conclusion that $\Lambda\left(
B\right)  =\mathbb{R}$. But, of the nontrivial candidates for $\mathbf{K}_{0}%
$, none can satisfy the added axioms%
\[
P_{0}JP_{0}\geq0
\]
and%
\[
P_{0}U_{B}\left(  t\right)  P_{0}=U_{B}\left(  t\right)  P_{0},\text{\qquad
for all }t\geq0.
\]
We postpone further details to a future paper on semidirect products.
\end{proof}

We refer to Section \ref{axb} for more detailed discussion on the $(ax+b)$-group.

\subsection{\label{Hilbert}The Hilbert Transform}

In a recent paper \cite{Seg98}, Segal obtains a positive energy representation
of the Poincar\'{e} group $P$ on a Hilbert space of ``complex'' spinors. The
construction is a renormalization of the usual Klein-Gordon inner product, by
use of an operator $J$ much like the one described above. In Fourier
variables, it is $J\colon\phi\left(  k\right)  \mapsto i\theta\left(
k\right)  \phi\left(  k\right)  $, where $k=\left(  k_{0},k_{1},k_{2}%
,k_{3}\right)  $ denotes a point in momentum space, and%
\[
\theta\left(  k\right)  =%
\begin{cases}
+1 & \text{if }k_{0}\geq0,  \\
-1 & \text{if }k_{0}<0.
\end{cases}
\]
Hence $J$ appears as the Hilbert transform with respect to time. One of the
corollaries of our results in Section \ref{SSS} is that, even in the general
case of reflection positivity for semisimple Lie groups, the appropriate $J$
must indeed be a generalized Hilbert transform, and the construction is tied
to the complementary series.

\subsection{\label{RemConsiderationsMar.1}One-Parameter Groups}

In addition to the examples from path-space
measures, there are those which arise directly from the theory of
representations of reductive Lie groups. While this applies completely
generally, the
$SL(2,\mathbb{R})\simeq SU(1,1)$
case is worked out in detail in Section
\textup{\ref{S:Ex}} below. When the representations there are restricted to a
suitable one-parameter subgroup in
$SL(2,\mathbb{R})$,
we arrive at the following basic
setup for the axiom system \textup{(\ref{eqConsiderations.1}%
)--(\ref{eqConsiderations.3}).}

Let $0<s<1$ be given. Let $\mathbf{H}:=\mathbf{H}(s)$ be the Hilbert space
given by the norm squared%
\[
\left\|  f\right\|  ^{2}:=\int_{\mathbb{R}}\int_{\mathbb{R}}\overline{f\left(
x\right)  }f\left(  y\right)  \left|  x-y\right|  ^{s-1}\,dx\,dy<\infty\,.
\]
With $s$ fixed, let $\mathbf{K}_{0}\subset\mathbf{H}$ be the subspace of
functions in $\mathbf{H}$ which
have compact support
in $\left(  -1,1\right)  $. Let
$\left\{  U\left(  t\right)  \right\}  _{t\in\mathbb{R}}$ be the unitary
one-parameter group given by%
\begin{equation}
U_{s}\left(  t\right)  f\left(  x\right)  =e^{\left(  s+1\right)  t}f\left(
e^{2t}x\right)  . \label{eqConsiderationsMar.17}%
\end{equation}
It follows from representation theory \textup{(}see Section \textup{\ref{S:Ex}%
)} that, for every $s\in\left(  0,1\right)  $, the spectrum of the group
$U_{s}\left(  \,\cdot\,\right)  $ in \textup{(\ref{eqConsiderationsMar.17})}
is continuous and is all of $\mathbb{R}$.

Let $J\colon\mathbf{H}\rightarrow\mathbf{H}$ be given by%
\[
Jf\left(  x\right)  =\left|  x\right|  ^{-s-1}f\left(  1/x\right)  \,.
\]
Then an elementary calculation shows that the axioms
\textup{(\ref{eqConsiderations.1})--(\ref{eqConsiderations.3}) are }satisfied
when $P_{0}$ denotes the projection onto $\mathbf{K}_{0}$. We will then form
the completion of $\mathbf{K}_{0}$ relative to the norm $\left\|
\,\cdot\,\right\|  _{\widehat{\mathbf{H}_{+}(s)}}$ given by%
\begin{equation}
\left\|  f\right\|  _{\widehat{\mathbf{H}_{+}(s)}}^{2}=%
\ip{f}{Jf}%
_{\mathbf{H}(s)}\,, \label{eqConsiderationsMar.18}%
\end{equation}
for $f\in\mathbf{K}_{0}$. While we are completing a space of functions, it
turns out that the elements in the completion are generally distributions
which may not be functions.

We stress this example since it is the simplest instance when the completion
\textup{(\ref{eqConsiderations.4})} is made explicit as a concrete space of
distributions. It follows from the representation-theoretic setup
\textup{(}Section \textup{\ref{S:Ex})} that%
\[
\mathbf{\hat{H}}_{+}
=\left(
\mathbf{K}_{0}/\mathbf{N}\right)  \sptilde
\]
is then the Hilbert space of distributions obtained from completion of
measurable functions on $\left(  -1,1\right)  $ completed relative to the new
norm squared from \textup{(\ref{eqConsiderationsMar.18}),} viz.,%
\begin{equation}
\int_{-1}^{1}\int_{-1}^{1}\overline{f\left(  x\right)  }f\left(  y\right)
\left|  1-xy\right|  ^{s-1}\,dx\,dy<\infty\,. \label{eqConsiderationsMar.19}%
\end{equation}

The fact that $\widehat{\mathbf{H}_{+}(s)}$, for $0<s<1$, is a space of
distributions is significant for the spectral theory. The Dirac delta
``function'' $\delta$ is in $\widehat{\mathbf{H}_{+}(s)}$ and is a ground
state vector. When the unitary one-parameter group
\textup{(\ref{eqConsiderationsMar.17})} is passed to $\widehat{\mathbf{H}%
_{+}(s)}$, we get \textup{(}for each $s$\textup{)} the positive selfadjoint
generator $H=H_{s}$ from \textup{(\ref{eqConsiderations.6}),} and we check
that $H\delta=\left(  1-s\right)  \delta$ with $1-s=$the bottom of the
spectrum of $H$. It is true in fact that the full spectrum of $H$ in
$\widehat{\mathbf{H}_{+}(s)}$ is $\left\{  2n+1-s\mid n=0,1,2,\dots\right\}
$, and that the spectrum is simple.

It can be shown that $\widehat{\mathbf{H}_{+}(s)}$, for $0<s<1$, from
\textup{(\ref{eqConsiderationsMar.18})--(\ref{eqConsiderationsMar.19})} in
fact are associated with measures on $\mathbb{R}$ as follows. The assertion is
that, for every $f\in\widehat{\mathbf{H}_{+}(s)}$,%
\[
C_{c}\left(  -1,1\right)  \ni\varphi\longmapsto%
\ip{f}{\varphi}%
_{\widehat{\mathbf{H}_{+}(s)}}%
\in \mathbb{C}
\]
extends to a Radon measure,%
\[
\mu_{f}\left(  \varphi\right)  =\int_{\mathbb{R}}\varphi\left(  x\right)
\,d\mu_{f}\left(  x\right)  \,,
\]
on $\mathbb{R}$. To see this, use the estimates%
\begin{align}
\left|
\ip{f}{\varphi}%
_{\widehat{\mathbf{H}_{+}(s)}}\right|  &\leq\left\|  f\right\|  _{\widehat
{\mathbf{H}_{+}(s)}}\left\|  \varphi\right\|  _{\widehat{\mathbf{H}_{+}(s)}%
}  , \label{eqConsiderationsApr.24} \\
\intertext{and}\left\|  \varphi\right\|  _{\widehat{\mathbf{H}_{+}(s)}%
}  &\leq\operatorname*{const}\times\sup_{x\in\left(  -1,1\right)  }\left|
\varphi\left(  x\right)  \right|  \,,\qquad\forall\,\varphi\in C_{c}\left(
-1,1\right)  \,. \label{eqConsiderationsApr.25} 
\end{align}

We will study the Hilbert spaces
$\mathbf{H}(s)$ and $\widehat{\mathbf{H}_{+}(s)}$ further in Section
\textup{\ref{S:Ex}} below, where we show among other things that $\left\{
\left(  d/dx\right)  ^{n}\delta\mid n=0,1,2,\dots\right\}  $ forms an
orthogonal basis in the reflection Hilbert space $\widehat{\mathbf{H}_{+}(s)}$
for all $s\in\left(  0,1\right)  $. This is a way to turn the Taylor expansion
into an orthogonal decomposition.

\section{\label{Setting}The General Setting}

The setup is a connected Lie group $G$ with an nontrivial involution, that is
a symmetry, $\tau\colon G\rightarrow G$. The differential $\tau\colon
\frak{{g}\rightarrow\frak{g}}$ is then an involution on the Lie algebra
$\frak{g}$ of $G$. Let $H:=G^{\tau}=\left\{  x\in G\mid\tau(x)=x\right\}  $
and $\frak{h}:=\{X\in\frak{g}\mid\tau(X)=X\}$. Then $\frak{h}$ is the Lie
algebra of $H$. Let $\frak{q}:=\{X\in\frak{g}\mid\tau(X)=-X\}$. Then
$\frak{q}$ is a vector space isomorphic to the tangent space $T_{x_{o}}(G/H)$,
where $x_{o}=eH\in G/H$.

The $c$-dual of $(\frak{g},\tau)$ is defined to be the Lie algebra
\[
\frak{g}^{c}:=\frak{h}+i\frak{q}\subset\frak{g}_{\mathbb{C}}%
\]
with the involution $\tau^{c}:=\tau|\frak{g}^{c}$. We let $G^{c}$ be the
simply connected Lie group with the Lie algebra $\frak{{g}^{c}}$. Then $\tau$
integrates to an involution $\tau^{c}\colon G^{c}\rightarrow G^{c}$, and
$H^{c}:=(G^{c})^{\tau^{c}}$ is connected, \cite{'O87}. Here are few examples
of triples $(G,G^{c},\tau)$:

\textbf{1. Compact Lie groups:} Let $G$ be a compact Lie group. Then $G^{c}$
is a simply connected reductive Lie group with Cartan involution $\tau^{c}$,
and every connected simply connected reductive Lie group can be constructed in
this way.

\textbf{2. The group case:} Let $H$ be a Lie group and let $G=H\times H$. We
identify $H$ with the diagonal $d(H)=\{(a,a)\in G\mid a\in H\}$. Define
$\tau(a,b):=(b,a)$. Then $G^{\tau}=H$. Furthermore the map
\[
G/H\ni(a,b)H\longmapsto ab^{-1}\in H
\]
is a diffeomorphism, intertwining the canonical action of $G$ on $G/H$ and the
action $(a,b)\cdot x=axb^{-1}$ on $H$. Denote the simply connected complex Lie
group with the Lie algebra $\frak{h}_{\mathbb{C}}$ by $H_{\mathbb{C}}$. Then
$G^{c}=H_{\mathbb{C}}$, and $\tau^{c}$ is the conjugation with respect to the
real form $\frak{h}$. Hence the analysis on Lie groups and their
complexification is a special case of symmetric space duality. The case
$H_{\mathbb{C}}=SL(2n,\mathbb{C})$ and $H=SU(n,n)$ was treated by R. Schrader
in \cite{Sch86}.

\textbf{3. Semidirect product with abelian normal subgroup:} Let $G=QH$ with
$Q$ and $H$ connected, $Q$ normal subgroup of $G$ and $Q\cap H=\{e\}$. Then we
can define an involution on $G$ by
\[
\tau(qh):=q^{-1}h\,.
\]
Thus $\frak{q}$ is the Lie algebra of $Q$. Define $\Phi\colon\frak{g}%
\rightarrow\frak{g}^{c}$ by
\[
\Phi(X+Y):=iX+Y\,,\quad X\in\frak{q}\,,\;Y\in\frak{h}\,.
\]
Then $\Phi$ is a Lie algebra isomorphism. Thus $G^{c}$ is the simply connected
covering group of $G$. A special case hereof is the $(ax+b)$-group and the
Heisenberg group, where $H$ is also abelian.

\textbf{4. The $(ax+b)$-group:} The $(ax+b)$-group is the group of
transformations $x\mapsto ax+b,\quad a>0,\,b\in\mathbb{R} $. We take $Q=
\mathbb{R}$ and $H=\mathbb{R}^{+}$. Then $\tau$ is given by $\tau
(a,b)\;=(a,-b)$.

\textbf{5. The Heisenberg group:} Let $H_{n}$ be the $(2n+1)$-dimensional
Heisenberg group. We write
\[
H_{n}=\left\{  h(\mathbf{x},\mathbf{y},z)=\left(
\begin{array}
[c]{lll}%
1 & \mathbf{x}^{t} & z\\
0 & I_{n-1} & \mathbf{y}\\
0 & 0 & 1
\end{array}
\right)  \,\,\mathbf{x},\mathbf{y}\in\mathbb{R}^{n},\,z\in\mathbb{R}\right\}
\simeq\mathbb{R}^{n}\times\mathbb{R}^{n}\times\mathbb{R}%
\]
with multiplication given by
\[
h(\mathbf{x},\mathbf{y},z)h(\mathbf{x}^{\prime},\mathbf{y}^{\prime},z^{\prime
})=h(\mathbf{x}+\mathbf{x}^{\prime},\mathbf{y}+\mathbf{y}^{\prime}%
,z+z^{\prime}+(\mathbf{x}\mid\mathbf{y}^{\prime}))\,.
\]
In particular $h(\mathbf{x},\mathbf{y},z)^{-1}=h(-\mathbf{x},-\mathbf{y}%
,(\mathbf{x}\mid\mathbf{y})-z)$. In this case we take $H=\{h(\mathbf{x}%
,0,0)\}\simeq\mathbb{R}^{n}$, which is abelian, and $Q=\{h(0,\mathbf{y},z)\in
H_{n}\mid\mathbf{y}\in\mathbb{R}^{n},\,z\in\mathbb{R}\}\simeq\mathbb{R}^{n+1}%
$. The involution is given by
\[
\tau(h(\mathbf{x},\mathbf{y},z)):=h(\mathbf{x},-\mathbf{y},-z)\,.
\]

Starting from the pair $(G,\tau)$ and a unitary representation $(\pi
,\mathbf{H}(\pi))$ of $G$ we need a compactible involution on the Hilbert
space $\mathbf{H}(\pi)$, that is a unitary linear map $J\colon\mathbf{H}%
(\pi)\rightarrow\mathbf{H}(\pi)$ intertwining $\pi$ and $\pi\circ\tau$. Thus
\[
J\pi(g)=\pi(\tau(g))J\qquad\forall g\in G\,.
\]
We will also need a semigroup $S$ such that $H\subset S$ or at least
$H\subset\overline{S}$ where the bar denotes topological closure. We shall
consider closed subspaces $\mathbf{K}_{0}\subset\mathbf{H}(\pi)$, where
$\mathbf{H}(\pi)$ is the Hilbert space of $\pi$, such that $\mathbf{K}_{0}$ is
invariant under $\pi(S^{o})$. Let $J\colon\mathbf{H}(\pi)\rightarrow
\mathbf{H}(\pi)$ be a unitary intertwining operator for $\pi$ and $\pi
\circ\tau$ such that $J^{2}=\operatorname*{id}$. We assume that $\mathbf{K}%
_{0}$ may be chosen such that $\Vert v\Vert_{J}^{2}:=%
\ip{v}{Jv}%
\geq0$ for all $v\in\mathbf{K}_{0}$. We will always assume our inner product
conjugate linear in the first argument. We form, in the usual way, the Hilbert
space $\mathbf{K}=\left(  \mathbf{K}_{0}/\mathbf{N}\right)  \sptilde$ by
dividing out with $\mathbf{N}=\{v\in\mathbf{K}_{0}\mid%
\ip{v}{Jv}%
=0\}$ and completing in the norm $\Vert\cdot\Vert_{J}$. (This is of course a
variation of the Gelfand-Naimark-Segal (GNS) construction.) With the
properties of $(G,\pi,\mathbf{H}(\pi),\mathbf{K})$ as stated, we show, using
the L\"uscher-Mack theorem, that the simply connected Lie group $G^{c}$ with
Lie algebra $\frak{g}^{c}=\frak{h}\oplus i\frak{q}$ carries a \textit{unitary}
representation $\pi^{c}$ on $\mathbf{K}$ such that $\{\pi^{c}(h\exp(iY))\mid
h\in H,Y\in C^{o}\}$ is obtained from $\pi$ by passing the corresponding
operators $\pi(h\exp Y)$ to the quotient $\mathbf{K}_{0}/\mathbf{N}$. To see
that this leads to a unitary representation $\pi_{c}$ of $G^{c}$ we use a
basic result of L\"uscher and Mack \cite{LM75} and in a more general context
one of Jorgensen \cite{Jor86}. In fact, when $Y\in C$, the selfadjoint
operator $d\pi(Y)$ on $\mathbf{K}$ has spectrum contained in $(-\infty,0]$. As
in
Lemma \ref{MNoughtZero},
we show that in the case where $C$ extends to
an $G^{c}$ invariant regular cone in $i\frak{g}^{c}=i\frak{h}\oplus\frak{q}$
and $\pi^{c}$ is injective, then each $\pi^{c}$ (as a unitary representation
of $G^{c}$) must be a direct integral of highest-weight representations of
$G^{c}$. The examples show that one can relax the condition in different ways,
that is one can avoid using the L\"uscher-Mack theorem by instead
constructing local representations and using only cones that are neither
generating nor $H$-invariant.

Assume now that $G$ is a semidirect product of $H$ and $N$ with $N$ normal and
abelian. Define $\tau\colon G\rightarrow G$ by $\tau(hn)=hn^{-1}$. Let $\pi
\in\hat{H}$ (the unitary dual) and extend $\pi$ to a unitary representation of
$G$ by setting $\pi(hn)=\pi(h)$. In this case, $G^{c}$ is locally isomorphic
to $G$, and $\pi$ gives rise to a unitary representation $\pi^{c}$ of $G^{c}$
by the formula $d\pi^{c}(X)=d\pi(X)$, $X\in\frak{h}$, and $d\pi^{c}%
|_{i\frak{q}}=0$. A special case of this is the $3$-dimensional Heisenberg
group, and the $(ax+b)$-group. In Sections \ref{Diagonal} and \ref{axb}, we
show that, if we induce instead a character of the subgroup $N$ to $G$, then
we have $\left(  \mathbf{K}_{0}/\mathbf{N}\right)  \sptilde=\{0\}$.

Our approach to the general representation correspondence $\pi\mapsto\pi^{c} $
is related to the integrability problem for representations of Lie groups (see
\cite{JoMo84}); but the present positivity viewpoint comes from
Osterwalder-Schrader positivity; see \cite{OsSc73, OsSc75}. In addition the
following other papers are relevant in this connection:
\cite{FOS83,Jor86,Jor87,KlLa83,Pra89,Sch86}.

\section{\label{Preliminaries}Preliminaries}

\setcounter{equation}{0}

The setting for the paper is a general Lie group $G$ with a nontrivial
involutive automorphism $\tau$.

\begin{definition}
\label{ReflectionSymmetric}A unitary representation $\pi$ acting on a Hilbert
space $\mathbf{H}(\pi)$ is said to be \textit{reflection symmetric} if there
is a unitary operator $J\colon\mathbf{H}(\pi)\rightarrow\mathbf{H}(\pi)$ such
that \begin{list}{}{\setlength{\leftmargin}{\customleftmargin}
\setlength{\itemsep}{0.5\customskipamount}
\setlength{\parsep}{0.5\customskipamount}
\setlength{\topsep}{\customskipamount}
\setlength{\parindent}{0pt}}
\item[\hss\llap{\rm R1)}] ${\displaystyle J^2 = \operatorname*{id}}$;
\item[\hss\llap{\rm R2)}] ${\displaystyle J\pi(g) = \pi(\tau(g))J\, ,
\quad g\in G}$.
\end{list}
\end{definition}

If (R1) holds, then $\pi$ and $\pi\circ\tau$ are equivalent. Furthermore,
generally from (R2) we have $J^{2}\pi(g) = \pi(g) J^{2}$. Thus, if $\pi$ is
irreducible, then we can always renormalize $J$ such that (R1) holds. Let $H =
G^{\tau}= \{g\in G\mid\tau(g) = g\}$ and let $\frak{h}$ be the Lie algebra of
$H$. Then $\frak{h} = \{ X\in\frak{g}\mid\tau(X) = X\}$. Define $\frak{q} =
\{Y\in\frak{g}\mid\tau(Y) = -Y\}$. Then $\frak{g} = \frak{h} \oplus\frak{q}$,
$[\frak{h},\frak{q}]\subset\frak{q}$ and $[\frak{q},\frak{q}]\subset\frak{h}$.

\begin{definition}
\label{Hyperbolic}A closed convex cone $C\subset\frak{q}$ is
\textit{hyperbolic} if $C^{o}\not =\emptyset$ and if $\operatorname{ad}X$ is
semisimple with real eigenvalues for every $X\in C^{o}$.
\end{definition}

We will assume the following for $(G, \pi,\tau, J)$: \begin{list}{}%
{\setlength{\leftmargin}{\customleftmargin}
\setlength{\itemsep}{0.5\customskipamount}
\setlength{\parsep}{0.5\customskipamount}
\setlength{\topsep}{\customskipamount}}
\item[\hss\llap{{\rm PR1)}}] $\pi
$ is reflection symmetric with reflection $J$;
\item[\hss\llap{{\rm PR2)}}] there is an $H$-invariant hyperbolic cone
$C\subset\frak{q}$ such that $S(C) = H\exp C$ is a closed semigroup and
$S(C)^o = H\exp C^o$ is diffeomorphic to $H\times C^o$;
\item[\hss\llap{{\rm PR3)}}] there is a subspace ${0}\not=
\mathbf{K}_0\subset\mathbf{H}(\pi
)$ invariant under $S(C)$ satisfying the positivity
condition
\[ \ip{v}{v}_J:= \ip{v}{J(v)} \ge0,\quad\forall v\in\mathbf{K}_0\, .\]
\end{list}

\begin{remark}
\label{WhatToAssume}In (PR3) we can always assume that $\mathbf{K}_{0}$ is
closed, as the invariance and the positivity pass over to the closure. In
(PR2) it is only necessary to assume that $\mathbf{K}_{0}$ is invariant under
$\exp C$, as one can always replace $\mathbf{K}_{0}$ by $\overline
{\left\langle \pi(H)\mathbf{K}_{0}\right\rangle }$, the closed space generated
by $\pi(H)\mathbf{K}_{0}$, which is $S(C)$-invariant, as $C$ is $H$-invariant.
For the exact conditions on the cone for (PR2) to hold see the original paper
by J. Lawson \cite{JL94}, or the monograph \cite[pp. 194 ff.]{HiNe93}.
\end{remark}

In some of the examples we will replace (PR2)and (PR3) by the following:
weaker conditions \begin{list}{}{\setlength{\leftmargin}{\customleftmargin}
\setlength{\itemsep}{0.5\customskipamount}
\setlength{\parsep}{0.5\customskipamount}
\setlength{\topsep}{\customskipamount}}
\item[\hss\llap{\rm PR2$^{\prime}$)}] $C$ is (merely) some nontrivial
cone in $\frak{q} $.
\item[\hss\llap{\rm PR3$^{\prime}$)}] There is a subspace
$0\not= \mathbf{K}_0\subset\mathbf{H} (\pi)$ invariant under $H$ and $\exp C$
satisfying the positivity condition from (PR3).
\end{list}
(See Section \ref{Diagonal} for further details.)

Since the operators $\{\pi(h)\mid h\in H\}$ commute with $J$, they clearly
pass to the quotient by
\[
\mathbf{N}:=\{v\in\mathbf{K}_{0}\mid%
\ip{v}{ Jv}%
=0\}
\]
and implement unitary operators on $\mathbf{K}:=\left(  \mathbf{K}%
_{0}/\mathbf{N}\right)  \sptilde$ relative to the inner product induced by
\begin{equation}%
\ip{u}{v}%
_{J}:=%
\ip{u}{J(v)}%
\,. \label{E:innerpr}%
\end{equation}
which will be denoted by the same symbol. Hence we shall be concerned with
passing the operators $\{\pi(\exp Y)\mid Y\in C\}$ to the quotient
$\mathbf{K}_{0}/\mathbf{N}$, and for this we need a basic Lemma.

In general, when $\left(  \mathbf{K}_{0},J\right)  $ is given, satisfying the
positivity axiom, then the corresponding composite quotient mapping
\[
\mathbf{K}_{0}\longrightarrow\mathbf{K}_{0}/\mathbf{N}\hooklongrightarrow
\left(  \mathbf{K}_{0}/\mathbf{N}\right)  \sptilde=:\mathbf{K}%
\]
is \emph{contractive} relative to the respective Hilbert norms. The resulting
(contractive) mapping will be denoted $\beta$. An operator $\gamma$ on
$\mathbf{H}$ which leaves $\mathbf{K}_{0}$ invariant is said to \emph{induce}
the operator $\tilde{\gamma}$ on $\mathbf{K}$ if $\beta\circ\gamma
=\tilde{\gamma}\circ\beta$ holds on $\mathbf{K}_{0}$. In general, an induced
operation $\gamma\mapsto\tilde{\gamma}$ may not exist; and, if it does,
$\tilde{\gamma}$ may fail to be bounded, even if $\gamma$ is bounded.

This above-mentioned operator-theoretic formulation of reflection positivity
has applications to the Feynman-Kac formula in mathematical physics, and there
is a considerable literature on that subject, with work by E. Nelson
\cite{Nel64,Nel73}, A. Klein and L.J. Landau \cite{Kle78,KlLa75,KlLa81}, B.
Simon, and W.B. Arveson \cite{Arv84}. Since we shall not use path space
measures here, we will omit those applications, and instead refer the reader
to the survey paper \cite{Arv84} (lecture 4) by W.B. Arveson. In addition to
mathematical physics, our motivation also derives from recent papers on
non-commutative harmonic analysis which explore analytic continuation of the
underlying representations; see, e.g., \cite{HOO91,BK98,Nee94,'O90a,O93,Ol82}.

\section{\label{Basic}A Basic Lemma}

\setcounter{equation}{0}

\begin{lemma}
\label{BasicLemma}

\begin{enumerate}
\item [\hss\llap{\rm1)}]Let $J$ be a period-$2$ unitary operator on a Hilbert
space $\mathbf{H}$, and let $\mathbf{K}_{0}\subset\mathbf{H}$ be a closed
subspace such that $%
\ip{v}{J(v)}%
\geq0$, $v\in\mathbf{K}_{0}$. Let $\gamma$ be an invertible operator on
$\mathbf{H}$ such that $J\gamma=\gamma^{-1}J$ and which leaves $\mathbf{K}%
_{0}$ invariant and has $(\gamma^{-1})^{\ast}\gamma$ bounded on $\mathbf{H}$.
Then $\gamma$ induces a bounded operator $\tilde{\gamma}$ on $\mathbf{K}%
=\left(  \mathbf{K}_{0}/\mathbf{N}\right)  \sptilde$, where $\mathbf{N}%
=\{v\in\mathbf{K}_{0}\mid%
\ip{v}{Jv}%
=0\}$, and the norm of $\tilde{\gamma}$ relative to the $J$-inner product in
$\mathbf{K}$ satisfies
\begin{equation}
\Vert\tilde{\gamma}\Vert\leq\Vert(\gamma^{-1})^{\ast}\gamma\Vert_{sp}^{1/2}\,,
\label{E:3.1}%
\end{equation}
where $\Vert\cdot\Vert_{sp}$ is the spectral radius.

\item[\hss\llap{\rm2)}] If we have a semigroup $S$ of operators on
$\mathbf{H}$ satisfying the conditions in \textup{(1)}, then
\begin{equation}
(\gamma_{1}\gamma_{2})\sptilde=\tilde{\gamma_{1}}\tilde{\gamma_{2}}%
\,,\quad\gamma_{1},\gamma_{2}\in S\,. \label{E:3.2}%
\end{equation}
\end{enumerate}
\end{lemma}

\begin{proof}
For $v\in\mathbf{K}_{0}$, $v\not =0$, we have
\begin{align*}
\Vert\gamma(v)\Vert_{J}^{2}  &  =%
\ip{\gamma(v)}{J\gamma(v)}%
\\
&  =%
\ip{\gamma(v)}{\gamma^{-1}J(v)}%
\\
&  =%
\ip{(\gamma^{-1})^{*}\gamma(v)}{J(v)}%
\\
&  =%
\ip{(\gamma^{-1})^{*}\gamma(v)}{v}%
_{J}\\
&  \leq\Vert(\gamma^{-1})^{\ast}\gamma(v)\Vert_{J}\Vert v\Vert_{J}\\
&  \leq\Vert((\gamma^{-1})^{\ast}\gamma)^{2}(v)\Vert_{J}^{1/2}\Vert v\Vert
_{J}^{1+1/2}\\
&  \mkern4.5mu\vdots\\
&  \leq\Vert((\gamma^{-1})^{\ast}\gamma)^{2^{n}}(v)\Vert_{J}^{1/2^{n}}\Vert
v\Vert_{J}^{1+1/2+\cdots+1/2^{n}}\\
&  \leq\left(  \Vert((\gamma^{-1})^{\ast}\gamma)^{2^{n}}\Vert\Vert
v\Vert\right)  ^{1/2^{n}}\Vert v\Vert_{J}^{2}\,.
\end{align*}
Since $\displaystyle{\lim_{n\rightarrow\infty}\Vert((\gamma^{-1})^{\ast}%
\gamma)^{2^{n}}\Vert^{1/2^{n}}=\Vert(\gamma^{-1})^{\ast}\gamma\Vert_{sp}}$,
and ${\displaystyle
\lim_{n\rightarrow\infty}\Vert v\Vert^{1/2^{n}}=1}$, the result follows.

By this we get
\[%
\ip{\gamma(v)}{J\gamma(v)}%
\leq\Vert(\gamma^{-1})^{\ast}\gamma\Vert_{sp}%
\ip
{v}{J(v)}%
\]
which shows that $\gamma(\mathbf{N})\subset\mathbf{N}$, whence $\gamma$ passes
to a bounded operator on the quotient $\mathbf{K}_{0}/\mathbf{N}$ and then
also on $\mathbf{K}$ satisfying the estimate stated in (1). If both the
operators in (\ref{E:3.2}) leave $\mathbf{N}$ invariant, so does $\gamma
_{1}\gamma_{2}$ and the operator induced by $\gamma_{1}\gamma_{2}$ is
$\tilde{\gamma_{1}}\tilde{\gamma_{2}}$ as stated.
\end{proof}

\begin{corollary}
\label{GammaContraction}Let the notation be as above and assume that $\gamma$
is unitary on $\mathbf{H}$. Then the constant on the right in
\textup{(\ref{E:3.1})} is one. Hence $\tilde{\gamma}$ is a contraction on
$\mathbf{K}$.
\end{corollary}

To understand the assumptions on the space $\mathbf{K}_{0}$, that is
positivity and invariance, we include the following which is based on an idea
of R.S. Phillips \cite{Phil}.

\begin{proposition}
\label{P:3.3}Let $\mathbf{H}$ be a Hilbert space and let $J$ be a period-$2$
unitary operator on $\mathbf{H}$. Let $S$ be a commutative semigroup of
unitary operators on $\mathbf{H}$ such that $S=S_{+}S_{-}$ with $S_{+}%
=\{\gamma\in S\mid J\gamma=\gamma J\}$ and $S_{-}=\{\gamma\in S\mid
J\gamma=\gamma^{-1}J\}$. Then $\mathbf{H}$ possesses a maximal positive and
invariant subspace, that is a subspace $\mathbf{K}_{0}$ such that $%
\ip{v}{J(v)}%
\geq0$, $v\in\mathbf{K}_{0}$, and $\gamma\mathbf{K}_{0}\subset\mathbf{K}_{0}$,
$\gamma\in S$.
\end{proposition}

\begin{proof}
The basic idea is contained in \cite[pp. 386 ff.]{Phil}.
\end{proof}

\begin{remark}
\label{RemBasicNew.1}
A nice application is to the case $\mathbf{H}=L^{2}\left(  X,m\right)  $
where $X$ is a Stone space. There is a $m$-a.e.-defined automorphism
$\theta\colon X\rightarrow X$ such that
\[
J\left( f\right) =f\circ \theta ,\qquad f\in L^{2}\left( X,m\right)
\]
and $S$ is represented by multiplication operators on
$L^{2}\left(
X,m\right)  $.
By \cite{Phil} we know that there are clopen subsets $A,B\subset X$
such that with $M_{0}:=\left\{ x\in X\mid \theta \left( x\right) =x\right\} $
and $M_{1}=X\setminus M_{0}$
we have $A,B\subset M_{1}$, $A\cap B=\varnothing $,
$A\cup B=M_{1}$ and $\theta \left( A\right) =B$.
Let $\mathbf{K}_{0}:=L^{2}\left( M_{0}\cup A\right) $.
Then $\mathbf{K}_{0}$ is a maximal
positive and invariant subspace.
\end{remark}

\begin{lemma}
\label{MNoughtZero}If $M_{0}\subset X$ is of measure zero, then the space
$\mathbf{K}$ is trivial, that is $%
\ip{f}{J(f)}%
=0$ for all $f\in\mathbf{K}_{0}$.
\end{lemma}

\begin{remark}
\label{AbelianSubspace}Assume that we have (PR1) and (PR2).
By \cite{KrNe96} there is an
abelian subspace $\frak{a}\subset\frak{q}$ such that
$C^{o}=\operatorname{Ad}(H)(C^{o}\cap\frak{a})$. Let $S_{A}=\exp(C^{o}%
\cap\frak{a})$. Then $S_{A}$ is an abelian semigroup, so one can use
Proposition \ref{P:3.3} to construct a maximal positive and invariant subspace
for $S_{A}$. But in general we can not expect this space to be invariant under
$S$.
\end{remark}

We read off from the basic Lemma the following Proposition:

\begin{proposition}
\label{Contractive}Let $\pi$ be a unitary representation of the group $G$.
Assume that $(\tau,J,C,\mathbf{K}_{0})$ satisfies the conditions
\textup{(PR1), (PR2}$^{\prime}$\textup{),} and \textup{(PR3}$^{\prime}%
$\textup{).} If $Y\in C$, then $\pi(\exp Y)$ induces a contractive selfadjoint
operator $\tilde{\pi}(\exp Y)$ on $\mathbf{K}$.
\end{proposition}

\begin{proof}
If $Y\in C$, then $\pi(\exp Y)\mathbf{K}_{0}\subset\mathbf{K}_{0}$, and
$\pi(\exp Y)$ is unitary on $\mathbf{H}(\pi)$. Thus
\begin{align*}%
\ip{\pi(\exp Y)u}{J(v)}%
&  =%
\ip{u}{\pi(\exp(-Y))J(v)}%
\\
&  =%
\ip{u}{J(\pi(\exp Y)v)}%
\,,
\end{align*}
proving that $\pi(\exp Y)$ is selfadjoint in the $J$-inner product. Since
$\pi(\exp Y)$ is unitary on $\mathbf{H}(\pi)$
\[
\Vert\pi(\exp Y)\Vert=\Vert\pi(\exp Y)\Vert_{sp}=1\,,
\]
and the contractivity property follows.
\end{proof}

\begin{lemma}
\label{HolomorphicRPlus}Let $\pi$ be a unitary representation of $G$ such that
$(\tau,J,C,\mathbf{K}_{0})$ satisfies the conditions \textup{(PR1),
(PR2}$^{\prime}$\textup{),} and \textup{(PR3}$^{\prime}$\textup{).} Then for
$Y\in C$ there is a selfadjoint operator $d\tilde{\pi}(Y)$ in $\mathbf{K}%
=\left(  \mathbf{K}_{0}/\mathbf{N}\right)  \sptilde$, with spectrum contained
in $(-\infty,0]$, such that
\[
\tilde{\pi}(\exp(tY))=e^{t\,d\tilde{\pi}(Y)},\quad t\in\mathbb{R}_{+}%
\]
is a contractive semigroup on $\mathbf{K}$. Furthermore the following hold:

\begin{enumerate}
\item [\hss\llap{\rm1)}]$t\mapsto e^{t\,d\tilde{\pi}(Y)}$ extends to a
continuous map $z\mapsto e^{z\,d\tilde{\pi}(Y)}$ on $\{z\in\mathbb{C}%
\mid\operatorname{Re}(z)\geq0\}$ holomorphic on the open right half-plane, and
such that $e^{(z+w)\,d\tilde{\pi}(Y)}=e^{z\,d\tilde{\pi}(Y)}e^{w\,d\tilde{\pi
}(Y)}$.

\item[\hss\llap{\rm2)}] There exists a one-parameter group of unitary
operators
\[
\tilde{\pi}\left(  \exp(itY)\right)  :=e^{it\,d\tilde{\pi}(Y)},\quad
t\in\mathbb{R}%
\]
on $\mathbf{K}$.
\end{enumerate}
\end{lemma}

\begin{proof}
The last statement follows by the spectral theorem. By construction
$\{\tilde{\pi}(\exp(tY))\mid t\in\mathbb{R}_{+}\}$ is a semigroup of
selfadjoint contractive operators on $\mathbf{K}$. The existence of the
operators $d\tilde{\pi}(Y)$ as stated then follows from a general result in
operator theory; see, e.g., \cite{Fr80} or \cite{KlLa81}.
\end{proof}

\begin{corollary}
\label{Id}Let the situation be as in the last corollary. If $Y\in C\cap-C$
then $e^{t\,d\tilde{\pi}(Y)}=\operatorname*{id}$ for all $t\in\mathbb{R}_{+}$.
In particular $d\tilde{\pi}(Y)=0$ for every $Y\in C\cap-C$.
\end{corollary}

\begin{proof}
This follows as the spectrum of $d\tilde{\pi}(Y)$ and $d\tilde{\pi}(-Y)$ is
contained in $(-\infty,0]$.
\end{proof}

We remark here that we have introduced the map $d\pi$ without using the space
of \textit{smooth vectors} for the representation $\pi$. Let us recall that a
vector $\mathbf{v}\in\mathbf{H}(\pi)$ is called \textit{smooth} if the map
\[
\mathbb{R}\ni t\longmapsto\hat{\mathbf{v}}(t):=\pi(\exp tX)\mathbf{v}%
\in\mathbf{H}(\pi)
\]
is smooth for all $X\in\frak{g}$. The vector is \textit{analytic} if the above
map is analytic. We denote by $\mathbf{H}^{\infty}(\pi)$ the space of smooth
vectors and by $\mathbf{H}^{\omega}(\pi)$ the space of analytic vectors.
Both $\mathbf{H}^{\omega}(\pi)$ and 
$\mathbf{H}^{\infty}(\pi)$
are
$G$-invariant dense subspaces of $\mathbf{H}%
(\pi)$. We define a representation of $\frak{g}$ on $\mathbf{H}^{\infty}(\pi)$
by
\[
d\pi(X)\mathbf{v}=\lim_{t\rightarrow0}\frac{\pi(\exp tX)\mathbf{v}-\mathbf{v}%
}{t}\,.
\]

Recall that if $\pi$ is infinite-dimensional, then $d\pi$ is a representation
of $\frak{g}$ by unbounded operator on $\mathbf{H} (\pi)$, but the analytic
vectors and the $C^{\infty}$-vectors form dense domains for $d\pi$; see
\cite{Nel59,Pou92,WaI72}.

The operator $d\tilde{\pi}(X)$ in the above statements is an extension of the
operator $d\tilde{\pi}(X)$ on the space of smooth vectors. This allows us to
use the same notation for those two objects. We extend this representation to
$\frak{g}_{{\mathbb{C}}}$ by complex linearity, $d\pi(X+iY)=d\pi
(X)+i\,d\pi(Y)$, $X,Y\in\frak{g}$. Let $U(\frak{g})$ denote the universal
enveloping algebra of $\frak{g}_{{\mathbb{C}}}$. Then $d\pi$ extends to a
representation on $U(\frak{g})$, again denoted by $d\pi$. The space
$\mathbf{H}^{\infty}(\pi)$ is a topological vector space in a natural way,
cf.\ \cite{WaI72}. Furthermore $\mathbf{H}^{\infty}(\pi)$ is invariant under
$G$ and $U(\frak{g})$. As $\pi(g)\pi(\exp(tX))\mathbf{v}=\pi(\exp
(t\operatorname{Ad}(g)X))\pi(g)\mathbf{v}$, we get
\[
\pi(g)\,d\pi(X)\mathbf{v}=d\pi(\operatorname{Ad}(g)X)\pi(g)\mathbf{v},
\qquad \mathbf{v}\in\mathbf{H}^{\infty}(\pi)
\]
for all $g\in G$ and all $X\in\frak{g}$. Define $Z^{\ast}=-\sigma(Z)$,
$Z\in\frak{g}_{{\mathbb{C}}}$,
where $\sigma$ is the conjugation $X+iY\mapsto X-iY$, $X,Y\in\frak{g}$.
Then a simple calculation shows that for the
densely defined operator $\pi(Z)$, $Z\in\frak{g}_{{\mathbb{C}}}$, we have
$\pi(Z)^{\ast}=\pi(Z^{\ast})$
on $\mathbf{H}^{\infty}(\pi)$.

When (R1--2) and (PR1--3) hold, and $Y\in C$, we showed in
Lemma
\ref{HolomorphicRPlus} that the operator $\tilde{\pi}\left(  Y\right)  $ is
selfadjoint in $\mathbf{K}=\left(  \mathbf{K}_{0}/\mathbf{N}\right)  \sptilde$
with spectrum in $\left[  0,\infty\right)  $. Once $\tilde{\pi}$ is identified
as a unitary representation of $G^{c}$, then $\tilde{\pi}\left(  iY\right)  $
is automatically a selfadjoint operator in $\mathbf{K}$ by \cite{NeSt59}, but
semiboundedness of the corresponding spectrum of $\tilde{\pi}\left(  Y\right)
$ only holds for $Y\in C$. Yet if $\tilde{\pi}$ is obtained, as in
Lemma
\ref{HolomorphicRPlus}, from a unitary representation $\pi$ of $G$ acting on
$\mathbf{H}$, then the spectrum of $\pi\left(  Y\right)  $ is contained in the
purely imaginary axis $i\mathbb{R}$, and yet $\tilde{\pi}\left(  Y\right)  $
has spectrum in $\left[  0,\infty\right)  \subset\mathbb{R}$. The explanation
is that the Hilbert spaces $\mathbf{H}$ and $\mathbf{K}$ are different for the
two representations $\pi$ and $\tilde{\pi}$.

\section{\label{S-Hrep}Holomorphic Representations}

The unitary representations that show up in the duality are direct integrals
of highest weight $G^{c}$-modules. Those representations can also be viewed as
\textit{holomorphic representations} of a semigroup related to an extension of
the $H$-invariant cone $C$ that we started with. We will therefore give a
short overview over this theory, while
referring to the forthcoming monograph \cite{Ne99} for more details.

The theory of highest weight modules and holomorphic representations will
always be related to the name of Harish-Chandra because of his fundamental
work on bounded symmetric domains and the holomorphic discrete series,
\cite{HCIV,HCV,HCVI}. Later Gelfand and Gindikin in \cite{GG}
proposed a new approach for studying the Plancherel formula for
semisimple Lie group $G$. Their idea was to consider functions in
$\mathbf{L}^{2}(G)$ as the sum of boundary values of holomorphic functions
defined on domains in $G_{{\mathbb{C}}}$. The first deep results in this
direction are due to Ol'shanskii \cite{Ol82} and Stanton \cite{Stanton}, who
realized the holomorphic discrete series of the group $G$ in a Hardy space of
a local tube domain containing $G$ in the boundary. The generalization to
noncompactly causal symmetric spaces was carried out in
\cite{HO95,HOO91,'OO88a,'OO88b}. This program was carried out for solvable
groups in \cite{HO} and for general groups in \cite{BK95,Nee94}.

Let $G_{{\mathbb{C}}}$ be a complex Lie group with Lie algebra $\frak{g}%
_{{\mathbb{C}}}$ and let $\frak{g}$ be a real form of $\frak{g}_{{\mathbb{C}}%
}$. We assume for simplicity that $G$, the analytic subgroup of
$G_{{\mathbb{C}}}$ with Lie algebra $\frak{g}$, is closed in $G_{{\mathbb{C}}%
}$. Let $C$ be a regular $G$-invariant cone in $\frak{g}$ such that the set
$S(C)=G\exp iC$ is a closed semigroup in $G_{{\mathbb{C}}}$. Moreover, we
assume that the map
\[
G\times C\ni(a,X)\longmapsto a\exp iX\in S(C)
\]
is a homeomorphism and even a diffeomorphism when restricted to $G\times
C^{o}$. Finally, we assume that there exists a real automorphism $\sigma$ of
$G_{{\mathbb{C}}}$ whose differential is the complex conjugation of
$\frak{g}_{{\mathbb{C}}}$ with respect to $\frak{g}$, that is $\sigma
(X+iY)=X-iY$ for $X,Y\in\frak{g}$. We notice that in this case $G_{{\mathbb{C}%
}}/G$ is a symmetric space and that the corresponding subspace $\frak{q}$ is
just $i\frak{g}$. Those hypotheses are also satisfied for Hermitian Lie groups
and also for some solvable Lie groups; cf.\ \cite{HO}. Define%
\begin{equation}
W(\pi):=\{X\in\frak{g}\mid\forall u\in\mathbf{H}^{\infty}(\pi)\,:\,\ip
{i\,d\pi(X)u}{u}\leq0\}\,. \label{Cpi}%
\end{equation}
Thus $W(\pi)$ is the set of elements of $\frak{g}$ for which $d\pi(iX)$ is
\textit{negative}. The elements of $W(\pi)$ are called \textit{negative
elements} for the representation $\pi$.

\begin{lemma}
\label{LemS-HrepNew.1}$W(\pi)$ is a closed $G$-invariant convex cone in
$\frak{g}$.
\end{lemma}

\begin{definition}
\label{DefS-HrepMar.1}Let $W$ be a $G$-invariant cone in $\frak{g}$. We denote
the set of all unitary representations $\pi$ of $G$ with $W\subset W(\pi)$ by
$\mathcal{A}(W)$. A unitary representation $\pi$ is called \textit{$W$%
}\emph{-admissible} if $\pi\in\mathcal{A}(W)$.
\end{definition}

The representations in $\mathcal{A}(W)$ will be studied in detail in Section
\ref{S-hwm} below. We show in Theorem \ref{AChw} that a $\rho\in
\mathcal{A}(W)$ which is irreducible is in fact a \emph{highest weight
representation}, and the corresponding $K^{c}$-weights are determined. The
representations are then identified as discrete summands in $\mathbf{L}%
^{2}\left(  G^{c}\right)  $.

Let $S$ be a semigroup with unit and let $\sharp\colon S\rightarrow S$ be a
bijective involutive antihomomorphism, that is%
\[
(ab)^{\sharp}=b^{\sharp}a^{\sharp}\qquad\text{and}\qquad\quad a^{\sharp\sharp
}=a
\]

We call $\sharp$ an \textit{involution} on the semigroup $S$, and we call the
pair $(S,\sharp)$ a \textit{semigroup with involution} or an
\textit{involutive semigroup}. For us the important examples are the
semigroups of the form $S(C)=H\exp C$ with $\gamma^{\sharp}=\tau(\gamma^{-1}%
)$. Another class of examples consists of the \textit{contractive semigroups}
on a Hilbert space $\mathbf{H}$. Let $S(\mathbf{H})=\{T\in B(\mathbf{H}%
)\mid\Vert T\Vert\leq1\}$. Denote by $T^{\ast}$ the adjoint of $T$ with
respect to the inner product on $\mathbf{H}$. Then $(S,\ast)$ is a semigroup
with involution.

\begin{definition}
Let $(S,\sharp)$ be a topological semigroup with involution: then a semigroup
homomorphism $\rho\colon S\rightarrow S(\mathbf{H})$ is called a contractive
representation of $(S,\sharp)$ if $\rho(g^{\sharp})=\rho(g)^{\ast}$ and $\rho$
is continuous with respect to the weak operator topology of $S(\mathbf{H})$. A
contractive representation is called irreducible\label{irr} if there is no
closed nontrivial subspace of $\mathbf{H}$ invariant under $\rho(S)$.
\end{definition}

\begin{definition}
\label{DefS-HrepNew.3}Let $\rho$ be a contractive representation of the
semigroup $S(W)=G\exp iW\subset G_{{\mathbb{C}}}$. Then $\rho$ is holomorphic
if the function $\rho\colon S(C)^{o}\rightarrow B(\mathbf{V})$ is holomorphic.
\end{definition}

The following lemma shows that, if a unitary representation of the group $G$
extends to a holomorphic representation of $S (C)$, then this extension is unique.

\begin{lemma}
\label{unique} If $f\colon S(W)\rightarrow S(\mathbf{H})$ is continuous and
$f|_{S(W)^{o}}$ is holomorphic such that $f|_{G}=0$, then $f=0$.
\end{lemma}

To construct a holomorphic extension $\rho$ of a representation $\pi$ we have
to assume that $\pi\in\mathcal{A}(W)$. Then for any $X\in W$, the operator
$i\,d\pi(X)$ generates a self adjoint contraction semigroup which we denote by%
\[
T_{X}(t)=e^{ti\,d\pi(X)}\,.
\]

For $s=g\exp iX\in S(C)$ we define%
\begin{equation}
\rho(s):=\rho(g)T_{X}(1) \label{extension}%
\end{equation}

\begin{theorem}
\label{Th:extension} $\rho$ is a contractive and holomorphic representation of
the semigroup $S(W)$. In particular, every representation $\pi\in
\mathcal{A}(W)$ extends uniquely to a holomorphic representation of $S(W)$
which is uniquely determined by $\pi$.
\end{theorem}

We will usually denote the holomorphic extension of the representation $\pi$
by the same letter. For the converse of Theorem \ref{Th:extension}, we remark
the following simple fact: Let $(S,\sharp)$ be a semigroup with involution and
let $\rho$ be a contractive representation of $S$. Let
\[
G(S):=\{s\in S\mid s^{\sharp}s=ss^{\sharp}=1\}
\]
Then $G(S)$ is a closed subgroup of $S$ and $\pi:=\rho|_{G(S)}$ is a unitary
representation of $G(S)$. Obviously,
\[
G\subset G(S(W))\,.
\]
Thus every holomorphic representation of $S(W)$ defines a unique unitary
representation of $G$ by restriction.

\begin{theorem}
\label{ThmS-HrepNew.6}Let $\rho$ be a holomorphic representation of $S(W)$.
Then $\rho|_{G}\in\mathcal{A}(W)$ and the $\rho$ agrees with the extension of
$\rho|_{G}$ to $S(W)$.
\end{theorem}

Two representations $\rho$ and $\pi$ of the semigroup $S(W)$ are said to (be
unitarily) equivalent if there exists a unitary isomorphism $U\colon
\mathbf{H}(\rho)\rightarrow\mathbf{H}(\pi)$ such that
\[
U\rho(s)=\pi(s)U\quad\forall s\in S(W)
\]
In particular, two contractive representations $\rho$ and $\pi$ of $S(W)$ are
equivalent if and only if $\rho|_{G}$ and $\pi|_{G}$ are unitarily equivalent.
We call a holomorphic contractive representation $\rho$ of $S(W)$
$W$-\textit{admissible} if $\rho|_{G}\in\mathcal{A}(W)$ and write $\rho
\in\mathcal{A}(W)$.

We denote by $\widehat{S(W)}$ the set of equivalence classes of irreducible
holomorphic representations of $S(W)$.

We say that a representation $\rho$ is bounded if
$\Vert\rho(s)\Vert\leq 1$
for all $s\in S(W)$. Note that this depends only on the unitary equivalence
class of $\rho$. We denote by $\widehat{S(W)}_{b}$ the subset in
$\widehat{S(W)}$ of bounded representations. Let $\rho$ and $\pi$ be
holomorphic representations of $S(W)$. Define a representation of $S(W)$ in
$\mathbf{H}(\rho)\hat{\otimes}\mathbf{H}(\pi)$ by%
\[
\lbrack\rho\otimes\pi](s):=\rho(s)\otimes\pi(s)
\]
Then $\rho\otimes\pi\in\mathcal{A}(W)$. We denote the representation
$s\mapsto\operatorname*{id}$ by $\iota$.

\begin{theorem}
[Neeb, Ol'shanskii \cite{HO95,Nee94,Ne99}]\label{DirInt} Let $\rho$ be a
holomorphic representation of $S(W)_{b}$.
Then there exists a Borel measure $\mu$ on $\widehat{S(W)}$ supported on
$\widehat{S(W)}_{b}$ and a direct integral of representations
\[
\left(  \int_{\widehat{S(W)}_{b}}^{\oplus}\rho_{\omega}\,d\mu(\omega
),\int_{\widehat{S(W)}_{b}}^{\oplus}\mathbf{H}(\omega)\,d\mu
(\omega)\right)
\]
such that:

\begin{enumerate}
\item \label{DirInt(1)}The representation $\rho$ is equivalent to
$\int_{\widehat{S(W)}_{b}}^{\oplus}\rho_{\omega}\,d\mu(\omega)$ .

\item \label{DirInt(2)}There exists a subset $N\subset\widehat{S(W)}_{b}$
such that $\mu(N)=0$ and if $\omega\in\widehat{S(W)}_{b}\setminus N$, then
$\rho_{\omega}$ is equivalent to $(\pi_{\omega}\otimes\iota,\mathbf{H}%
(\omega)\hat{\otimes}\mathbf{L}(\omega))$, where $\pi_{\omega}\in\omega$ and
$\mathbf{L}(\omega)$ is a Hilbert space.

\item \label{DirInt(3)}If $\omega\in\widehat{S(W)}_{b}$ then set
$n(\omega):=\dim\mathbf{L}(\omega)$. Then $n$ is a $\mu$-measurable function
from $\widehat{S(W)}_{b}$ to the extended positive axis $[0,\infty]$ which
is called the multiplicity function.
\end{enumerate}
\end{theorem}

\begin{proof}
See \cite{Ne99}, Theorem XI.6.13.
\end{proof}

\section{\label{LM}The L\"uscher-Mack Theorem}

\setcounter{equation}{0}

We use reference \cite{HiNe93} for the L\"uscher-Mack Theorem, but
\cite{FOS83}, \cite{GoJo83}, \cite{Jor86}, \cite{Jor87}, \cite{JoMo84},
\cite{KlLa83}, \cite{LM75}, and \cite{Sch86} should also be mentioned in this
connection. We have two ways of making the connection between the unitary
representations of $G$ and those of $G^{c}$: one is based on the
L\"uscher-Mack principle, and the other on the notion of local
representations from Jorgensen's papers \cite{Jor86} and \cite{Jor87}.

Let $\pi$, $C$, $\mathbf{H}(\pi)$, $J$ and $\mathbf{K}_{0}$ be as before. We
have proved that the operators
\[
\{\pi(h\exp(Y))\mid h\in H,Y\in C\}
\]
pass to the space $\mathbf{K}=\left(  \mathbf{K}_{0}/\mathbf{N}\right)
\sptilde$ such that $\tilde{\pi}(h)$ is unitary on $\mathbf{K}$, and
$\tilde{\pi}(\exp Y)$ is contractive and selfadjoint on $\mathbf{K}$. As a
result we arrive at selfadjoint operators $d\tilde{\pi}(Y)$ with spectrum in
$(-\infty,0]$ such that for $Y\in C$, $\tilde{\pi}(\exp Y)=e^{d\tilde{\pi}%
(Y)}$ on $\mathbf{K}$. As a consequence of that we notice that
\[
t\longmapsto e^{t\,d\tilde{\pi}(Y)}%
\]
extends to a continuous map on $\{z\in{{\mathbb{C}}}\mid\operatorname{Re}%
(z)\geq0\}$ holomorphic on the open right half plane $\{z\in{{\mathbb{C}}}%
\mid\operatorname{Re}(z)>0\}$. Furthermore,
\[
e^{(z+w)\,d\tilde{\pi}(Y)}=e^{z\,d\tilde{\pi}(Y)}e^{w\,d\tilde{\pi}(Y)}\,.
\]
As $\mathbf{K}$ is a unitary $H$-module we know that the $H$-analytic vectors
$\mathbf{K}^{\omega}(H)$ are dense in $\mathbf{K}$. Thus $\mathbf{K}%
_{oo}:=S(C^{o})\mathbf{K}^{\omega}(H)$ is dense in $\mathbf{K}$. We notice
that for $u\in\mathbf{K}_{oo}$ and $X\in C^{o}$ the function $t\mapsto
\tilde{\pi}(\exp tX)u$ extends to a holomorphic function on an open
neighborhood of the right half-plane. This and the Campbell-Hausdorff formula
are among the main tools used in proving the following Theorem of L\"uscher
and Mack \cite{LM75}. We refer to \cite[p. 292]{HiNe93} for the proof. Our
present use of Lie theory, cones, and semigroups will follow standard
conventions (see, e.g., \cite{FHO93,Hel62,JL94,WaI72,Yos91}): the exponential
mapping from the Lie algebra $\frak{g}$ to $G$ is denoted $\exp$, the adjoint
representation of $\frak{g}$, $\operatorname{ad}$, and that of $G$ is denoted
$\operatorname{Ad}$.

\begin{theorem}
[L\"uscher-Mack \cite{LM75}]\label{LuscherMack} Let $\rho$ be a strongly
continuous contractive representation of $S(C)$ on the Hilbert space
$\mathbf{H}$ such that $\rho(s)^{\ast}=\rho(\tau(s)^{-1})$. Let $G^{c}$ be the
connected, simply connected Lie group with Lie algebra $\frak{g}^{c}%
=\frak{h}\oplus i\frak{q}$. Then there exists a continuous unitary
representation $\rho^{c}\colon G^{c}\rightarrow\mathrm{U}(\mathbf{H})$,
extending $\rho$, such that for the differentiated representations $d\rho$ and
$d\rho^{c}$ we have:

\begin{enumerate}
\item [\hss\llap{\rm1)}]\label{LuscherMack(1)}$d\rho^{c}(X)=d\rho
(X)\,\quad\forall X\in\frak{h}$.

\item[\hss\llap{\rm2)}] \label{LuscherMack(2)}$d\rho^{c}(iY)=i\,d\rho
(Y)\,\quad\forall Y\in C$.
\end{enumerate}
\end{theorem}

We apply this to our situation to get the following theorem:

\begin{theorem}
\label{PiCIrreducible}Assume that $(\pi,C,\mathbf{H},J)$ satisfies
\textup{(PR1)--(PR3).} Then the following hold:

\begin{enumerate}
\item [\hss\llap{\rm1)}]\label{PiCIrreducible(1)}$S(C)$ acts via
$s\mapsto\tilde{\pi}(s)$ by contractions on $\mathbf{K}$.

\item[\hss\llap{\rm2)}] \label{PiCIrreducible(2)}Let $G^{c}$ be the simply
connected Lie group with Lie algebra $\frak{g}^{c}$. Then there exists a
unitary representation $\tilde{\pi}^{c}$ of $G^{c}$ such that $d\tilde{\pi
}^{c}(X)=d\tilde{\pi}(X)$ for $X\in\frak{h}$ and $i\,d\tilde{\pi}%
^{c}(Y)=d\tilde{\pi}(iY)$ for $Y\in C$.

\item[\hss\llap{\rm3)}] \label{PiCIrreducible(3)}The representation
$\tilde{\pi}^{c}$ is irreducible if and only if $\tilde{\pi}$ is irreducible.
\end{enumerate}
\end{theorem}

\begin{proof}
(1) and (2) follow by the L\"uscher-Mack theorem and Proposition 3.6, as the
resulting representation of $S$ is obviously continuous.

(3) Let $\mathbf{L}$ be a $G^{c}$-invariant subspace in $\mathbf{K}$. Then
$\mathbf{L}$ is $\tilde{\pi}(H)$ invariant. Let $Y\in C^{o}$, $u\in
\mathbf{L}^{\omega}$ and $v\in\mathbf{L}^{\perp}$. Define $f\colon
\{z\in{{\mathbb{C}}}\mid\operatorname{Re}(z)\geq0\}\rightarrow{{\mathbb{C}}}$
by
\[
f(z):=%
\ip{v}{e^{zd\tilde{\pi}(Y)}u}%
_{J}\,.
\]
Then $f$ is holomorphic in $\{z\in{{\mathbb{C}}}\mid\operatorname{Re}(z)>0\}$,
and $f(it)=0$ for every (real) $t$. Thus $f$ is identically zero. In
particular $f(t)=0$ for every $t>0$. Thus
\[
0=%
\ip{v}{e^{td\tilde{\pi}(Y)}u}%
_{J}=%
\ip{v}{\tilde{\pi}(\exp tY)u}%
_{J}\,.
\]
As $S^{o}=H\exp C^{o}$ it follows that $\tilde{\pi}(S^{o})(\mathbf{L}^{\omega
})\subset(\mathbf{L}^{\perp})^{\perp}=\mathbf{L}$. By continuity we get
$\tilde{\pi}(S)\mathbf{L}\subset\mathbf{L}$. Thus $\mathbf{K}$ is reducible as
an $S$-module.

The other direction follows in exactly the same way.
\end{proof}

We notice now that $-iC\subset W(\tilde{\pi}^{c})$. Thus $W(\tilde{\pi}^{c})$
is non-trivial and contain the $-\tau$-stable and $G$-invariant cone generated
by $-iC$, i.e. $-i\cdot\overline{\operatorname{conv}\{\operatorname{Ad}%
(G)C\}}\subset W(\tilde{\pi}^{c})$. But in general $W(\tilde{\pi}^{c})$ is
neither generating nor pointed. It even does not have to be $-\tau$-invariant.
In fact, the Lie algebra of the $(ax+b)$-group, and the Heisenberg group, do
not have \emph{any} pointed, generating, invariant cones.

\begin{lemma}
\label{KernelPiC}$W(\tilde{\pi}^{c})\cap-W(\tilde{\pi}^{c})=\ker(\tilde{\pi
}^{c})$.
\end{lemma}

\begin{proof}
This is obvious from the spectral theorem.
\end{proof}

\begin{lemma}
\label{IdealGC}$\frak{g}_{1}^{c}:=W(\tilde{\pi}^{c})-W(\tilde{\pi}^{c})$ is an
ideal in $\frak{g}^{c}$. Furthermore, $[\frak{q},\frak{q}]\oplus
i\frak{q}\subset\frak{g}_{1}^{c}$.
\end{lemma}

\begin{proof}
Let $X\in\frak{g}^{c}$. Then, as $W(\tilde{\pi}^{c})$ is invariant by
construction, we conclude that
\[
e^{t\operatorname{ad}(X)}\left(  W(\tilde{\pi}^{c})-W(\tilde{\pi}^{c})\right)
\subset W(\tilde{\pi}^{c})-W(\tilde{\pi}^{c}),\,t\in\mathbb{R}\,.
\]
By differentiation at $t=0$, it follows that $[X,\frak{g}_{1}^{c}%
]\subset\frak{g}_{1}^{c}$. This shows that $\frak{g}_{1}^{c}$ is an ideal in
$\frak{g}^{c}$. The last part follows as $C$ is generating (in $\frak{q}$).
\end{proof}

It is not clear if $\frak{g}_{1}^{c}$ is $\tau$-stable. To get a $\tau$-stable
subalgebra one can replace $W(\tilde{\pi}^{c})$ by the cone generated by
$-\operatorname{Ad}(G)C\subset W(\tilde{\pi}^{c})$ or by the maximal $G$- and
$-\tau$-stable cone $W(\tilde{\pi}^{c})\cap(-\tau(W(\tilde{\pi}^{c})))$ in
$W(\pi^{c})$.

We have now the following important consequence of the Neeb-Ol'shanskii theorem:

\begin{theorem}
\label{HighestWeight}Let the analytic subgroup $G_{1}^{c}$ of $G^{c}$
corresponding to $\frak{g}_{1}^{c}$ be as described, and let $W(\tilde{\pi
}^{c})$ be the corresponding module. Then $\tilde{\pi}^{c}|_{G_{1}^{c}}$ is a
direct integral of irreducible representations in $\mathcal{A}(W)$.
\end{theorem}

\section{\label{SSS}Bounded Symmetric Domains}

\setcounter{equation}{0}

We have seen that the representations $\pi^{c}$ that we can produce using the
duality are direct integrals of holomorphic representations of suitable
subsemigroups of $G_{\mathbb{C}}^{c}$ (or a subgroup). Those on the
other hand only exist if there is a $G^{c}$ invariant cone in $\frak{g}$. We
will discuss the case of simple Lie group $G^{c}$ in some detail here. We
refer to \cite{JO97}, chapter 2, and the references therein for proofs. In the
duality $G\leftrightarrow G^{c}$ it will be the group $G^{c}$ that has
holomorphic representations. Therefore we will start using the notation
$G^{c}$ for Hermitian groups.

\begin{theorem}
[Kostant]\label{ThmKostant}Suppose that $\mathbf{V}$ is a finite-dimensional
real vector space. Let $L$ be a connected reductive subgroup of $GL(\mathbf{V}%
)$ acting irreducibly on $\mathbf{V}$. Let $G^{c}=L^{\prime}$ be the
commutator subgroup of $L$. Further let $K^{c}$ be a maximal compact
subgroup of $G^{c}$. Then the following properties are equivalent:

\begin{enumerate}
\item \label{ThmKostant(1)}There exists a regular $L$-invariant closed cone in
$\mathbf{V}$.

\item \label{ThmKostant(2)}The $G^{c}$-module $\mathbf{V}$ is spherical.
\end{enumerate}
\end{theorem}

Let $C\subset\mathbf{V}$ be a regular $L$-invariant cone. Then $\mathbf{V}$ is
spherical as a $G^{c}$-module. Let $K^{c}$ be a maximal compact subgroup of
$G^{c}$. A $K^{c}$-invariant vector $u_{K^{c}}$ can be constructed in the
following way: Let $u\in C^{o}$, the interior of $C$, be arbitrary. Define
\[
u_{K^{c}}=\int_{K^{c}}k\cdot u\,dk\,.
\]
Then $u_{K^{c}}\in C^{o}$ is $K^{c}$-invariant. Suppose that the group $L$
acts on $\mathbf{V}$. Let $\operatorname{Cone}_{L}(\mathbf{V})$ denote the set
of regular $L$-invariant cones in $\mathbf{V}$.

\begin{theorem}
[Vinberg]\label{ThmVinberg}Let $L$, $G^{c}$, and $\mathbf{V}$ be as in the
theorem of Kostant. Then the following properties are equivalent:

\begin{enumerate}
\item \label{ThmVinberg(1)}$\operatorname{Cone}_{L}(\mathbf{V})\not =%
\emptyset$;

\item \label{ThmVinberg(2)}The $G^{c}$-module $\mathbf{V}$ is spherical;

\item \label{ThmVinberg(3)}There exists a ray in $\mathbf{V}$ through $0$
which is invariant with respect to some minimal parabolic subgroup $P$ of
$G^{c}$.
\end{enumerate}

If those conditions hold, every invariant pointed cone in $\mathbf{V}$ is regular.
\end{theorem}

For the next theorem, see \cite{Paneitz81,Paneitz84}.

\begin{theorem}
[Paneitz, Vinberg]\label{ThmPV} Let $G^{c}$ be a
connected semisimple Lie group and $(\mathbf{V,\pi)}$ a real
finite-dimensional irreducible $G^{c}$-module such that $\operatorname{Cone}%
_{G^{c}}(\mathbf{V})\not =\emptyset$.  Let $\theta$ be a Cartan involution on
$G^{c}$. Choose an inner product on $\mathbf{V}$ such that $\pi(x)^{\ast}%
=\pi(\theta(x)^{-1})$ for all $x\in G^{c}$. Then there exists a unique up to
multiplication by $(-1)$ invariant cone $C_{\min}\in\operatorname{Cone}%
_{G^{c}}(\mathbf{V})$ given by
\[
C_{\min}=\operatorname{conv}(\pi(G^{c})u)\cup\left\{  0\right\}
=\overline{\operatorname{conv}\{\pi(G^{c})(\mathbb{R}^{+}v_{K^{c}})\}}\,,
\]
where $u$ is a highest weight vector, $v_{K^{c}}$ is a nonzero $K^{c}$-fixed
vector unique up to scalar multiple, and $(u,v_{K^{c}})>0$. The unique
\textup{(}up to
multiplication by $(-1)$\textup{)} maximal cone is given by
\[
C_{\max}=C_{\min}^{\ast}:=\left\{  w\in\mathbf{V}\mid\forall v\in C_{\min
}\,:\,(w,v)\geq0\right\}  \,.
\]
\end{theorem}

Assume now that $G^{c}$ is a connected simple Lie group. Then $G^{c}$ acts on
$\frak{g}^{c}$ by the adjoint action. Let $K^{c}\subset G^{c}$ be a maximal
almost compact subgroup. Then by Kostant's Theorem we have
$\operatorname{Cone}_{G^{c}}(\frak{g}^{c})\not =\emptyset$ if and only if
there exists a $Z^{0}\in\frak{g}^{c}$ which is invariant under
$\operatorname{Ad}(K^{c})$. Let $\frak{k}^{c}$ be the Lie algebra of $K^{c}$.
Then $[\frak{k}^{c},Z^{0}]=0$. Hence $\mathbb{R}Z^{0}+\frak{k}^{c}%
\not =\frak{g}^{c}$ is a Lie algebra containing $\frak{k}^{c}$. But
$\frak{k}^{c}$ is maximal in $\frak{g}^{c}$. Hence $Z^{0}\in\frak{k}^{c}$.
Similarly it follows that $\frak{z}_{\frak{g}^{c}}(Z^{0})=\frak{k}^{c}$.
Finally the Theorem of Paneitz and Vinberg implies that $\operatorname{Cone}%
_{G^{c}}(\frak{g}^{c})\not =\emptyset$ if and only if the center of
$\frak{k}^{c}$ is one dimensional. In that case we can normalize the element
$Z^{0}$ such that $\operatorname{ad}(Z^{0})$ has eigenvalues $0,i,-i$. Let
$\frak{t}$ be a Cartan subalgebra of $\frak{g}^{c}$ containing $Z^{0}$. Then
$\frak{k}^{c}\frak{\subset z}_{\frak{g}^{c}}(Z^{0})\subset\frak{k}^{c}$. Hence
$\frak{t}$ is contained in $\frak{k}^{c}$. For $\alpha\in\frak{t}_{\mathbb{C}%
}^{\ast}$ let
\[
\frak{g}_{{\mathbb{C}\alpha}}^{c}:=\left\{  X\in\frak{g}_{\mathbb{C}}^{c}%
\mid\forall Z\in\frak{t}_{\mathbb{C}}\,:\,[Z,X]=\alpha(Z)X\right\}  \,.
\]
It is well known that $\dim\frak{g}_{{\mathbb{C}\alpha}}^{c}\leq1$ for all
$\alpha\not =0$ and
\[
\frak{g}_{\mathbb{C}}^{c}=\frak{t}_{\mathbb{C}}\oplus\bigoplus_{\alpha
\in\Delta}\frak{g}_{{\mathbb{C}\alpha}}^{c}%
\]
where $\Delta=\left\{  \alpha\in\frak{t}_{\mathbb{C}}^{\ast}\setminus\left\{
0\right\}  \mid\frak{g}_{{\mathbb{C}\alpha}}^{c}\not =\left\{  0\right\}
\right\}  $. We notice that $\alpha(\frak{t})\subset i\mathbb{R}$ for all
$\alpha\in\Delta$ as $\frak{t}\subset\frak{k}$ and $a\,d(X)$ is skew-symmetric
for all $X\in\frak{k}$. Let $\theta\colon\frak{g}^{c}\frak{\rightarrow
}\frak{g}^{c}$ be the Cartan involution corresponding to $\frak{k}^{c}$. We do
denote the corresponding involution $\theta\otimes1$ on $\frak{g}_{\mathbb{C}%
}^{c}$ and the integrated involution on $G^{c}$ by the same letter. Then
$\frak{g}^{c}=\frak{k}^{c}\oplus\frak{p}^{c}$ where $\frak{p}^{c}%
\frak{=}\left\{  X\in\frak{g}^{c}\mid\theta(X)=-X\right\}  $. Thus
$\frak{g}_{\mathbb{C}}^{c}=\frak{k}_{\mathbb{C}}^{c}\oplus\frak{p}%
_{\mathbb{C}}^{c}$. Let $\Delta_{c}=\left\{  \alpha\in\Delta\mid\alpha
(Z^{0})=0\right\}  $ and $\Delta_{p}=\left\{  \alpha\in\Delta\mid\alpha
(Z^{0})=\pm i\right\}  $. As $\frak{z}_{\frak{g}_{\mathbb{C}}}(Z^{0}%
)=\frak{k}_{\mathbb{C}}$ we get%
\[
\Delta_{c}=\left\{  \alpha\mid\frak{g}_{{\mathbb{C}\alpha}}^{c}\subset
\frak{k}_{\mathbb{C}}^{c}\right\}  \qquad\text{and}\qquad\Delta_{p}=\left\{
\alpha\mid\frak{g}_{{\mathbb{C}\alpha}}^{c}\subset\frak{p}_{\mathbb{C}}%
^{c}\right\}  \,.
\]
Choose a positive system $\Delta^{+}$ in $\Delta$ such that $\Delta_{p}%
^{+}:=\left\{  \alpha\in\Delta_{p}\mid\alpha(Z^{0})=i\right\}  \subset
\Delta^{+}$. Then $\Delta^{+}=\Delta_{c}^{+}\cup\Delta_{p}^{+}$ and
$\Delta_{c}^{+}$ is a positive system in $\Delta_{c}$. For $\Gamma
\subset\Delta$ let $\frak{g}_{\mathbb{C}}^{c}(\Gamma):=\bigoplus_{\alpha
\in\Gamma}\frak{g}_{\mathbb{C\alpha}}^{c}$. Then%
\begin{align*}
\frak{p}^{+}:=  &  \left\{  X\in\frak{g}_{\mathbb{C}}^{c}\mid\lbrack
Z^{0},X]=iX\right\}  =\frak{g}_{\mathbb{C}}^{c}(\Delta_{p}^{+})\,;\\
\frak{p}^{-}:=  &  \left\{  X\in\frak{g}_{\mathbb{C}}^{c}\mid\lbrack
Z^{0},X]=-iX\right\}  =\frak{g}_{\mathbb{C}}^{c}(-\Delta_{p}^{+})\,.
\end{align*}
Furthermore $\frak{p}^{+}$ and $\frak{p}^{-}$ are abelian subalgebras with
$\frak{p}_{\mathbb{C}}=\frak{p}^{+}\oplus\frak{p}^{-}$. Let $P^{\pm}%
:=\exp(\frak{p}^{\pm})$ and $K_{\mathbb{C}}^{c}=\exp(\frak{k}_{\mathbb{C}}%
^{c})$. Both $P^{+}$ and $P^{-}$ are simply connected closed abelian subgroups
of $G_{\mathbb{C}}^{c}$. Hence $\exp\colon\frak{p}^{\pm}\rightarrow P^{\pm}$
is a diffeomorphism. Let%
\begin{equation}
\zeta=(\exp|_{\frak{p}^{+}})^{-1}\colon P^{+}\longrightarrow\frak{p}%
^{+}\,\text{.} \label{E:zeta}%
\end{equation}
The set $P^{+}K_{\mathbb{C}}^{c}P^{-}$ is open and dense in $G_{\mathbb{C}%
}^{c}$, $G^{c}\subset P^{+}K_{\mathbb{C}}^{c}P^{-}$, $G^{c}K_{\mathbb{C}}%
^{c}P^{-}$ is open in $G_{\mathbb{C}}^{c}$, and $G^{c}\cap K_{\mathbb{C}}%
^{c}P^{-}=K^{c}$. Thus $G^{c}/K^{c}$ is holomorphically equivalent to an
open submanifold $D$ of the complex flag manifold $X_{\mathbb{C}%
}=G_{\mathbb{C}}^{c}/K_{\mathbb{C}}^{c}P^{-}$. Furthermore the map
$pK_{\mathbb{C}}^{c}P^{-}\mapsto\zeta(p)$ induces a biholomorphic map---also
denoted by $\zeta$---of $G^{c}/K^{c}$ onto a bounded symmetric domain
$\Omega_{\mathbb{C}}\subset\frak{p}^{-}\simeq\mathbb{C}^{\dim(G^{c}/K^{c})}$.

For $x\in P^{+}K_{\mathbb{C}}^{c}P^{-}$ we can write in a unique way%
\begin{equation}
x=p(x)k_{\mathbb{C}}(x)q(x) \label{E:P+KP-projections}%
\end{equation}
with $p(x)\in P^{+}$, $k_{\mathbb{C}}(x)\in K_{\mathbb{C}}^{c}$ and $q(x)\in
P^{-}$. For $g\in G_{\mathbb{C}}^{c}$ and $Z\in\frak{p}^{+}$ we introduce the
following notations when ever they make sense:%
\begin{align*}
g\cdot Z  &  =\zeta(p(g\exp(Z))\in\frak{p}^{+}\\
j(g,Z)  &  =k_{\mathbb{C}}(g\exp Z)\in K_{\mathbb{C}}^{c}\,.
\end{align*}
If $Z\in\Omega_{\mathbb{C}}$ and $g\in G^{c}$ then $g\cdot Z$ is defined and
$g\cdot Z\in\mathbb{\Omega}_{\mathbb{C}}$. Furthermore $(g,Z)\mapsto g\cdot Z$
defines an action of $G^{c}$ on $\mathbb{\Omega}_{\mathbb{C}}$ such that
$\zeta\colon G^{c}/K^{c}\rightarrow\Omega_{\mathbb{C}}$ is a $G^{c}$-map. The
map $j$ is \textit{the universal automorphic factor} and it satisfies the
following:%
\begin{align}
j(k,Z)  &  =k,\nonumber\\
j(p,Z)  &  =1,\nonumber\\
j(ab,Z)  &  =j(a,b\cdot Z)j(b,Z), \label{E:Multipl}%
\end{align}
if $k\in K_{\mathbb{C}}^{c}$, $Z\in\frak{p}^{+}$, $p\in P^{+}$, and $a,b\in
G_{\mathbb{C}}^{c}$ are such that the expressions above are defined.

Define%
\[
S(\Omega_{\mathbb{C}}):=\left\{  \gamma\in G_{\mathbb{C}}^{c}\mid\gamma
^{-1}\cdot\Omega_{\mathbb{C}}\subset\Omega_{\mathbb{C}}\right\}
\]
and%
\[
S(\Omega_{\mathbb{C}})^{o}:=\left\{  \gamma\in G_{\mathbb{C}}^{c}\mid
\gamma^{-1}\cdot\overline{\Omega_{\mathbb{C}}}\subset\Omega_{\mathbb{C}%
}\right\}
\]
where $\overline{\left\{  \Omega_{\mathbb{C}}\right\}  }$ stands for the
topological closure of $\Omega_{\mathbb{C}}$ in $\frak{p}^{+}$. Then
$S(\Omega_{\mathbb{C}})$ is a closed semigroup in $\frak{g}_{\mathbb{C}}^{c}$
of the form%
\[
S(\Omega_{\mathbb{C}})=G^{c}\exp(iC_{\max})\,
\]
where $C_{\max}$ is the maximal $G^{c}$-invariant cone in $\frak{g}%
_{\mathbb{C}}^{c}$ containing $-Z^{0}$. Furthermore $S(\Omega_{\mathbb{C}%
})^{o}$ is the topological interior of $S(\Omega_{\mathbb{C}})$ and%
\[
S(\Omega_{\mathbb{C}})^{o}=S(C_{\max}^{o})=G^{c}\exp(iC_{\max}^{o})\,.
\]
We refer to \cite{HO95} or \cite{HiNe93} for all of this.

\section{\label{S-hwm}Highest Weight Modules}

Our notion of reflection for unitary representations leads to the class of
representations in $\mathcal{A}(W)$ of Definition \ref{DefS-HrepMar.1}, and in
the present section we analyze these representations more closely. The
analysis is based in large part on \cite{HO95}, and involves results of (among
others) M.~Davidson and R.~Fabec \cite{MDRF96}, K.-H. Neeb
\cite{Nee94,Ne94,Ne99}, Harish-Chandra \cite{HC65,HC70}, Ol'shanskii
\cite{Ol82}, R.J. Stanton \cite{Stanton}, Wallach \cite{WaI72}, and H.~Rossi
and M.~Vergne \cite{VR76}.

We have seen that the interesting representations are those in $\mathcal{A}%
(W)$ where $W$ is an invariant cone in $\frak{g}^{c}$. It turns out that the
irreducible representations in $\mathcal{A}(W)$ are \textit{highest weight
representations. }A $(\frak{g}^{c},K^{c})$-module\label{gK-modul} is a complex
vector space $\mathbf{V}$ such that

\begin{enumerate}
\item [1)]$\mathbf{V}$ is a $\frak{g}^{c}$-module.

\item[2)] $\mathbf{V}$ carries a representation of $K^{c}$, and the span of
$K^{c}\cdot v$ is finite-dimensional for every $v\in\mathbf{V}$.

\item[3)] For $v\in\mathbf{V}$ and $X\in\frak{k}^{c}$ we have
\[
X\cdot v=\lim_{t\rightarrow0}\frac{\exp(tX)\cdot v-v}{t}.
\]

\item[4)] For $Y\in\frak{g}^{c}$ and $k\in K^{c}$ the following holds for
every $v\in\mathbf{V}$:
\[
k\cdot(X\cdot v)=(\operatorname{Ad}(k)X)\cdot\lbrack k\cdot v]\,.
\]
\end{enumerate}

Note that (3) makes sense, as $K^{c}\cdot v$ is contained in a finite
dimensional vector space and this space contains a unique Hausdorff topology
as a topological vector space. The $(\frak{g}^{c},K^{c})$-module is called
\textit{admissible}\label{st-ad2} if the multiplicity of every irreducible
representation of $K^{c}$ in $\mathbf{V}$ is finite. If $(\pi,\mathbf{V})$ is
an irreducible unitary representation of $G^{c}$, then the space of $K^{c}%
$-finite elements in $\mathbf{V}$, denoted by $\mathbf{V}_{K^{c}}$, is an
admissible $(\frak{g}^{c},K^{c})$-module.

Let $\frak{t}$ be a Cartan subalgebra of $\frak{k}^{c}$ and $\frak{g}^{c}$ as
in the last section.

\begin{definition}
\label{DefNewS-hwm.1}Let $\mathbf{V}$ be a $(\frak{g}^{c},K^{c})$-module. Then
$\mathbf{V}$ is a highest-weight module if there exists a nonzero element
$v\in\mathbf{V}$ and a $\lambda\in\frak{t}_{\mathbb{C}}^{\ast}$ such that

\begin{enumerate}
\item [1)]\label{DefNewS-hwm.1(1)}$X\cdot v=\lambda(X)v$ for all $X\in
\frak{t}$.

\item[2)] \label{DefNewS-hwm.1(2)}There exists a positive system $\Delta^{+}$
in $\Delta$ such that $\frak{g}_{\mathbb{C}}^{c}(\Delta^{+})\cdot v=0$.

\item[3)] \label{DefNewS-hwm.1(3)}$\mathbf{V}=U(\frak{g}^{c})\cdot v$.
\end{enumerate}

The element $v$ is called a primitive element of weight $\lambda$.
\end{definition}

Let $W\in\operatorname{Cone}_{G^{c}}(\frak{g}^{c})$ and $\pi\in\mathcal{A}(W)$
irreducible. We assume that $-Z^{o}\in W^{o}$. Then $\mathbf{V}_{K^{c}}$ is an
irreducible admissible $(\frak{g}^{c},K^{c})$-module, and
\[
\mathbf{V}_{K^{c}}=\bigoplus_{\lambda\in\frak{t}_{\mathbb{C}}^{\ast}%
}\mathbf{V}_{K^{c}}(\lambda)
\]
where $\mathbf{V}_{K^{c}}(\lambda)=\mathbf{V}_{K^{c}}(\lambda,\frak{t}%
_{\mathbb{C}})$. Let $v\in\mathbf{V}_{K^{c}}(\lambda)$ be nonzero. Let
$\alpha\in\Delta_{p}^{+}$ and let $X\in\frak{p}_{\alpha}^{+}\setminus\{0\}$.
Then
\[
X^{k}\cdot v\in\mathbf{V}_{K^{c}}(\lambda+k\alpha).
\]
In particular,
\[
-iZ^{0}\cdot(X^{k}\cdot v)=[-i\lambda(Z^{0})+k]v.
\]

This yields the following lemma.

\begin{lemma}
\label{LemNewS-hwm.2}Let the notation be as above. Then the following holds:

\begin{enumerate}
\item \label{LemNewS-hwm.2(1)}$-i\lambda(Z^{0})\leq0$.

\item \label{LemNewS-hwm.2(2)}There exists a $\lambda$ such that $\frak{p}%
^{+}\cdot\mathbf{V}_{K^{c}}(\lambda)=\{0\}$.
\end{enumerate}
\end{lemma}

\begin{lemma}
\label{LemNewS-hwm.3}Let $\mathbf{W}^{\lambda}$ be the $K^{c}$-module
generated by $\mathbf{V}_{K^{c}}(\lambda)$. Then $\mathbf{W}^{\lambda}$ is
irreducible, $\mathbf{V}_{K^{c}}=U(\frak{p}^{-})\mathbf{W}^{\lambda}$ and the
multiplicity of $\mathbf{W}^{\lambda}$ in $\mathbf{V}_{K^{c}}$ is one.
\end{lemma}

Let $\alpha\in\Delta_{p}^{+}$ then there exists a unique element $H_{\alpha
}\in i\frak{t}\cap\lbrack\frak{g}_{\mathbb{C\alpha}},\frak{g}_{\mathbb{C}%
\alpha}]$ such that $\alpha(H_{\alpha})=2$. Let $\mu$ be the highest weight of
$\mathbf{W}^{\lambda}$ with respect to $\Delta_{c}^{+}$ and let $v^{\lambda}$
be a nonzero highest weight vector. Then $v^{\lambda}$ is a primitive element
with respect to the positive system $\Delta_{c}^{+}\cup\Delta_{p}^{+}$.

\begin{theorem}
\label{AChw} Let $\rho\in\mathcal{A}(W)$ be irreducible. Then the
corresponding $(\frak{g}^{c},K^{c})$-module is a highest-weight module and
equals $U(\frak{p}^{-})\mathbf{W}^{\lambda}$. In particular, every weight of
$\mathbf{V}_{K^{c}}$ is of the form
\[
\nu-\sum_{\alpha\in\Delta(\frak{p}^{+},\frak{t}_{\mathbb{C}})}n_{\alpha}%
\alpha\,.
\]
Furthermore, $\left\langle \nu,H_{\alpha}\right\rangle \leq0$ for all $\alpha\in\Delta_{p}^{+}$.
\end{theorem}

The $K^{c}$-representation $\pi^{\lambda}$ on $\mathbf{W}^{\lambda}$ is called
the \emph{minimal}\textit{ }$K^{c}$-type of $\mathbf{V}$ and $\mathbf{V}%
_{K^{c}}$. The multiplicity of $\pi^{\lambda}$ in $\mathbf{V}$ is one. We
recall how to realize highest-weight modules in a space of holomorphic
functions on $G^{c}/K^{c}$. We follow here the geometric construction by M.
Davidson and R. Fabec \cite{MDRF96}. For a more general approach, see
\cite{Ne94,Ne99}. To explain the method we start with the example
$G^{c}=SU(1,1)=\left\{  \left(
\begin{array}
[c]{cc}%
\alpha & \beta\\
\bar{\beta} & \bar{\alpha}%
\end{array}
\right)  \biggm|\left|  \alpha\right|  ^{2}-\left|  \beta\right|
^{2}=1\right\}  $. We set $X=\left(
\begin{array}
[c]{cc}%
0 & 1\\
0 & 0
\end{array}
\right)  $, $Y=\left(
\begin{array}
[c]{cc}%
0 & 0\\
1 & 0
\end{array}
\right)  $ and $H:=H_{1}=\left(
\begin{array}
[c]{cc}%
1 & 0\\
0 & -1
\end{array}
\right)  $. Then $Z^{0}=\frac{i}{2}H$ and $\frak{p}^{+}=\mathbb{C}X$,
$\frak{k}_{\mathbb{C}}^{c}=\mathbb{C}H$ and $\frak{p}^{-}={\ \mathbb{C}}Y$. We
use this to identify those spaces with $\mathbb{C}$. Let $\gamma=\left(
\begin{array}
[c]{cc}%
a & b\\
c & d
\end{array}
\right)  \in SL(2,\mathbb{C})\,$. Then%
\begin{align*}
\left(
\begin{array}
[c]{cc}%
a & b\\
c & d
\end{array}
\right)   &  =\left(
\begin{array}
[c]{cc}%
1 & z\\
0 & 1
\end{array}
\right)  \left(
\begin{array}
[c]{cc}%
\gamma & 0\\
0 & \gamma^{-1}%
\end{array}
\right)  \left(
\begin{array}
[c]{cc}%
1 & 0\\
y & 1
\end{array}
\right) \\
&  =\left(
\begin{array}
[c]{cc}%
\gamma+\gamma^{-1}z & \gamma^{-1}z\\
\gamma^{-1}y & \gamma^{-1}%
\end{array}
\right)  \,.
\end{align*}
Hence $P^{+}K_{\mathbb{C}}^{c}P^{-}=\left\{  \left(
\begin{array}
[c]{cc}%
a & b\\
c & d
\end{array}
\right)  \biggm|d\not =0\right\}  $ and if $x=\left(
\begin{array}
[c]{cc}%
a & b\\
c & d
\end{array}
\right)  \in P^{+}K_{\mathbb{C}}^{c}P^{-}$, then%
\[
p(x)=b/d\,,\quad k_{\mathbb{C}}(x)=d^{-1}\quad\text{and}\quad
q(x)=c/d\,.
\]
Thus
\begin{align}
\zeta(xK_{\mathbb{C}}^{c}P^{-})  &  =b/d\\
x\cdot z  &  =\frac{az+b}{cz+d}\\
j(x,z)  &  =(cz+d)^{-1}\,. \label{E:jsu}%
\end{align}
To identify $\Omega_{\mathbb{C}}$ we notice that on $SU(1,1)$ we have
$\zeta(x)=\beta/\bar{\alpha}$. Hence $G^{c}/K^{c}\simeq D=\{z\in\mathbb{C}%
\mid|z|<1\}$. The finite-dimensional holomorphic representations of
$K_{\mathbb{C}}^{c}$ are the characters
\[
\chi_{n}(\exp ziH)=e^{inz}\,.
\]
In particular, $d\chi_{n}(Z^{0})=in/2$ or
\[
-i\,d\chi_{n}(Z^{0})=\frac{n}{2}\,.
\]
Let $(\pi,\mathbf{V})$ be a unitary highest-weight representation of $SU(1,1)$
and assume that $(\pi,\mathbf{V})\in\mathcal{A}(W)$. Then $n\leq0$ by Lemma
\ref{LemNewS-hwm.2} and Theorem \ref{AChw}. Let $\mathbf{V}(n)$ be the
one-dimensional space of $\chi_{n}$-isotropic vectors. Then
\[
\mathbf{V}_{K^{c}}=\bigoplus_{k\in\,\mathbb{N}}\mathbf{V}(n-2k),
\]
and the spaces $\mathbf{V}(m)$ and $\mathbf{V}(k)$ are orthogonal if
$m\not =k$.

Let $\sigma$ be the conjugation of $\frak{s}\frak{l}(2,\mathbb{C})$ with
respect to $SU(1,1)$. Then $\sigma$ is given by%
\[
\sigma\left(  \left(
\begin{array}
[c]{cc}%
a & b\\
c & d
\end{array}
\right)  \right)  =\left(
\begin{array}
[c]{cc}%
-\bar{a} & \bar{c}\\
\bar{b} & -\bar{a}%
\end{array}
\right)
\]
so that $\sigma(X)=Y$.
Since
$\pi(T)^{\ast}=-\pi(\sigma(T))$ for all
$T\in\frak{s}\frak{l}(2,\mathbb{C})$ we get%
\[
\pi(Y)^{\ast}=\pi(-X).
\]
Finally, it follows from $[Y,X]=-H$ that for $v\in V(n)$:
\begin{align*}
\Vert\pi(Y)^{k}v\Vert^{2}  &  =%
\ip{\pi(Y)^{k}v}{\pi(Y)^{k}v}%
\\
&  =%
\ip{\pi((-X)^{k}Y^{k})v}{v}%
\end{align*}

\begin{lemma}
\label{LemNewS-hwm.5}Let the notation be as above. Then
\begin{align*}
\pi(-X)^{k}\pi(Y)^{k}v  &  =(-1)^{k}k!\frac{\Gamma(n+1)}{\Gamma(n-k+1)}\,v\\
&  =(-n)_{k}v
\end{align*}
where $(a)_{k}=a(a+1)\cdots(a+k-1)$.
\end{lemma}

As
\[
\sum_{k=0}^{\infty}(-n)_{k}\frac{|z^{2}|^{k}}{k!}=(1-|z|^{2})^{n},
\]
(see\ \cite{GR80}) converges if and only if $|z|<1$, it follows that
\[
q_{(zX)}v:=\sum_{k=0}^{\infty}\overline{z}^{k}\frac{Y^{k}v}{n!}%
\]
converges if and only if $zX\in{\Omega}_{\mathbb{C}}$.

Let now $G^{c}$ be arbitrary. Let $\sigma\colon\frak{g}_{\mathbb{C}}%
^{c}\rightarrow\frak{g}_{\mathbb{C}}^{c}$ be the conjugation with respect to
$\frak{g}^{c}$. We use the notation from earlier in this section. Using the
usual $\frak{sl}(2,\mathbb{C})$ reduction, we get the following theorem.

\begin{theorem}
[Davidson-Fabec]\label{ThDF1} Let $T\in\frak{p}^{+}$. Define $q_{T}%
\colon\mathbf{W}^{\lambda}\rightarrow\mathbf{V}$ by the formula
\[
q_{T}v:=\sum_{n=0}^{\infty}\frac{\sigma(T)^{n}v}{n!}\,.
\]

\begin{enumerate}
\item [1)]If $v\not =0$, then the series that defines $q_{T}$ converges in the
Hilbert space $\mathbf{V}$ if and only if $T\in{\Omega}_{\mathbb{C}}$.

\item[2)] Let $\pi_{\lambda}$ be the representation of $K^{c}$ on
$\mathbf{W}^{\lambda}$. Let%
\begin{equation}
J_{\lambda}(g,Z):=\pi_{\lambda}(j(g,Z))\,. \label{E:jlambda}%
\end{equation}
Then
\begin{equation}
\pi(g)v=q_{g\cdot0}J_{\lambda}(g,0)^{\ast-1}v \label{E:pilambda}%
\end{equation}
for $g\in G^{c}$ and $v\in\mathbf{W}^{\lambda}$.
\end{enumerate}
\end{theorem}

It follows that the span of the $q_{Z}\mathbf{W}^{\lambda}$ with $Z\in{\Omega
}_{{\mathbb{C}}}$ is dense in $\mathbf{V}$, since $\mathbf{V}$ is assumed to
be irreducible. Define $Q\colon{\Omega}_{\mathbb{C}}\times{\Omega}%
_{\mathbb{C}}\rightarrow GL(\mathbf{W}^{\lambda})$ by%
\[
Q(W,Z)=q_{W}^{\ast}q_{Z}\,.
\]
\label{QWZ}

Then the following theorem holds.

\begin{theorem}
[Davidson-Fabec]\label{S:DF1}Let the notation be as above. Then the following hold:

\begin{enumerate}
\item \label{S:DF1(1)}$Q(W,Z)=J_{\lambda}(\exp(-\sigma(W)),Z)^{\ast-1}$.

\item \label{S:DF1(2)}$Q(W,Z)$ is holomorphic in the first variable and
antiholomorphic in the second variable.

\item \label{S:DF1(3)}$%
\ip{v}{Q(W,Z)u}%
=%
\ip{q_{W}v}{q_{Z}u}%
$ for all $u,v\in\mathbf{W}_{\lambda}$.

\item \label{S:DF1(4)}$Q$ is a positive-definite reproducing kernel.

\item \label{S:DF1(5)}$Q(g\cdot W,g\cdot Z)=J_{\lambda}(g,W)Q(W,Z)J_{\lambda
}(g,Z)^{\ast}$.
\end{enumerate}
\end{theorem}

For $Z\in{\Omega}_{\mathbb{C}}$ and $u\in\mathbf{W}^{\lambda}$, let
$F_{Z,u}\colon{\Omega}_{\mathbb{C}}\rightarrow\mathbf{W}^{\lambda}$ be the
holomorphic function
\begin{equation}
F_{Z,u}(W):=Q(W,Z)u \label{E:FQ}%
\end{equation}
and define
\begin{equation}%
\ip{F_{T,w}}{F_{Z,u}}%
_{Q}:=%
\ip{w}{Q(W,Z)u}%
. \label{E:Innpr}%
\end{equation}

Let $\mathbf{H}({\Omega}_{\mathbb{C}},\mathbf{W}^{\lambda})$ be the completion
of the span of $\{F_{Z,u}\mid Z\in{\Omega}_{\mathbb{C}},u\in\mathbf{W}%
^{\lambda}\}$ with respect to this inner product. Then $\mathbf{H}({\Omega
}_{\mathbb{C}},\mathbf{W}^{\lambda})$ is a Hilbert space consisting of
$\mathbf{W}^{\lambda}$-valued holomorphic functions. Define a representation
of $G^{c}$ in $\mathbf{H}({\Omega}_{\mathbb{C}},\mathbf{W}^{\lambda})$ by
\begin{equation}
(\rho(g)F)(W):=J_{\lambda}(g^{-1},W)^{-1}F(g^{-1}\cdot W)\,. \label{E:intrep}%
\end{equation}

Then $\rho$ is a unitary representation of $G^{c}$ in $\mathbf{H}({\Omega
}_{\mathbb{C}},\mathbf{W}^{\lambda})$ called the \textit{geometric
realization}\label{geometric} of $(\pi,\mathbf{V})$.

\begin{theorem}
[Davidson-Fabec]\label{ThmS-hwmNew.8}The map $q_{Z}v\mapsto F_{Z,v}$ extends
to a unitary intertwining operator $U$ between $(\pi,\mathbf{V})$ and
$(\rho,\mathbf{H}({\Omega}_{\mathbb{C}},\mathbf{W}^{\lambda}))$. It can be
defined globally by%
\[
\lbrack Uw](Z)=q_{Z}^{\ast}w,\quad w\in\mathbf{V},\,Z\in{\Omega}_{\mathbb{C}%
}\,.
\]
\end{theorem}

As the theorem stands, it gives a geometric realization for every unitary
highest-weight module. What it does not do is give a natural analytic
construction of the inner product on $\mathbf{H}({\Omega}_{\mathbb{C}%
},\mathbf{W}^{\lambda})$. This is known only for some special representations,
e.g., the \textit{holomorphic discrete series} of the group $G^{c}$
\cite{HCVI,DS,TI}
or
symmetric spaces of Hermitian type \cite{'OO88a,'OO88b}. At this point we will
only discuss the holomorphic discrete series, which was constructed by
Harish-Chandra in \cite{HCVI}, in particular Theorem 4 and Lemma 29. For that,
let $\rho=\frac{1}{2}\sum_{\alpha\in\Delta^{+}}\alpha$ and let $\mu$ denote
the highest weight of the representation of $K^{c}$ on $\mathbf{W}^{\lambda}$.
For $f,g\in\mathbf{H}({\Omega}_{\mathbb{C}},\mathbf{W}^{\lambda})$, let $\mu$
be the $G^{c}$-invariant measure on $\mathbb{\Omega}_{\mathbb{C}}$ and
\[%
\ip{g}{f}%
_{\lambda}:=\int_{G^{c}/K^{c}}%
\ip{g(Z)}{Q(Z,Z)^{-1}f(Z)}%
_{\mathbf{W}^{\lambda}}\,d\mu\,.
\]

\begin{theorem}
[Harish-Chandra \cite{HC65,HC70}]\label{s-HC} Assume that
\[
\left\langle \mu + \rho ,H_{\alpha}\right\rangle
<0\text{\qquad for all }\alpha\in\Delta_{p}^{+}.
\]
Then $%
\ip{g}{f}%
_{\lambda}$ is finite for $f,g\in\mathbf{H}({\Omega}_{\mathbb{C}}%
,\mathbf{W}^{\lambda})$ and there exists a positive constant $c_{\lambda}$
such that
\[%
\ip{g}{f}%
_{Q}=c_{\lambda}%
\ip{g}{f}%
_{\lambda}\,.
\]
Moreover, $(\rho,\mathbf{H}({\Omega}_{\mathbb{C}},\mathbf{W}_{\lambda}))$ is
unitarily equivalent to a discrete summand in $\mathbf{L}^{2}(G^{c})$.
\end{theorem}

The Theorem of Harish-Chandra relates some of the unitary highest weight
modules to the discrete part of the Plancherel measure. It was shown by
Ol'shanskii \cite{Ol82} and Stanton \cite{Stanton} that this ``holomorphic''
part of the discrete spectrum can
be
realized as a Hardy space of holomorphic
functions on a local tube domain. Those results were generalized to symmetric
spaces of \emph{Hermitian type} (or \emph{compactly causal symmetric spaces})
in a series of papers \cite{'OO88a,'OO88b,HOO91,OO96b,BO98}

The last theorem shows in particular that the corresponding highest weight
modules are unitary. It was shown by Wallach \cite{WaI72} and Rossi and Vergne
\cite{VR76} that those are not all the unitary highest weight modules. The
problem is to decide for which representations of $K^{c}$ the reproducing
kernel $Q(Z,W)$ is positive definite. We refer to \cite{EHW83,Ja83} for the
classification of unitary highest weight modules. We will from now on assume
that the representation of $K^{c}$ is a character $\chi_{\lambda}$ where
$\lambda\in i\frak{t}^{\ast}$ is trivial on $\frak{t}$ $\cap\lbrack
\frak{k},\frak{k}]$. Choose a maximal set $\left\{  \gamma_{1},\dots
,\gamma_{r}\right\}  $ of long strongly orthogonal roots in $\Delta_{p}^{+}$.
This can be done by putting $r=\operatorname*{rank}(G^{c}/K^{c})$ and then
choosing $\gamma_{r}$ to be a maximal root in $\Delta_{p}^{+}$, $\gamma_{r-1}$
maximal in $\{\gamma\in\Delta_{p}^{+}\mid\gamma\,\ $strongly orthogonal to
$\gamma_{r}\}$, etc. Let $H_{j}:=H_{\gamma_{j}}$ and
\begin{equation}
\frak{a}=i\bigoplus
_{j}\mathbb{R}H_{j}\subset\frak{t}. \label{eqfraktura}
\end{equation}
By the theorem of Moore (see\ \cite{He78}%
) we know that the roots in $\Delta_{p}$ restricted to $\frak{a}$, are given
by $\pm\frac{1}{2}(\gamma_{i}+\gamma_{j})$, $1\leq i\leq j\leq r$ and possibly
$\frac{1}{2}\gamma_{j}$. The root spaces for $\gamma_{j}$ are all
one-dimensional and the root spaces $\frak{g}_{\pm\frac{1}{2}(\gamma
_{i}+\gamma_{j})}$, $1\leq i<j\leq r$, have all the common dimension $d$.

\begin{theorem}
[Vergne-Rossi \cite{VR76}, Wallach \cite{Wal92}]\label{S:VRW} Assume that
$G^{c}$ is simple. Let $\lambda_{0}\in\frak{a}^{\ast}$ be such that
$\left\langle \lambda_{0},H_{r}\right\rangle =1$. Let $\gamma=\left\langle
\lambda_{0},Z^{0}\right\rangle $ and let
\[
L_{\mathrm{pos}}:=-\frac{\gamma(r-1)d}{2}\,.
\]
For $\nu-\rho<L_{\mathrm{pos}}$ there exists a irreducible unitary highest
weight representation $(\rho_{\nu},\mathbf{K}_{\nu})$ of $G^{c}$ with
one-dimensional minimal $K^{c}$-type $\nu-\rho$.
\end{theorem}

\begin{proof}
By \cite[pp.~41--42]{VR76} (see also \cite{NW79}), $(\rho_{\nu},\mathbf{K}%
_{\lambda})$ exists if $\left\langle \nu-\rho,H_{r}\right\rangle \leq
-\frac{(r-1)d}{2}$. But $\nu-\rho=\left\langle \nu-\rho,H_{r}\right\rangle
\lambda_{0}$. Hence $\left\langle \nu-\rho,2Z^{0}\right\rangle =\left\langle
\nu-\rho,H_{r}\right\rangle \left\langle \lambda_{0},2Z^{0}\right\rangle
=\gamma\left\langle \nu-\rho,H_{r}\right\rangle $.
\end{proof}

We will later specialize this to the case where $G^{c}/K^{c}$ is a tube type
domain, that is biholomorphically equivalent to $\mathbb{R}^{n}+i\Omega$, where
$n=\dim_{\,\mathbb{C}}G^{c}/K^{c}$, and $\Omega$ is a self-dual regular cone in
$\mathbb{R}^{n}$. If we assume $G^{c}$ simple, then this is exactly the case
if $G^{c}$ is locally isomorphic to one of the groups: $SU(n,n)$, $SO^{\ast
}(4n)$, $Sp(n,\mathbb{R)}$, $SO_{o}(n,2)$ and $E_{7(-25)}$. In this case we
have%
\begin{equation}
Z^{0}=\frac{i}{2}\sum_{j=1}^{r}H_{j} \label{eqS-hwmNew.6}%
\end{equation}
and%
\begin{align}
\Delta_{p}^{+}  &  =\left\{  \gamma_{i},\frac{1}{2}(\gamma_{k}+\gamma
_{j})\biggm| 1\leq i,j,k\leq r,\,\,j\leq k\right\} \label{eqS-hwmNew.7}\\
\Delta_{c}^{+}  &  =\left\{  \frac{1}{2}(\gamma_{k}-\gamma_{j})\biggm|
1\leq\,j\leq k\leq r\right\}  \,. \label{eqS-hwmNew.8}%
\end{align}
In this case we have:

\begin{lemma}
\label{CThw} Assume that $G^{c}/K$ is a tube-type domain, then $\gamma=r$ and
$\nu-\rho\leq L_{\mathrm{pos}}$ if and only if $\nu\leq r$.
\end{lemma}

\begin{proof}
If $G^{c}/K^{c}$ is of Cayley type then $2Z^{0}=i\sum_{j=1}^{r}H_{j}$ and
$\gamma_{j}=\gamma_{r}-\sum n_{\alpha}\alpha$, $\alpha\in\Delta_{c}^{+}$,
$n_{\alpha}\geq0$. Thus $\left\langle \nu-\rho,Z^{0}\right\rangle
=r\left\langle \nu-\rho,H_{r}\right\rangle $. We also have (see\ \cite{OO96b}%
)
\[
\rho=\frac{1}{2}\left(  1+\frac{(r-1)d}{2}\right)  (\gamma_{1}+\cdots
+\gamma_{r})\,.
\]
{}From this the theorem follows.
\end{proof}

\begin{remark}
\label{RemS-hwmNew.13}Let us remind the reader that we have only described
here the continuous part of the unitary spectrum. There are also finitely many
discrete points, the so-called \textit{Wallach set}, giving rise to unitary
highest weight representations.
\end{remark}

\begin{remark}
\label{R:oppo}Let $\sigma\colon\frak{g}_{\mathbb{C}}^{c}\rightarrow
\frak{g}_{\mathbb{C}}^{c}$ be the conjugation with respect to $\frak{g}^{c}$.
Thus $\sigma(X+iY)=X-iY$, $X,Y\in\frak{g}^{c}$. Then $\sigma(\frak{p}%
^{+})=\frak{p}^{-}$ and $\sigma(\frak{k}_{\mathbb{C}}^{c})=\frak{k}%
_{\mathbb{C}}^{c}$. In this section we viewed $Q(W,Z)$ as a function on
$\Omega_{\mathbb{C}}\times\Omega_{\mathbb{C}}$, holomorphic in the first
variable and antiholomorphic in the second variable. In many applications it
is better to view $Q$ as a function on $\Omega_{\mathbb{C}}\times\sigma
(\Omega_{\mathbb{C}})$, \ holomorphic in both variables.
\end{remark}

\section{\label{S:Ex}An Example: $SU(1,1)$}

The simplest case of a non-trivial reflection positivity is the case
$G=SL(2,\mathbb{R})$ and $G^{c}=SU(1,1)$. In this case
\[
G^{c}/K^{c}=D:=\left\{  z\in\mathbb{C}\mid\left|  z\right|  <1\right\}
\]
and $SU(1,1)$ acts by
\[
\left(
\begin{array}
[c]{cc}%
a & b\\
\bar{b} & \overline{a}%
\end{array}
\right)  \cdot z=\frac{az+b}{\bar{b}z+\overline{a}}\,.
\]
Let $\sigma\colon D\rightarrow D$ be complex conjugation, $z\mapsto\bar{z}$
and let $\tau\colon SL(2,\mathbb{C})\rightarrow SL(2,\mathbb{C)}$ be the
involution given by%
\begin{equation}
\tau\left(  \left(
\begin{array}
[c]{cc}%
a & b\\
c & d
\end{array}
\right)  \right)  =\left(
\begin{array}
[c]{cc}%
0 & 1\\
1 & 0
\end{array}
\right)  \left(
\begin{array}
[c]{cc}%
a & b\\
c & d
\end{array}
\right)  \left(
\begin{array}
[c]{cc}%
0 & 1\\
1 & 0
\end{array}
\right)  =\left(
\begin{array}
[c]{cc}%
d & c\\
b & a
\end{array}
\right)  \,. \label{suinvo}%
\end{equation}
Then $\sigma(g\cdot0)=\tau(g)\cdot0$ for $g\in SU(1,1,)$. We have%
\[
H_{\mathbb{C}}=\left\{  \left(
\begin{array}
[c]{cc}%
z & w\\
w & z
\end{array}
\right)  \biggm|z^{2}-w^{2}=1\right\}  \,.
\]
and%
\begin{equation}
H=\pm\left\{  h_{t}=\left(
\begin{array}
[c]{cc}%
\cosh(t) & \sinh(t)\\
\sinh(t) & \cosh(t)
\end{array}
\right)  \biggm|t\in\mathbb{R}\right\}  =SU(1,1)\cap SL(2,\mathbb{R})\,,\,
\label{HSU}%
\end{equation}
$H/\left\{  \pm I\right\}  =(-1,1)$ and $G=SL(2,\mathbb{R})$. Knowing that the
representations of $\widetilde{SU(1,1)}$ that we can get are highest weight
modules, we see by looking at the infinitesimal character of those
representations, that we have to start with the complementary series
representations of $G=SL(2,\mathbb{R})$. They are constructed in the following
way. Let $P$ be the parabolic subgroup
\[
P:=\left\{  p(a,x)=\left(
\begin{matrix}
a & x\\
0 & a^{-1}%
\end{matrix}
\right)  \biggm|a\in\mathbb{R}^{\ast},\,x\in\mathbb{R}\right\}  =G\cap
K_{\mathbb{C}}^{c}P^{+}\,.
\]
For $s\in\mathbb{C}$, let $\pi_{s}$ be the representation of $G$ acting by
$[\pi_{s}(a)f](b)=f(a^{-1}b)$ on the space $\mathbf{H}_{s}$ of functions
$f\colon G\rightarrow\mathbb{C}$,
\[
f(gp(a,x))=|a|^{-(s+1)}f(g)\,,\quad\int_{SO(2)}|f(k)|^{2}\,dk<\infty\,,
\]
and with inner product
\[%
\ip{f}{g}%
=\int_{SO(2)}\overline{f(k)}g(k)\,dk\,,
\]
that is $\pi_{s}$ is the principal series representation of $G$ with parameter
$s$. A simple calculation shows that the pairing%
\begin{equation}
\mathbf{H}_{s}\times\mathbf{H}_{-s}\ni(f,g)\longmapsto\int_{SO(2)}%
f(k)g(k)\,dk\in\mathbb{C} \label{paring}%
\end{equation}
is invariant under the group action. The representations $\pi_{s}$ are unitary
in the above Hilbert-space structure as long as $s\in i\mathbb{R}$. Let%
\[
\bar{N}=\left\{  \bar{n}_{t}=\left(
\begin{array}
[c]{cc}%
1 & 0\\
t & 1
\end{array}
\right)  \biggm|t\in\mathbb{R}\right\}  \,.
\]
Then
\[
\bar{n}_{t}p(\gamma,x)=\left(
\begin{array}
[c]{cc}%
\gamma &  x\\
\gamma t & \gamma^{-1}+xt
\end{array}
\right)  =\left(
\begin{array}
[c]{cc}%
a & b\\
c & d
\end{array}
\right)  \,.
\]
Hence as we have seen before $\bar{N}P=\left\{  \left(
\begin{array}
[c]{cc}%
a & b\\
c & d
\end{array}
\right)  \biggm|a\not =0\right\}  $ and%
\begin{equation}
\gamma=a,\,\,t=c/a\,. \label{E:proje}%
\end{equation}
In particular we have%
\[
\left(
\begin{array}
[c]{cc}%
\cosh(t) & \sinh(t)\\
\sinh(t) & \cosh(t)
\end{array}
\right)  =\left(
\begin{array}
[c]{cc}%
1 & 0\\
\tanh(t) & 1
\end{array}
\right)  \left(
\begin{array}
[c]{cc}%
\cosh(t) & 0\\
0 & 1/\cosh(t)
\end{array}
\right)  \left(
\begin{array}
[c]{cc}%
1 & \tanh(t)\\
0 & 1
\end{array}
\right)  \,.
\]
Thus $HP/P\simeq(-1,1)$, but we notice that this is not the realization in
$\frak{p}^{+}$ but in $\frak{p}^{-}$. By (\ref{E:proje}) this is expressed in
the action of $G$ by%
\begin{equation}
\left(
\begin{array}
[c]{cc}%
a & b\\
c & d
\end{array}
\right)  \cdot_{\text{$\operatorname*{opp}$}}z=\frac{dz+c}{bz+a}\,.
\label{E:newact}%
\end{equation}
Notice that this is the usual action twisted by $\tau$, that is
\begin{equation}
g\cdot_{\text{$\operatorname*{opp}$}}z=\tau(g)\cdot z, \label{E:tautwist}%
\end{equation}
where $\cdot$ stands for the usual action. By identify $\bar{N}$ with
$\mathbb{R}$ using $\bar{n}_{t}\mapsto t$, we can realize the principal series
representations as acting on functions on $\mathbb{R}$ (compare to
(\ref{E:newact})):%
\begin{align}
\pi_{s}\left(  \left(
\begin{array}
[c]{cc}%
a & b\\
c & d
\end{array}
\right)  \right)  f(t)  &  =\left|  d-bt\right|  ^{-s-1}f\left(  \frac
{-c+at}{d-bt}\right)  \,\label{eqEx.21}\\%
\ip{f}{g}%
&  =\frac{1}{\pi}\int_{-\infty}^{\infty}\overline{f(t)}\,g(t)(1+t^{2}%
)^{1+\operatorname{Re}(s)}\,dt\, \label{eqEx.22}%
\end{align}
and the pairing in (\ref{paring}) is simply%
\begin{equation}%
\ip{f}{g}%
=\frac{1}{\pi}\int_{-\infty}^{\infty}\overline{f(t)}\,g(t)\,dt\,.
\label{pairing2}%
\end{equation}
For defining the complementary series we need the intertwining operator
$A_{s}\colon\mathbf{H}_{s}\rightarrow\mathbf{H}_{-s}$ defined by
\begin{equation}
A_{s}(f)(g)=\frac{1}{\pi}\int_{-\infty}^{\infty}f(gw\bar{n}_{y})\,dy
\label{eqEx.24}%
\end{equation}
for $\operatorname{Re}s\geq0$ and then generally by analytic continuation.
Here $w$ is the Weyl group element $w=\left(
\begin{matrix}
0 & 1\\
-1 & 0
\end{matrix}
\right)  $ and $\bar{n}_{y}=\left(
\begin{matrix}
1 & 0\\
y & 1
\end{matrix}
\right)  $. In our realization of the representation on $\mathbb{R}$ we get
\begin{equation}
A_{s}f(x)=\frac{1}{\pi}\int_{-\infty}^{\infty}f(y)|x-y|^{s-1}\,dy\,.
\label{eqEx.25}%
\end{equation}
By (\ref{pairing2}) the bilinear form
\begin{equation}%
\ip{f}{A_{s}g}%
=\frac{1}{\pi^{2}}\int_{-\infty}^{\infty}\int_{-\infty}^{\infty}%
\overline{f(x)}g(y)|x-y|^{s-1}\,dx\,dy \label{eqEx.26}%
\end{equation}
is $G$-invariant and actually the representation $\pi_{s}$ is unitary for
$0<s<1$. Define%
\[
Jf(t)=\left|  t\right|  ^{-s-1}f(1/t)
\]
or on the group level $Jf(a):=f(\tau(a)w^{-1})=f(\tau(aw))$. The map
$J\colon\mathbf{H}_{s}\rightarrow\mathbf{H}_{s}$ intertwines $\pi_{s}$ and
$\pi_{s}\circ\tau\simeq\pi_{s}$, $J^{2}=1$, and
\[
A_{s}(J(g))(x)=\frac{1}{\pi^{2}}\int_{-\infty}^{\infty}g(y)|1-xy|^{s-1}%
\,dy\,.
\]
Hence
\begin{equation}%
\ip{f}{g}%
_{J}=%
\ip{f}{A_{s}Jg}%
=\frac{1}{\pi^{2}}\int_{-\infty}^{\infty}\int_{-\infty}^{\infty}%
\overline{f(x)}g(y)|1-xy|^{s-1}\,dy\,dx\,. \label{eqEx.27}%
\end{equation}
We recall now the reproducing kernel $Q(w,z)$ from Section \ref{S-hwm},
corresponding to the lowest $K^{c}$-type $s-1$. In our case%
\begin{equation}
Q(w,z)=(1-w\bar{z})^{s-1}\,. \label{eqEx.28}%
\end{equation}
which by Theorem \ref{CThw} is positive if and only if $s\leq1$. It follows
that $%
\ip{\,\cdot\,}{\,\cdot\,}%
_{J}$ is positive definite on the space of functions supported on the
$H$-orbit $(-1,1)$. Let $\mathbf{K}_{0}$ be the closure of $\mathcal{C}%
_{c}^{\infty}(-1,1)$. Notice that the above inner product is defined on
$\mathcal{C}_{c}^{\infty}(-1,1)$ for every $s$ as we only integrate over
compact subsets of $(-1,1)$. As we are using the realization in $\frak{p}^{-}$
we define the semigroup now by
\begin{equation}
S=S_{-}(\Omega):=\left\{  \gamma\in SL(2,\mathbb{R})\mid\gamma\cdot
(-1,1)\subset(-1,1)\right\}  \,. \label{eqEx.29}%
\end{equation}
Then $S$ is a closed semigroup of the form $H=\exp(C)$ where $C$ is the
$H$-invariant cone generated by $\left(
\begin{array}
[c]{cc}%
-1 & 0\\
0 & 1
\end{array}
\right)  $. Let us remark here, that if $\gamma=\left(
\begin{array}
[c]{cc}%
a & b\\
c & d
\end{array}
\right)  \in S$, then $d-bt>0$ for $\left|  t\right|  <1$. The semigroup $S$
acts on $\mathbf{K}_{0}$ and by the L\"uscher-Mack Theorem \ref{LM} we get
an highest weight module for $\widetilde{SU(1,1)}$, which is irreducible as we
will see in a moment.

We also know (see\ \cite{HO95}) that $S=H\exp C$ is a closed semigroup and
that $\gamma I\subset I$, and actually $S$ is exactly the semigroup of
elements in $SL(2,\mathbb{R})$ that act by contractions on $I$. Hence $S$ acts
on $\mathbf{K}$. By a theorem of L\"uscher and Mack \cite{HiNe93,LM75}, the
representation of $S$ on $\mathbf{K}$ extends to a representation of $G^{c}$,
which in this case is the universal covering of $SU(1,1)$ that is locally
isomorphic to $SL(2,\mathbb{R})$. According to Theorem \ref{HighestWeight} the
resulting representation is a direct integral of highest weight
representations. We notice that this defines a representation of
$SL(2,\mathbb{R})$ if and only if certain integrality conditions hold; see
\cite{JoMo84}. The question then arises to identify this direct integral and
construct an explicit intertwining operator into the corresponding space of
holomorphic functions on $D$.

We notice first that the kernel $(y,x)\mapsto Q(y,x)=(1-yx)^{s-1}$ is the
reproducing kernel of the irreducible highest weight representation given by%
\begin{equation}
\lbrack\rho_{s}(g)f](z)=(-cz+a)^{s-1}f\left(  \frac{dz-b}{-cz+a}\right)
\,,\quad g=\left(
\begin{array}
[c]{cc}%
a & b\\
c & d
\end{array}
\right)  \,. \label{eqEx.30}%
\end{equation}
In particular $Q(y,x)$ extends to a holomorphic function on $D\times\sigma(D)$
according to Remark \ref{R:oppo}. For $f\in\mathcal{C}_{c}^{\infty}(-1,1)$
define%
\begin{equation}
Uf(z):=\frac{1}{\pi}\int_{-1}^{1}f(u)(1-zu)^{s-1}\,du=\frac{1}{\pi}\int
_{-1}^{1}f(u)Q(z,u)\,du\,. \label{eqEx.31}%
\end{equation}
By simple calculation, using (\ref{E:tautwist}), we get for $\gamma=\left(
\begin{array}
[c]{cc}%
a & b\\
c & d
\end{array}
\right)  \in S_{-}(\Omega)$:%
\begin{align*}
U(\pi_{s}(\gamma)f)(z)  &  =\frac{1}{\pi}\int_{-1}^{1}(d-bt)^{-s-1}f\left(
\frac{-c+at}{at-c}\right)  (1-zt)^{s-1}\,dt\\
&  =\frac{1}{\pi}\int_{-1}^{1}f(u)(cu+d)^{s-1}\left(  1-z\frac{du+c}%
{bu+a}\right)  ^{s-1}\,du\\
&  =\frac{1}{\pi}\int_{-1}^{1}f(u)\left(  bu+a-dzu-cz\right)  ^{s-1}\,du\\
&  =\frac{1}{\pi}\int_{-1}^{1}f(u)\left(  -cz+a-(dz-b)u\right)  ^{s-1}\,du\\
&  =(-cz+a)^{s-1}\frac{1}{\pi}\int_{-1}^{1}f(u)\left(  1-\frac{dz-b}%
{-cz+a}\right)  ^{s-1}\,du\\
&  =\rho_{s}(\gamma)Uf(z)\,,
\end{align*}
where the respective representations are given by (\ref{eqEx.21}) and
(\ref{eqEx.30}). Here the last equality follows from (\ref{E:jsu}),
(\ref{E:jlambda}) and (\ref{E:intrep}). As $\rho_{s}$ is irreducible it
follows that either $U$ is surjective or identically zero. Using that $Q(z,u)$
is the reproducing kernel for the representation $\rho_{s}$ we get for $f$ and
$g$ with compact support:%
\begin{align*}%
\ip{Uf}{Ug}%
&  =\frac{1}{\pi^{2}}\int_{-1}^{1}\int_{-1}^{1}\overline{f(u)}g(v)\,%
\ip{Q(\cdot,u)}{Q(\cdot,v)}%
\,dv\,du\\
&  =\frac{1}{\pi^{2}}\int_{-1}^{1}\int_{-1}^{1}\overline{f(u)}%
g(v)Q(u,v)\,dv\,du\\
&  =%
\ip{f}{g}%
\,.
\end{align*}
It follows that $U$ is a unitary isomorphism.

We can describe $U$ in a different way using the representation $\rho_{s}$
instead of the reproducing kernel. Let $\openone$ be the constant function
$z\mapsto1$. Then%
\begin{equation}
\lbrack\rho_{s}(g)\openone](z)=J_{s}(g^{-1},z)=(-cz+a)^{s-1}\,.
\label{eqEx.32}%
\end{equation}
We therefore get%
\begin{align}
\int_{H}f(h)\rho_{s}(h)\openone(z)\,dh  &  =\frac{1}{\pi}\int_{-\infty
}^{\infty}\left[  \cosh(t)^{-s-1}f(\tanh(t))\right]  (-\sinh(t)z+\cosh
(t))^{s-1}\,dt\label{E:R*one}\\
&  =\frac{1}{\pi}\int_{-1}^{1}f(u)(1-uz)^{s-1}\,du\nonumber\\
&  =Uf(z)\,.\nonumber
\end{align}
We will meet the transform in (\ref{E:R*one}) again in the generalization of
the Bargmann transform in Section \ref{S:Bargmann}. That shows that the
Bargmann transform introduced in \cite{OO96} is closely related to the
reflection positivity and the Osterwalder-Schrader duality.

In summary, we have the representation $\pi_{s}$ from (\ref{eqEx.21}) acting
on the Hilbert space $\widehat{\mathbf{H}_{+}(s)}$ of distributions obtained
from completion with respect to%
\[
\int_{-1}^{1}\int_{-1}^{1}\overline{f(x)}\,f(y)\left(  1-xy\right)
^{s-1}\,dx\,dy\,,
\]
and the unitarily equivalent representation $\rho_{s}$ from (\ref{eqEx.30}).
The operator $U$ from (\ref{eqEx.31}) intertwines the two. Moreover $U$ passes
to the distributions on $\left(  -1,1\right)  $, in the completion
$\widehat{\mathbf{H}_{+}(s)}$, and we have%
\begin{equation}
U\left(  \delta^{\left(  n\right)  }\right)  =\frac{\left(  s-1\right)
\left(  s-2\right)  \cdots\left(  s-n\right)  }{\pi}\,z^{n},\label{eqEx.34}%
\end{equation}
where, for $n=0,1,2,\dots$, $\delta^{\left(  n\right)  }=\left(  d/dx\right)
^{n}\delta$ are the derivatives of the Dirac ``function'',
defined by
\begin{equation}
\left\langle f,\delta^{\left(  n\right)  }\right\rangle
=\left\langle \left( -1\right)^{n}f^{\left(  n\right)  },\delta\right\rangle
=\left( -1\right)^{n}f^{\left(  n\right)  }\left( 0\right) ,\label{eqExApr.35}
\end{equation}
where $f$ is a test function. Furthermore,
$z^{n}$ are the
monomials in the reproducing kernel Hilbert space $\mathbf{H}(s)$
corresponding to the complex kernel $\left(  1-\bar{z}w\right)  ^{s-1}$. This
Hilbert space consists of analytic functions $f\left(  z\right)  =\sum
_{n=0}^{\infty}C_{n}z^{n}$, defined in $D=\left\{  z\in\mathbb{C}\mid\left|
z\right|  <1\right\}  $, and satisfying%
\[
\sum_{n=0}^{\infty}\left|  C_{n}\right|  ^{2}\frac{1}{\left|  \binom{s-1}%
{n}\right|  }<\infty,
\]
where the $\binom{s-1}{n}$ refers to the (fractional) binomial coefficients.
This sum also defines the norm in $\mathbf{H}(s)$. For every $w\in D$, the
function $u_{w}\left(  z\right)  :=\left(  1-\bar{w}z\right)  ^{s-1}$ is in
$\mathbf{H}(s)$, and for the inner product, we have%
\[%
\ip{u_{w_{1}}}{u_{w_{2}}}%
_{\mathbf{H}(s)}=\left(  1-\bar{w}_{1}w_{2}\right)  ^{s-1}\,.
\]
So $\mathbf{H}(s)$ is indeed a reproducing kernel Hilbert space, as it follows
that the values $f\left(  z\right)  $, for $f\in\mathbf{H}(s)$ and $w\in D$,
are given by the inner products%
\[
f\left(  w\right)  =%
\ip{u_{w}}{f}%
_{\mathbf{H}(s)}\,.
\]
Since $u_{w}\left(  z\right)  =\sum_{n=0}^{\infty}\binom{s-1}{n}\bar{w}%
^{n}z^{n}$, we conclude that the monomials $z^{n}$ form an orthogonal basis in
$\mathbf{H}(s)$, and it follows from (\ref{eqEx.34}) that the distributions
$\delta^{\left(  n\right)  }$, $n=0,1,2,\dots$, form an orthogonal basis in
the Hilbert space $\widehat{\mathbf{H}_{+}(s)}$, and that%
\[
\left\|  \delta^{\left(  n\right)  }\right\|  _{\widehat{\mathbf{H}_{+}(s)}%
}^{2}=
\frac{n!\left(  1-s\right)  \left(  2-s\right)  \cdots\left(  n-s\right)  }
{\pi^{2\mathstrut}}.
\]

\section{\label{S:Sssp}Reflection Symmetry for Semisimple Symmetric Spaces}

The main results in this section are Theorems \ref{PR} and \ref{S:Posref}.
They are stated for non-compactly causal symmetric spaces, and the proofs are
based on our Basic Lemma and the L\"uscher-Mack theorem. At the end of the
section we show that results from Jorgensen's paper \cite{Jor86} lead to an
extension of the scope of the two theorems.

We now generalize the construction from the last section to a bigger class of
semisimple symmetric pairs. We restrict ourself to the case of characters
induced from a maximal parabolic subgroup, which leads to highest weight
modules with one-dimensional lowest $K^{c}$-type. This is meant as a
simplification and not as a limitation of our method. An additional source of
inspiration for the present chapter is the following series of papers:
\cite{Nel59,'O90a,'OO89a,'OO88a,'OO88b,OO96,JO97,KlLa83,OsSc73,Pra89,Sch86}.

Assume that $G^{c}/K^{c}%
=D\simeq\Omega_{\mathbb{C}}\subset\frak{p}^{+}$ is a bounded symmetric domain
with $G^{c}$ simply connected and simple. Let $\theta^{c}$ be the Cartan
involution on $G^{c}$ corresponding to $K^{c}$.
Let $\sigma\colon D\rightarrow D$ be a conjugation, that is a non-trivial
order two antiholomorphic map. Those involutions were classified in
\cite{HJ75,HJ78}, see also \cite{HO95,'O90b,'O90a}. Then $\sigma$ defines an
involution on the group $I_{o}(D)$, the connected component of holomorphic
isometries of $D$, by%
\[
\tau(f)(Z)=\sigma\left(  f(\sigma(Z))\right)  \,.
\]
But $I_{o}(D)$ is locally isomorphic to $G^{c}$, see \cite{He78}, Chapter
VIII. Hence $\tau$ defines an involution on $G^{c}$ and $\frak{g}^{c}$. Let
$H^{c}=G^{c\tau}$, and $\frak{h=}\left\{  X\in\frak{g}^{c}\mid\tau
(X)=X\right\}  $ and $\frak{q}^{c}=\left\{  X\in\frak{g}^{c}\mid
\tau(X)=-X\right\}  $. Then $\frak{g}^{c}=\frak{h}\oplus\frak{q}^{c}$. We
define%
\[
\frak{g}:=\frak{h}\oplus i\frak{q}^{c}%
\]
and $\frak{q=}\left\{  X\in\frak{g}\mid\tau(X)=-X\right\}  =i\frak{q}^{c}$.
Then $(\frak{g},\tau)$ is a symmetric pair. Let $G_{\mathbb{C}}$ be a simply
connected Lie group with Lie algebra $\frak{g}_{\mathbb{C}}$ and let $G\subset
G_{\mathbb{C}}$ be the connect Lie group with Lie algebra $\frak{g}$ . Then
$\tau$ integrates to an involution on $G$. Let $H=G^{\tau}=\left\{  a\in
G\mid\tau(a)=a\right\}  $. Then $G/H$ is a symmetric space. The involution
$\theta^{c}$ integrates to an involution on $G$ and $\theta
:=\tau\theta^{c}$ is a Cartan involution on $G$ that commutes with $\tau$. Let
$K$ be the corresponding maximal almost compact subgroup. Denote the
corresponding Cartan decomposition as usually by $\frak{g}=\frak{k}%
\oplus\frak{p}$.

As $\tau$ is antiholomorphic it follows that $\tau(Z^{0})=-Z^{0}$, where
$Z^{0}$ is a central element in $\frak{k}^{c}$ with eigenvalues $0,i,-i$%
.\ Hence $G^{c}/H^{c}$ is a symmetric space of \emph{Hermitian type}, in the
sense of \cite{'OO88a}. Those spaces are now usually called \emph{compactly
causal symmetric spaces} because those are exactly the symmetric spaces such
that $\frak{q}$ contains a regular $H$-invariant cone $C$ with $C^{o}%
\cap\frak{k\not =\emptyset}$. The minimal cone is given by%
\[
C_{\min}^{c}=\mathbb{R}^{+}\cdot\overline{\operatorname{conv}\left\{
\operatorname{Ad}(H)Z^{0}\right\}  }\,.
\]
The dual spaces $G/H$ are exactly the \emph{non-compactly causal symmetric
spaces}. Those are the symmetric spaces containing $H$-invariant regular cones
$C$ such that $C^{o}\cap\frak{p}\not =\emptyset$. We use \cite{HO95} as a
standard reference to the causal symmetric spaces.

\begin{example}
[Cayley type spaces]\label{ExaSsspNew.1}A special case of the above
construction is when $G^{c}/K^{c}$ is a tube type domain. Let $\mathbf{c}$ be
a Cayley transform from the bounded realization of $G^{c}/K^{c}$ to the
unbounded realization. This can be done by choosing $\mathbf{c}%
=\operatorname{Ad}(\exp(\frac{\pi i}{2}Y^{0}))$ where $\operatorname{ad}%
(Y^{0})$ has eigenvalues $0,1,-1$. Then $\operatorname{Ad}(\mathbf{c}%
)^{4}=\operatorname*{id}$ and $\tau=\operatorname{Ad}(\mathbf{c)}^{2}%
(G^{c})=G^{c}$. Hence $\tau$ is an involution on $G^{c}$. It is also well
known that $\tau(Z^{0})=-Z^{0}$. Hence $\tau$ defines a conjugation on $D$.
The symmetric spaces $G/H$ are the symmetric spaces of Cayley type. We have
$Y^{0}\in\frak{h}$ is central and $\frak{z}_{\frak{g}}(Y^{0})=\frak{h}$.
Furthermore $\operatorname{Ad}(\mathbf{c})$ is an isomorphism $\frak{g}%
^{c}\simeq\frak{g}$. The spaces that we get from this construction are locally
isomorphic to one of the following symmetric spaces, where we denote by the
subscript $+$ the group of elements having positive determinant:
$Sp(n,\mathbb{R})/GL(n,\mathbb{R})_{+}$, $SU(n,n)/GL(n,\mathbb{C})_{+}$,
$SO^{\ast}(4n)/SU^{\ast}(2n)\mathbb{R}_{+}$, $SO(2,k)/SO(1,k-1)\mathbb{R}_{+}$
and $E_{7(-25)}/E_{6(-26)}\mathbb{R}_{+}$.
\end{example}

\begin{example}
\label{ExaSsspNew.2}Assume that $H$ is a connected Lie group such that
$H/K_{H}$, $K_{H}$ a maximal compact subgroup of $H$, is a bounded symmetric
domain. Let $G^{c}=H\times H$ and $D=H/H_{K}\times\overline{H/H_{K}}$, where
the bar denotes opposite complex structure. Let $\tau(d,c)=(c,d)$. Then $\tau$
is a conjugation with fixed-point set the diagonal. The corresponding
involution on $G^{c}$ is $\tau(a,b)=(b,a)$. Thus $G^{c\tau}%
=\operatorname*{diagonal}\simeq H$. Identify $G^{c\tau}$ with $H$. Then
$G/H\ni(a,b)H\mapsto ab^{-1}\in H$ is an isomorphism. In this case $G$ is
locally isomorphic to $H_{\mathbb{C}}$ and the involution $\tau$ on $\frak{g}$
is the conjugation with respect to the real form $\frak{h}\subset\frak{g}$.
Let $H_{1}$ be the corresponding analytic subgroup. Then $H_{1}$ is locally
isomorphic to $H$ and the symmetric space we are looking at is $G_{\mathbb{C}%
}/H_{1}$.
\end{example}

We will need the following facts. Let $X^{0}=-iZ^{0}\in\frak{q}\cap\frak{p}$.
Then $X^{0}$ is $H\cap K$-invariant,%
\begin{equation}
\frak{z}_{\frak{g}}(X^{0})=\frak{k}_{\mathbb{C}}^{c}\cap\frak{g=k}\cap
\frak{h}\oplus\frak{p}\cap\frak{q}, \label{E:ZXo}%
\end{equation}
Let $\frak{n}:=\frak{p}^{+}\cap\frak{g}$, $\frak{\bar{n}}:=\frak{p}^{-}%
\cap\frak{g}$, and $\frak{p}_{\max}:=(\frak{k}_{\mathbb{C}}^{c}\oplus
\frak{p}^{+})\cap\frak{g}$. Then $\frak{p}_{\max}$ is a maximal parabolic
subgroup of $\frak{g}$ of the form $\frak{p}_{\max}=\frak{m}\oplus
\mathbb{R}X^{0}\oplus\frak{n}$, where $\frak{m}=\left\{  X\in\frak{k}%
\cap\frak{h}\oplus\frak{p}\cap\frak{q}\mid B(X,X^{0})=0\right\}  $, $B$ the
Killing form on $\frak{g}$. We have $H\cap K=Z(X^{0})$. Let $A:=\exp
(\frak{a})$, $N:=\exp(\frak{n})$, $\bar{N}:=\exp(\frak{\bar{n}})$. $M_{0}$ the
analytic subgroup of $G$ corresponding to $\frak{m}$, and $M:=(H\cap K)M_{0}$.
Then $M$ is a closed and $\tau$-stable subgroup of $G$, $M\cap A=\left\{
1\right\}  $, $MA=Z_{G}(A)$, and $P_{\max}=N_{G}(\frak{p}_{\max})=MAN$. We
have
\[
\frak{g}=\frak{h}\oplus\frak{p}_{\max}\,.
\]
Let $\Omega=\tau(\Omega_{\mathbb{C}})\cap\frak{g\subset\bar{n}}$. Then by
\cite{FHO93}:

\begin{lemma}
\label{HPminOpenG}$HP_{\mathrm{min}}$ is open in $G$ and $HP_{\max}%
=\exp(\Omega)P_{\max}\subset$ $\bar{N}P_{\mathrm{max}}$.
\end{lemma}

Let $\frak{a}=\mathbb{R}X^{0}$ and $A:=\exp(\frak{a})$. We need to fix the
normalization of measures before we discuss the generalized principal series
representations. Let the measure $da$ on $A$ be given by
\[
\int_{A}f(a)\,da=\frac{1}{\sqrt{2\pi}}\int_{-\infty}^{\infty}\,f(a_{t}%
)\,dt,\quad a_{t}=\exp2tX^{0}\,.
\]
Then Fourier inversion holds without any additional constants. We fix the
Lebesgue measure $dX$ on $\bar{\frak{n}}$ such that, for $d\bar{n}=\exp(dX)$,
we then have
\[
\int_{\bar{\frak{n}}}a(\bar{n})^{-2\rho}\,d\bar{n}=1\,.
\]
Here $\rho(X)=\frac{1}{2}\operatorname*{tr}(${$\operatorname{ad}$%
}$(X))|_{\frak{n}}$ as usual, and $a(g)\in A$, $g\in G$, is determined by
$g\in KMa(g)N$. The Haar measure on compact groups will usually be normalized
to have total measure one. The measure on $N$ is $\theta(d\bar{n})$. We fix a
Haar measure $dm$ on $M$, and $dg$ on $G$ such that
\[
\int_{G}f(g)\,dg=\int_{K}\int_{M}\int_{A}\int_{N}f(kman)a^{2\rho
}\,dn\,da\,dm\,dk\,,\qquad f\in\mathcal{C}_{c}^{\infty}(G).
\]
Then we can normalize the Haar measure $dh$ on $H$ such that for
$f\in\mathcal{C}_{c}^{\infty}(G)$, $\operatorname*{supp}(f)\subset
HP_{\mathrm{max}}$, we have, see \cite{'O87}:
\[
\int_{G}f(g)\,dg=\int_{H}\int_{M}\int_{A}\int_{N}f(hman)a^{2\rho
}\,dn\,da\,dm\,dh\,.
\]
The invariant measure $d\dot{x}$ on $G/H$ is then given by
\[
\int_{G}f(x)\,dx=\int_{G/H}\int_{H}f(xh)\,dh\,d\dot{x},\quad f\in
\mathcal{C}_{c}(G)
\]
and similarly for $K/H\cap K$.

\begin{lemma}
\label{L:Int} Let the measures be normalized as above. Then the following hold:

\begin{enumerate}
\item [\hss\llap{\rm1)}]\label{L:Int(1)}Let $f\in\mathcal{C}_{c}(\bar{N}MAN)$.
Then
\[
\int_{G}f(g)\,dg=\int_{\bar{N}}\int_{M}\int_{A}\int_{N}f(\bar{n}man)a^{2\rho
}\,d\bar{n}\,dm\,da\,dn\,.
\]

\item[\hss\llap{\rm2)}] \label{L:Int(2)}Let $f\in\mathcal{C}_{c}(\bar{N})$.
For $y\in\bar{N}MAN$ write $y=\bar{n}(y)m_{\bar{N}}(y)a_{\bar{N}}(y)n_{\bar
{N}}(y)$. Let $x\in G$. Then
\[
\int_{\bar{N}}f(\bar{n}(x\bar{n}))a_{\bar{N}}(x\bar{n})^{-2\rho}\,d\bar
{n}=\int_{\bar{N}}f(\bar{n})\,d\bar{n}\,.
\]

\item[\hss\llap{\rm3)}] \label{L:Int(3)}Write, for $g\in G$,
$g=k(g)m(g)a(g)n(g)$ according to $G=KMAN$. Let $h\in\mathcal{C}(K/H\cap K)$.
Then
\[
\int_{K/H\cap K}h(\dot{k})\,d\dot{k}=\int_{\bar{N}}h(k(\bar{n})H\cap
K)a(\bar{n})^{-2\rho}\,d\bar{n}\,.
\]

\item[\hss\llap{\rm4)}] \label{L:Int(4)}Let $h\in\mathcal{C}(K/H\cap K)$ and
let $x\in G$. Then
\[
\int_{K/H\cap K}f(k(xk)H\cap K)a(xk)^{-2\rho}\,d\dot{k}=\int_{K/H\cap K}%
f(\dot{k})\,d\dot{k}%
\]

\item[\hss\llap{\rm5)}] \label{L:Int(5)}Assume that $\operatorname*{supp}%
(f)\subset H/H\cap K\subset K/H\cap K$. Then
\[
\int_{K/H\cap K}f(\dot{k})\,d\dot{k}=\int_{H/H\cap K}f(k(h)H\cap
K)a(h)^{-2\rho}\,d\dot{h}\,.
\]

\item[\hss\llap{\rm6)}] \label{L:Int(6)}Let $f\in\mathcal{C}_{c}(\bar{N})$.
Then
\[
\int_{\bar{N}}f(\bar{n})\,d\bar{n}=\int_{H/H\cap K}f(\bar{n}(h))a_{\bar{N}%
}(h)^{-2\rho}\,d\bar{n}\,.
\]

\item[\hss\llap{\rm7)}] \label{L:Int(7)}For $x\in HP_{\mathrm{max}}$ write
$x=h(x)m_{H}(x)a_{H}(x)n_{H}(x)$ with $h(x)\in H$, $m_{H}(x)\in M$,
$a_{H}(x)\in A$, and $n_{H}(x)\in N$. Let $f\in\mathcal{C}_{c}^{\infty
}(H/H\cap K)$ and let $x\in G$ be such that $xHP_{\mathrm{max}}\subset
HP_{\mathrm{max}}$. Then
\[
\int_{H/H\cap K}f(h(xh)H\cap K)a_{H}(xh)^{-2\rho}\,d\dot{h}=\int_{H/H\cap
K}f(\dot{h})\,d\dot{h}%
\]
\end{enumerate}
\end{lemma}

Identify $\frak{a}_{\mathbb{C}}^{\ast}$ with $\mathbb{C}$ by $\frak{a}%
_{\mathbb{C}}^{\ast}\ni\nu\mapsto2\nu(X^{0})\in\mathbb{C}\,$. Then $\rho$
corresponds to $\dim\frak{n}$. For $\nu\in\frak{a}_{\mathbb{C}}^{\ast}$, let
$\mathcal{C}^{\infty}(\nu)$ be the space of $\mathcal{C}^{\infty}$-functions
$f\colon G\rightarrow\mathbb{C}$ such that, for $a_{t}=\exp t(2X^{0})$, we
have
\[
f(gma_{t}n)=e^{-(\nu+\rho)t}f(g)=a_{t}^{-(\nu+\rho)}f(g)\,.
\]
Define an inner product on $\mathcal{C}^{\infty}(\nu)$ by
\[%
\ip{f}{g}%
_{\nu}:=\int_{K}\,\overline{f(k)}g(k)\,dk=\int_{K/H\cap K}\,\overline
{f(k)}g(k)\,d\dot{k}\,.
\]
Then $\mathcal{C}^{\infty}(\nu)$ becomes a pre-Hilbert space. We denote by
$\mathbf{H}(\nu)$ the completion of $\mathcal{C}^{\infty}(\nu)$. Define
$\pi(\nu)$ by
\[
\lbrack\pi(\nu)(x)f](g):=f(x^{-1}g),\quad x,g\in G,\quad f\in\mathcal{C}%
^{\infty}(\nu)\,.
\]
Then $\pi(\nu)(x)$ is bounded, so it extends to a bounded operator on
$\mathbf{H}(\nu)$, which we denote by the same symbol. Furthermore $\pi(\nu)$
is a continuous representation of $G$ which is unitary if and only if $\nu\in
i\mathbb{R}$. By \cite{Pou92} we have $\mathbf{H}(\nu)^{\infty}=\mathcal{C}%
^{\infty}(\nu)$. We can realize $\mathbf{H}(\nu)$ as $\mathbf{L}^{2}(K/H\cap
K)$ and as $\mathbf{L}^{2}(\bar{N},a(\bar{n})^{2\operatorname{Re}(\nu)}%
\,d\bar{n})$ by restriction (see Lemma \ref{L:5.14}). In the first realization
the representation $\pi(\nu)$ becomes
\[
\lbrack\pi(\nu)(x)f](k)=a(x^{-1}k)^{-\nu-\rho}f(k(x^{-1}k))
\]
and in the second
\[
\lbrack\pi(\nu)(x)f](\bar{n})=a_{\bar{N}}(x^{-1}\bar{n})^{-\nu-\rho}f(\bar
{n}(x^{-1}\bar{n}))\,.
\]

\begin{lemma}
\label{L:invpar} The pairing
\[
\mathbf{H}(\nu)\times\mathbf{H}(-\bar{\nu})\ni(f,g)\longmapsto%
\ip{f}{g}%
_{\nu}:=\int_{K}\overline{f(k)}g(k)\,dk=\int_{K/H\cap K}\overline
{f(k)}g(k)\,d\dot{k}%
\]
is $G$-invariant, i.e.
\[%
\ip{\pi(\nu)(x)f}{g}%
_{\nu}=%
\ip{f}{\pi(-\bar{\nu})(x^{-1})g}%
_{\nu}\,.
\]
\end{lemma}

\begin{remark}
\label{ImaginaryNu}We notice that if $\nu$ is purely imaginary, that is
$-\bar{\nu}=\nu$, the above shows that $(\pi(\nu),\mathbf{H}(\nu))$ is then unitary.
\end{remark}

\begin{lemma}
\label{L:5.14}{\ Let the notation be as above.}

\begin{enumerate}
\item [\hss\llap{\rm1)}]\label{L:5.14(1)}On $\bar{N}$ the invariant pairing $%
\ip{\,\cdot\,}{\,\cdot\,}%
_{\nu}$ is given by
\[%
\ip{f}{g}%
_{\nu}=\int_{\bar{N}}\overline{f(\bar{n})}g(\bar{n})\,d\bar{n}\,,\quad
f\in\mathbf{H}(\nu),\,g\in\mathbf{H}(-\bar{\nu})\,.
\]

\item[\hss\llap{\rm2)}] \label{L:5.14(2)}Let $\mathbf{H}_{H}(\nu)$ be the
closure of $\{f\in\mathcal{C}^{\infty}(\nu)\mid\operatorname*{supp}(f)\subset
HP_{\mathrm{max}}\}$. Then $\mathbf{H}_{H}(\nu)\ni f\mapsto f|_{H}%
\in\mathbf{L}^{2}(H/H\cap K,a(h)^{2\operatorname{Re}(\nu)}\,d\dot{h})$ is an isometry.

\item[\hss\llap{\rm3)}] \label{L:5.14(3)}Let $f\in\mathbf{H}(\nu)$,
$g\in\mathbf{H}(-\bar{\nu})$ and assume that $\operatorname*{supp}(fg)\subset
HP_{\mathrm{max}}$. Then
\[%
\ip{f}{g}%
_{\nu}=\int_{H/H\cap K}\overline{f(h)}g(h)\,d\dot{h}\,.
\]
\end{enumerate}
\end{lemma}

Let us assume, from now on, that there exists an element $w\in K$ such that
{$\operatorname{Ad}$}$(w)(X^{0})=-X^{0}$. In particular such an element exists
if $-1$ is in the Weyl group $W(\frak{a}_{q}):=N_{K}(\frak{a}_{q}%
)/Z_{K}(\frak{a}_{q})$, where $\frak{a}_{q}\subset\frak{p}\cap\frak{q}$ is
maximal abelian (and then maximal abelian in $\frak{p}$ and $\frak{q})$. This
is always the case if $G/H$ is a Cayley type space because $G$ is then a
Hermitian groups which implies that $\theta$ is an inner automorphism. The
element $w$ does also exists if $G$ is a complexification of one of the groups
$\frak{s}\frak{p}(n,\mathbb{R})$, $\frak{s}\frak{u}(n,n)$, $\frak{s}%
\frak{o}^{\ast}(4n)$, $\frak{s}\frak{o}(2,k)$ and $\frak{e}_{7(-25)}$, see
\cite{JO97}, Lemma 5.20.

Let us work out more explicitly the details for the representations of the
Cayley-type spaces in order to compare the existence of $(\rho_{\nu
},\mathbf{K}_{\nu})$ to the existence of the complementary series, see\ Lemma
\ref{CayleyConstant} below:

For ``big'' $\nu$ define the intertwining operator $A(\nu)\colon\mathbf{H}%
(\nu)\rightarrow\mathbf{H}(-\nu)$ by the converging integral%
\[
\lbrack A(\nu)f](x):=\int_{N}f(gnw)\,dn\,,
\]
see \cite{KS80,Wal92}. The map $\nu\mapsto A(\nu)$ has an analytic
continuation to a meromorphic function on $\frak{a}_{\mathbb{C}}^{\ast}$
intertwining $\pi(\nu)$ and $\pi(-\nu)$. Using Lemma \ref{L:invpar} we define
a new invariant bilinear form on $\mathcal{C}_{c}^{\infty}(\nu)$ by
\[%
\ip{f}{g}%
_{\nu}:=%
\ip{f}{A(\nu)g}%
_{\nu}.
\]
If there exists a (maximal) constant $R>0$ such that the invariant bilinear
form $%
\ip{\,\cdot\,}{\,\cdot\,}%
_{\nu}$ is positive definite for $|\nu|<R$, then we complete $\mathcal{C}%
_{c}^{\infty}(\nu)$ with respect to this new inner product, but denote the
resulting space by the same symbol $\mathbf{H}(\nu)$ as before. We call the
resulting unitary representations \textit{the complementary series}. Otherwise
we set $R=0$. We compare here the constant $R$ for the Cayley type symmetric
space, see \cite{OZ95,S92,S93}, and the constant $L_{\mathrm{pos}}+\rho$ from
Theorem \ref{CThw}. We notice that in all cases $R\leq L_{\mathrm{pos}}+\rho$.

\begin{lemma}
\label{CayleyConstant}For the Cayley-type symmetric spaces the constants $R$
and $L_{\mathrm{pos}}+\rho$ are given by the following table:
\[
\begin{tabular}
[c]{|lc|c|c|}\hline
\multicolumn{1}{|c}{$\frak{g}$}
&  & $R$ & $L_{\mathrm{pos}}+\rho$\\\hline
$\frak{su}(2n+1,2n+1)$ &  & $2n+1$ & $2n+1$\\
$\frak{su}(2n,2n)$ &  & $0$ & $2n$\\
$\frak{so}^{\ast}(4n)$ &  & $n$ & $2n$\\
$\frak{sp}(2n,\mathbb{R})$ &  & $n$ & $2n$\\
$\frak{sp}(2n+1,\mathbb{R})$ &  & $0$ & $2n+1$\\
$\frak{so}(4n+2,2)$ &  & $2$ & $2$\\
$\frak{so}(2n+1,2)$ &  & $1$ & $2$\\
$\frak{so}(4n,2)$ &  & $0$ & $2$\\
$E_{7(-25)}$ &  & $3$ & $3$\\\hline
\end{tabular}
\]
\end{lemma}

\begin{lemma}
\label{PhiUnimodular}$w^{-1}\tau(\bar{N})w=\bar{N}$, and $\varphi\colon\bar
{N}\ni\bar{n}\mapsto w^{-1}\tau(\bar{n})w\in\bar{N}$ is unimodular.
\end{lemma}

\begin{proof}
The first claim follows as {$\operatorname{Ad}$}$(w)$ and $\tau$ act by $-1$
on $\frak{a}$, and thus map $N$ onto $\bar{N}$, and $\bar{N}$ onto $N$. The
second follows as we can realize $\varphi^{2}$ by conjugation by an element in
$M\cap K$.
\end{proof}

\begin{lemma}
\label{UnitaryIsomorphism}For $f\in\mathbf{H}(\nu)$ let $J(f)(x):=f(\tau
(xw))$. Then the following properties hold:

\begin{enumerate}
\item [\hss\llap{\rm1)}]\label{UnitaryIsomorphism(1)}$J(f)(x)=f(\tau
(x)w^{-1})$.

\item[\hss\llap{\rm2)}] \label{UnitaryIsomorphism(2)}$J(f)\in\mathbf{H}(\nu)$
and $A(\nu)J=JA(\nu)$.

\item[\hss\llap{\rm3)}] \label{UnitaryIsomorphism(3)}$J\colon\mathbf{H}%
(\nu)\rightarrow\mathbf{H}(\nu)$ is a unitary isomorphism.

\item[\hss\llap{\rm4)}] \label{UnitaryIsomorphism(4)}$J^{2}=\operatorname*{id}%
$.

\item[\hss\llap{\rm5)}] \label{UnitaryIsomorphism(5)}For $x\in G$ we have
$J\circ\pi(\nu)(x)=\pi(\nu)(\tau(x))\circ J$.
\end{enumerate}
\end{lemma}

Notice that
\begin{equation}
\lbrack A(\nu)J](f)(x):=\int_{\bar{N}}\,f(\tau(x)\bar{n})\,d\bar{n}
\label{E:AJ}%
\end{equation}
for $\operatorname{Re}\lambda$ ``big''. By simple calculation we get:

\begin{lemma}
\label{L:AnuJ} Assume that $G/H$ is non-compactly causal. Then $A(\nu)J$
intertwines $\pi(\nu)$, and $\pi(-\nu)\circ\tau$ if $A(\nu)J$ has no pole at
$\nu$.
\end{lemma}

Equation (\ref{E:AJ}) and Lemma \ref{L:AnuJ} show that even if each one of the
operators $A(\nu)$ and $J$ does not exist, the operator $A(\nu)J\colon
\mathbf{H}(\nu)\rightarrow\mathbf{H}(-\nu)$, $[A(\nu)J]\circ\pi(\nu)=[\pi
(-\nu)\circ\tau][A(\nu)J]$, will always exist. The next theorem shows that the
intertwining operator $A(\nu)J$ is a convolution operator which kernel
$y,x\mapsto a_{\bar{N}}(\tau(y)^{-1}x)^{\nu-\rho}$. The importance of that is
that $\bar{N}MAN=(P^{-}K_{\mathbb{C}}^{c}P^{+})\cap G$ and hence%
\begin{equation}
a_{\bar{N}}(x)^{\nu-\rho}=k_{\mathbb{C}}(\tau(x))^{-\nu+\rho}
\label{E:aNbarkC}%
\end{equation}
where $k_{\mathbb{C}}$ denotes the $K_{\mathbb{C}}^{c}-$projection from
(\ref{E:P+KP-projections}). Here we have used that $\tau(X^{0})=-X^{0}$. The
reflection positivity then reduces to the problem to determine those $\nu$ for
which this kernel is positive semidefinite.

\begin{theorem}
\label{ANuJf}\ 

\begin{enumerate}
\item \label{ANuJf(1)}Let $f\in\mathcal{C}^{\infty}(\nu)$. Then
\[
\lbrack A(\nu)J](f)(\bar{n})=\int_{\bar{N}}f(x)a_{\bar{N}}(\tau(\bar{n}%
)^{-1}x)^{\nu-\rho}\,dx\,.
\]

\item \label{ANuJf(2)}If $\operatorname*{supp}(f)\subset HP_{\mathrm{max}}$,
then for $h\in H$
\[
\lbrack A(\nu)J](f)(h)=\int_{H/H\cap K}f(x)a_{\bar{N}}(h^{-1}x)^{\nu-\rho
}\,d\dot{x}\,.
\]
\end{enumerate}
\end{theorem}

\begin{proof}
We may assume that $\nu$ is big enough such that the integral defining
$A(\nu)$ converges. The general statement follows then by analytic
continuation. We have
\begin{align*}
\lbrack A(\nu)J]f(\bar{n})  &  =\int_{\bar{N}}Jf(\bar{n}wx)\,dx\\
&  =\int f(\tau(\bar{n})w^{-1}\tau(x)w)\,dx\\
&  =\int f(\tau(\bar{n})x)\,dx\\
&  =\int f(\bar{n}(\tau(\bar{n})x))a_{\bar{N}}(\tau(\bar{n})x)^{-(\nu+\rho
)}\,dx\,.
\end{align*}
Now $a_{\bar{N}}(\tau(\bar{n})x)=a_{\bar{N}}(\tau(\bar{n})^{-1}\bar{n}%
(\tau(\bar{n})x))^{-1}$. By Lemma \ref{L:Int} we get
\begin{align*}
\lbrack A(\nu)J]f(\bar{n})  &  =\int f(\bar{n}(\tau(\bar{n})x))a_{\bar{N}%
}(\tau(\bar{n})^{-1}\bar{n}(\tau(\bar{n})x))^{\nu-\rho}a_{\bar{N}}(\tau
(\bar{n})x)^{-2\rho}\,dx\\
&  =\int f(x)a_{\bar{N}}(\tau(\bar{n})^{-1}x)^{\nu-\rho}\,dx.
\end{align*}

The second statement follows in the same way. Let $\mathcal{C}_{c}^{\infty
}(\Omega,\nu)$ be the subspace of $\mathcal{C}^{\infty}(\nu)$ consisting of
functions $f$ with support in $\exp\Omega P_{\max}$, and with
$\operatorname*{supp}(f|_{\Omega})$ compact.
\end{proof}

\begin{lemma}
\label{fgJ}Let $f,g\in\mathcal{C}_{c}^{\infty}(\Omega,\nu)$. Then the
following holds:

\begin{enumerate}
\item \label{fgJ(1)}$%
\ip{f}{g}%
_{J}=\int_{\bar{N}}\int_{\bar{N}}\overline{f(x)}g(y)a_{\bar{N}}(\tau
(x)^{-1}y)^{\nu-\rho}\,dx\,dy\,$.

\item \label{fgJ(2)}$%
\ip{f}{g}%
_{J}=\int_{H/H\cap K}\int_{H/H\cap K}\overline{f(h)}g(k)a_{\bar{N}}%
(h^{-1}k)^{\nu-\rho}\,dh\,dk\,$.
\end{enumerate}
\end{lemma}

Let $\sigma\colon\frak{g}_{\mathbb{C}}^{c}\rightarrow\frak{g}_{\mathbb{C}}%
^{c}$ be the conjugation with respect to $\frak{g}^{c}$. Then $\sigma
|_{\frak{g}}=\tau|_{\frak{g}}$. This, the fact that $\Omega=\sigma
(\Omega_{\mathbb{C}})\cap\frak{g}$, and Theorem \ref{S:DF1}, part
(\ref{S:DF1(1)}), implies that%
\begin{equation}
(X,Y)\longmapsto a_{\bar{N}}(\tau(\exp(X))^{-1}\exp(Y))^{\nu-\rho}%
=Q(\sigma(Y),\sigma(X))=:Q_{\sigma}\left( Y,X\right) ,\label{eqApr.38}
\end{equation}
where $Q$ is the reproducing kernel of the representation $\rho_{\nu}$, see
(\ref{E:aNbarkC}). Notice the twist by $\sigma$ that comes from the fact that
we are using the realization of $\Omega$ inside $\frak{p}^{-}$ whereas we have
realized $\Omega_{\mathbb{C}}$ inside $\frak{p}^{+}$. It follows that the
integral kernel is positive definite. Hence $%
\ip{\,\cdot\,}{\,\cdot\,}%
_{J}$ is positive definite on the space of functions supported in $HP_{\max}$.
We let $\mathbf{K}_{0}$ be the completion of $\mathcal{C}_{c}^{\infty}%
(\Omega,\nu)$.

What is still needed for the application of the L\"uscher-Mack Theorem is
the invariant cone $C\subset\frak{q}$ and the semigroup $S$. For that we
choose $C_{\min}$ as the minimal invariant cone in $\frak{q}$ containing
$X^{0}$, that is $C_{\min}$ is generated by $\operatorname{Ad}(H)X^{0}$. We
notice that%
\[
\operatorname{Ad}(\exp(tX^{0}))X=e^{-t}X\,,\qquad\forall X\in\frak{\bar{n}%
}\text{.}%
\]
Hence $\operatorname{Ad}(\exp tX^{0})$ acts by contractions on $\Omega$ if
$t>0$. Let
\[
S(\Omega):=\{g\in G\mid gH\subset HP_{\mathrm{max}}\}=\left\{  g\in G\mid
g\cdot\Omega\subset\Omega\right\}  \,.
\]
Then $S(\Omega)$ is a closed semigroup invariant under $s\mapsto s^{\sharp
}:=\tau(s)^{-1}$. It follows by construction that $S(\Omega)\subset HP_{\max}%
$. We remark the following results:

\begin{lemma}
Let $C=C_{\mathrm{max}}$ be the maximal pointed generating cone in $\frak{q}$
containing $X^{0}$. Then the following hold:

\begin{enumerate}
\item  Let $t>0$ and $Y\in\Omega$. Then $\exp tX^{0}\in S$ and $\exp
tX^{0}\cdot Y=e^{-t}Y$.

\item $S(\Omega)=H\exp C_{\mathrm{max}}$.
\end{enumerate}
\end{lemma}

\begin{proof}
(1) is a simple calculation. For (2) see \cite{HiNe93} and \cite{HO95}.
\end{proof}

\begin{corollary}
\label{CorSsspNew.17}The semigroup $S(C_{\min})=H\exp(C_{\min})$ acts by
contractions on $\Omega$.
\end{corollary}

\begin{lemma}
\label{CSInvariant}Let $s\in S(\Omega)$ and $f\in\mathcal{C}_{c}^{\infty
}(\Omega,\nu)$. Then $\pi(\nu)(s)f\in\mathcal{C}_{c}^{\infty}(\Omega,\nu)$,
that is $\mathcal{C}_{c}^{\infty}(\Omega)$ is $S(\Omega)$-invariant.
\end{lemma}

\begin{proof}
Let $f\in\mathcal{C}_{c}^{\infty}(\Omega)$ and $s\in S$. Then $\pi
(\nu)(s)f(x)=f(s^{-1}x)\not =0$ only if $s^{-1}x\in\operatorname*{supp}%
(f)\subset HP_{\mathrm{max}}$. Thus $\operatorname*{supp}(\pi(\nu)(s)f)\subset
s\operatorname*{supp}(f)\subset sHP_{\mathrm{max}}\subset HP_{\mathrm{max}}$.
\end{proof}

We still assume that $G^{c}$ is simple. Let $(\rho_{\nu},\mathbf{K}_{\nu})$ be
as above. Let $%
\openone
\in\mathbf{K}_{\nu}(\lambda-\rho)$ be the constant function $Z\mapsto1$. Then
$\left\|
\openone
\right\|  =1$. Let $H^{c}:=\left(  G^{c}\right)  ^{\tau}$. Then $H^{c}$ is
connected. Let $\tilde{H}$ be the universal covering of $H^{c}$ and then also
$H_{o}$. We notice that
\[
H^{c}/H^{c}\cap K^{c}=H/H\cap K\,.
\]
Denote the restriction of $\rho_{\nu}$ to $H^{c}$ by $\rho_{\nu,H}$. We can
lift $\rho_{\nu,H}$ to a representation of $\tilde{H}$ also denoted by
$\rho_{\nu,H}$. We let $C=C_{\mathrm{min}}$ be the minimal $H$-invariant cone
in $\frak{{q}}$ generated by $X^{0}$. We denote by $\tilde{C}=\tilde
{C}_{\mathrm{min}}$ the minimal $G^{c}$-invariant cone in $i\frak{{g}^{c}}$
with
$\tilde{C}\cap\frak{q}=\operatorname*{pr}{}_{\frak{q}}(\tilde{C})=C$,
where $\operatorname*{pr}_{\frak{q}}\colon\frak{g}\rightarrow\frak{q}$ denotes
the orthogonal projection (see\ \cite{HO95,'O90b}). As $L_{\mathrm{pos}}\leq0$
it follows that $\rho_{\lambda}$ extends to a holomorphic representation of
the universal semigroup $\Gamma(G^{c},\tilde{C})$ corresponding to $G^{c}$ and
$\tilde{C}$
(see \cite{HiNe93,KNO98}).
Let $G_{1}^{c}$ be the analytic subgroup of
$G_{\mathbb{C}}$ corresponding to the Lie algebra $\frak{g}^{c}$. Let $H_{1}$
be the analytic subgroup of $G_{1}^{c}$ corresponding to $\frak{h}$. Then---as
we are assuming that $G\subset G_{\mathbb{C}}$---we have $H_{1}=H_{o}$. Let
$\kappa\colon G^{c}\rightarrow G_{1}^{c}$ be the canonical projection and let
$Z_{H}=\kappa^{-1}(Z_{G_{1}^{c}}\cap H_{o})$. Then $\rho_{\nu}$ is trivial on
$Z_{H}$ as $\nu-\rho$ is trivial on $\exp([\frak{k}^{c},\frak{k}^{c}])\supset
H^{c}\cap K^{c}$. Thus $\rho_{\nu}$ factors to $G^{c}/Z_{H}$, and to
$\Gamma(G^{c},\tilde{C})/Z_{H}$. Notice that $(G^{c}/Z_{H})_{o}^{\tau}$ is
isomorphic to $H_{o}$. Therefore we can view $H_{o}$ as subgroup of
$G^{c}/Z_{H}$, and $S_{o}(C)=H_{o}\exp C$ as a subsemigroup of $\Gamma
(G^{c},\tilde{C})/Z_{H}$. In particular $\tau_{\nu}(s)$ is defined for $s\in
S_{o}(C)$. This allows us to write $\rho_{\nu}(h)$ or $\rho_{\nu,H}(h)$ for
$h\in H_{o}$. Using (\ref{E:aNbarkC}) and Lemma \ref{LemNewS-hwm.2} we get
\[
a_{\bar{N}}(h)^{\nu-\rho}=\left\langle
\openone
,\rho_{\nu,H}(h)%
\openone
\right\rangle \,.
\]

In particular we get that $(h,k)\mapsto a_{\bar{N}}(h^{-1}k)^{\nu-\rho}$ is
positive semidefinite if $\nu-\rho\le L_{\mathrm{pos}}$.

Let us now consider the case $G=H_{\mathbb{C}}$ and $G^{c}=\tilde{H}%
\times\tilde{H}$. Denote the constant $L_{\mathrm{pos}}$ for $\tilde{H}$ by
$S_{\mathrm{pos}}$ and denote, for $\mu\leq S_{\mathrm{pos}}$, the
representation with lowest $\tilde{H}\cap\tilde{K}$-type $\mu$ by $(\tau_{\mu
},L_{\mu})$. Let $\bar{\tau}_{\mu}$ be the conjugate representation. Recall
that we view $\tilde{H}$ as a subset of $G^{c}$ by the diagonal embedding
\[
\tilde{H}\ni h\longmapsto(h,h)\in\Delta(G^{c}):=\{(x,x)\in G^{c}\mid
x\in\tilde{H}\}\,.
\]
The center of $\frak{k}^{c}$ is two dimensional (over $\mathbb{R}$) and
generated by $i(X^{0},X^{0})$ and $i(X^{0},-X^{0})$. We choose $Z^{0}%
=i(X^{0},-X^{0})$. Then $\frak{p}^{+}=\frak{n}\times\bar{\frak{n}}$. Let
$\openone$
again be a lowest weight vector of norm one. Denote the corresponding vector
in the conjugate Hilbert space by
$\overline{\openone}$.
Then for $h\in\tilde{H}$:%
\begin{align*}%
\ip{\openone\otimes\overline{\openone}}{\tau_{\mu}\otimes\bar{\tau}_{\mu
}(h,h)\openone\otimes\overline{\openone}}%
&  =%
\ip{\openone}{\tau_{\lambda}(h)\openone}%
\overline{%
\ip{\openone}{\tau_{\lambda}(h)\openone}%
}\\
&  =\left|
\ip{\openone}{\tau_{\lambda}(h)\openone}%
\right|  ^{2}\\
&  =a_{\bar{N}}(h)^{2\mu}%
\end{align*}

Thus we define in this case $L_{\mathrm{pos}}:=2S_{\mathrm{pos}}$. As before
we notice that
$\tau_{\nu}\otimes\bar{\tau}_{\nu}(h,h)\openone\otimes\overline{\openone}$ is well
defined on $H$. We now have a new proof that $%
\ip{\,\cdot\,}{\,\cdot\,}%
$ is positive definite on $\Omega$.

\begin{lemma}
\label{ahaloN} For $\nu-\rho\leq L_{\mathrm{pos}}$ there exists an unitary
irreducible highest weight representation $(\rho_{\nu},\mathbf{K}_{\nu})$ of
$G^{c}$ and a lowest $K^{c}$-type vector $%
\openone
$ of norm one such that for every $h\in H$%
\[
a_{\bar{N}}(h)^{\nu-\rho}=%
\ip{\openone}{\tau_{\lambda}(h)\openone}%
\,.
\]
Hence the kernel
\[
(H\times H)\ni(h,k)\longmapsto a_{\bar{N}}(k^{-1}h)^{\nu-\rho}\in\mathbb{R}%
\]
is positive semidefinite. In particular $%
\ip{\,\cdot\,}{\,\cdot\,}%
_{J}$ is positive semidefinite on $\mathbf{C}_{c}^{\infty}(\Omega,\nu)$ for
$\nu-\rho\leq L_{\mathrm{pos}}$.
\end{lemma}

The Basic Lemma and the L\"uscher-Mack Theorem, together with the above, now
imply the following Theorem:

\begin{theorem}
[Reflection Symmetry for Complementary Series]\label{PR} Let $G/H$%
\linebreak be a non-compactly causal symmetric space such that there exists a
$w\in K$ such that $\operatorname{Ad}(w)|_{\frak{a}}=-1$. Let $\pi(\nu)$ be a
complementary series such that $\nu\leq L_{\mathrm{pos}}$. Let $C$ be the
minimal $H$-invariant cone in $\frak{q}$ such that $S(C)$ is contained in the
contraction semigroup of $HP_{\mathrm{max}}$ in $G/P_{\mathrm{max}}$. Let
$\Omega$ be the bounded realization of $H/H\cap K$ in $\bar{\frak{n}}$. Let
$J(f)(x):=f(\tau(x)w^{-1})$. Let $\mathbf{K}_{0}$ be the closure of
$\mathcal{C}_{c}^{\infty}(\Omega,\nu)$ in $\mathbf{H}(\nu)$. Then the
following hold:

\begin{enumerate}
\item [\hss\llap{\rm1)}]\label{PR(1)}$(G,\tau,\pi(\nu),C,J,\mathbf{K}_{0})$
satisfies the positivity conditions \textup{(PR1)--(PR2).}

\item[\hss\llap{\rm2)}] \label{PR(2)}$\pi(\nu)$ defines a contractive
representation $\tilde{\pi}(\nu)$ of $S(C)$ on $\mathbf{K}$ such that
\[
\tilde{\pi}(\nu)(\gamma)^{\ast}=\tilde{\pi}(\nu)(\tau(\gamma)^{-1}).
\]

\item[\hss\llap{\rm3)}] \label{PR(3)}There exists a unitary representation
$\tilde{\pi}^{c}$ of $G^{c}$ such that

\begin{enumerate}
\item [\hss\llap{\rm i)}]\label{PR(3)(1)}$d\tilde{\pi}(\nu)^{c}(X)=d\tilde
{\pi}(\nu)(X)\,\quad\forall X\in\frak{h}$.

\item[\hss\llap{\rm ii)}] \label{PR(3)(2)}$d\tilde{\pi}(\nu)^{c}%
(iY)=i\,d\tilde{\pi}(\nu)(Y)\,\quad\forall Y\in C$.
\end{enumerate}
\end{enumerate}
\end{theorem}

We remark that this Theorem includes the results of R. Schrader for the
complementary series of $SL(2n,\mathbb{C})$ \cite{Sch86}. In a moment we will
show that actually $\tilde{\pi}(\nu)^{c}\simeq\rho_{\nu}$, where $\rho_{\nu}$
is the irreducible unitary highest weight representation of $G^{c}$ such that
\[
a(h)^{\nu-\rho}=%
\ip{\openone}{\rho_{\nu}(h)\openone}%
\]
as before. {}From now on we assume that $\nu-\rho\leq L_{\mathrm{pos}}$. We
notice that the inner product $%
\ip{\,\cdot\,}{A(\nu)J(\,\cdot\,)}%
$ makes sense independent of the existence of $w$. Let $\mathbf{K}_{0}$ be the
completion of $\mathcal{C}_{c}^{\infty}(\Omega,\nu)$ in the norm $%
\ip{\,\cdot\,}{A(\nu)J(\,\cdot\,)}%
$. Let $\mathbf{N}$ be the space of vectors of zero length and let
$\mathbf{K}$ be the completion of $\mathbf{K}_{0}/\mathbf{N}$ in the induced
norm. First of all we have to show that $\pi(\nu)(\gamma)$ passes to a
continuous operator $\tilde{\pi}(\nu)(\gamma)$ on $\mathbf{K}$ such that
$\tilde{\pi}(\nu)(\gamma)^{\ast}=\tilde{\pi}(\nu)(\tau(\gamma)^{-1})$. For
that we recall that
\begin{equation}
H/H\cap K=H_{o}/H_{o}\cap K=\Omega\label{eqSsspNew.14}%
\end{equation}
so we my replace the integration over $H$ in $%
\ip{f}{A(\nu)Jf}%
_{\nu}$ with integration over $H_{o}$. Motivated by (\ref{E:R*one})\ we define
for $f\in\mathcal{C}_{c}^{\infty}(\Omega,\nu)$
\begin{align}
U(f)=\rho_{\nu}(f)%
\openone
&  :=\int_{H_{o}}f(h)\rho_{\nu}(h)%
\openone
\,dh\,\label{eqSsspNew.15}\\
&  =\int_{H_{o}}a_{\bar{N}}(h)^{-\nu-\rho}f(h\cdot0)J_{\nu}(h,\cdot)^{-1}\,dh.
\label{E:U}%
\end{align}

\begin{lemma}
\label{L:posref} Assume that $\nu-\rho\leq L_{\mathrm{pos}}$. Let $\rho_{\nu
},\,\mathbf{K}_{\nu}$ and $%
\openone
$ be as specified in Lemma \textup{\ref{ahaloN}} and let $f,g\in
\mathcal{C}_{c}^{\infty}(\Omega,\nu)$ and $s\in S(C)$. Then the following hold:

\begin{enumerate}
\item [\hss\llap{\rm1)}]\label{posref(1)}$%
\ip{f}{[A(\nu)J](g)}%
_{\nu}=%
\ip{Uf}{Ug}%
$.

\item[\hss\llap{\rm2)}] \label{posref(2)}$U(\pi(\nu)(s)f)=\rho_{\nu}(s)U(f)$.

\item[\hss\llap{\rm3)}] \label{posref(3)}$\pi_{\nu}(s)$ passes to a
contractive operator $\frak{\pi(\nu)}(s)$ on $\mathbf{K}$ such that
$\tilde{\pi}(\nu)(s)^{\ast}=\tilde{\pi}(\nu)(\tau(s)^{-1})$.
\end{enumerate}
\end{lemma}

\begin{proof}
(1) Let $f$ and $g$ be as above. Then
\begin{align*}%
\ip{f}{[A(\nu)J](g)}%
&  =\int_{H_{o}/H_{o}\cap K}\int_{H_{o}/H_{o}\cap K}\overline{f(h)}%
g(k)a_{\bar{N}}(h^{-1}k)^{\nu-\rho}\,dh\,dk\\
&  =\int_{H_{o}/H_{o}\cap K}\int_{H_{o}/H_{o}\cap K}\overline{f(h)}g(k)%
\ip{\openone}{\rho_{\nu}(h^{-1}k)\openone}%
\,dh\,dk\\
&  =\int_{H_{o}/H_{o}\cap K}\int_{H_{o}/H_{o}\cap K}\overline{f(h)}g(k)%
\ip{\rho_{\nu}(h)\openone}{\rho_{\nu}(k)\openone}%
\,dh\,dk\\
&  =%
\ip{Uf}{Ug}%
.
\end{align*}
This proves (1).

(2) This follows from Lemma \ref{L:Int},7) and the following calculation:
\begin{align*}
U(\pi_{\nu}(s)f)%
\openone
&  =\int f(s^{-1}h)\rho_{\nu}(h)%
\openone
\,dh\\
&  =\int f(h(s^{-1}h))a_{H}(s^{-1}h)^{-(\nu+\rho)}\rho_{\nu}(h)%
\openone
\,dh\\
&  =\int f(h(s^{-1}h))a_{H}(sh(s^{-1}))^{\nu-\rho}\rho_{\nu}(h)a_{H}%
(s^{-1}h)^{-2\rho}%
\openone
\,dh\\
&  =\int f(h)\rho_{\nu}(sh)%
\openone
\,dh\\
&  =\rho_{\nu}(s)U(f)\,,
\end{align*}
where we have used that
\[
\rho_{\nu}(sh)%
\openone
=a_{H}(sh)^{\nu-\rho}\rho_{\nu}(h(sh))%
\openone
\,.
\]

(3) By (1) and (2) we get:
\begin{align*}
\Vert\pi_{\nu}(s)f\Vert_{J}^{2}  &  =\Vert\rho_{\nu}(s)U(f)\Vert^{2}\\
&  \leq\Vert U(f)\Vert^{2}\\
&  =%
\ip{f}{[A(\nu)J]f}%
_{\nu}\qquad\left(  =\Vert f\Vert_{J}^{2}\right)  \,.
\end{align*}
Thus $\pi_{\nu}(s)$ passes to a contractive operator on $\mathbf{K}$. That
$\tilde{\pi}_{\nu}(s)^{\ast}=\tilde{\pi}_{\nu}(\tau(s)^{-1})$ follows from
Lemma \ref{L:AnuJ}.
\end{proof}

\begin{theorem}
[Identification Theorem]\label{S:Posref} Assume that $G/H$ is non-compactly
causal and that $\nu-\rho\leq L_{\mathrm{pos}}$. Let $\rho_{\nu}$,
$\mathbf{K}_{\nu}$ and $%
\openone
\in\mathbf{K}_{\nu}$ be as in Lemma \textup{\ref{ahaloN}.} Then the following hold:

\begin{enumerate}
\item [\hss\llap{\rm1)}]\label{S:Posref(1)}There exists a continuous
contractive representation $\frak{\pi(\nu)}$ of $S_{o}(C)$ on $\mathbf{K}$
such that
\[
\frak{\pi(\nu)}(s)^{\ast}=\frak{\pi(\nu)}(\tau(s)^{-1})\,,\quad\forall s\in
S_{o}(C)\,.
\]

\item[\hss\llap{\rm2)}] \label{S:Posref(2)}There exists a unitary
representation $\tilde{\pi}(\nu)^{c}$ of $G^{c}$ such that

\begin{enumerate}
\item [\hss\llap{\rm i)}]\label{S:Posref(2)(1)}$d\tilde{\pi}(\nu
)^{c}(X)=d\tilde{\pi}(\nu)(X)\,\quad\forall X\in\frak{h}$.

\item[\hss\llap{\rm ii)}] \label{S:Posref(2)(2)}$d\tilde{\pi}(\nu
)^{c}(iY)=i\,d\tilde{\pi}(\nu)(Y)\,\quad\forall Y\in C$.
\end{enumerate}

\item[\hss\llap{\rm3)}] \label{S:Posref(3)}The map
\[
\mathcal{C}_{c}^{\infty}(\Omega,\nu)\ni f\longmapsto U(f)\in\mathbf{K}_{\nu}%
\]
extends to an isometry $\mathbf{K}\simeq\mathbf{K}_{\nu}$ intertwining
$\tilde{\pi}(\nu)^{c}$ and $\rho_{\nu}$. In particular $\tilde{\pi}(\nu)^{c}$
is irreducible and isomorphic to $\rho_{\nu}$.
\end{enumerate}
\end{theorem}

\begin{proof}
(1) follows from Lemma \ref{L:posref} as obviously $\tilde{\pi}(\nu
)(sr)=\tilde{\pi}(\nu)(s)\tilde{\pi}_{\nu}(r)$.

(2) This follows now from the Theorem of L\"uscher-Mack.

(3) By Lemma \ref{L:posref} we know that $f\mapsto U(f)$ defines an isometric
$S_{o}(C)$-intertwining operator. Let $f\in\mathcal{C}_{c}^{\infty}(\Omega
,\nu)$. Differentiation and the fact that $\rho_{\nu}$ is holomorphic gives

\begin{enumerate}
\item [\hss\llap{\rm i)}]$U(d\tilde{\pi}(\nu)^{c}(X)f)=d\rho_{\nu
}(X)U(f)\,,\quad\forall X\in\frak{h}$.

\item[\hss\llap{\rm ii)}] $U(i\,d\tilde{\pi}(\nu)^{c}(Y)f)=i\,d\rho_{\nu
}(Y)U(f)\,,\quad\forall Y\in C$.
\end{enumerate}

But those are exactly the relations that define $\tilde{\pi}(\nu)^{c}$. The
fact that $\frak{h}\oplus iC$ generates $\frak{g}^{c}$ implies that $f\mapsto
U(f)$ induces an $\frak{g}^{c}$-intertwining operator intertwining $\tilde
{\pi}(\nu)^{c}$ and $\rho_{\nu}$. As both are also representations of $G^{c}$,
it follows that this is an isometric $G^{c}$-map. In particular as it is an
isometry by Lemma \ref{L:posref}, part 1, it follows that the map $U$ is an
isomorphism. This proves the theorem.
%
%
%
%
%
\end{proof}

We will now explain another view of the above results using local
representations instead of the L\"uscher-Mack Theorem. Let $\frak{a}_{p}$ be
a maximal abelian subspace of $\frak{p}$ containing $X^{0}$. Then
$\frak{a}_{p}$ is contained in $\frak{q}$. Let $\Delta(\frak{g},\frak{a})$ be
the set of roots of $\frak{a}$ in $\frak{g}$. We choose a positive system such
that $\Delta_{+}=\left\{  \alpha\mid\alpha(X^{0})=1\right\}  \subset\Delta
^{+}(\frak{g},\frak{a})$. Choose a maximal set of long strongly orthogonal
roots $\gamma_{1},\dots,\gamma_{r}$, $r=rank(H/H\cap K)$. Choose $X_{j}%
\in\frak{g}_{\gamma_{j}}$ such that with $X_{-j}=\tau(X_{j})$ we have
$[X_{j},X_{-j}]=H_{j}:=H_{\gamma_{j}}$. Then by Theorem 5.1.8 in \cite{HO95}
we have%
\[
\Omega={\operatorname{Ad}}(H\cap K)\left\{  \sum_{j=1}^{r}t_{j}X_{-j}\biggm
|-1<t_{j}<1,\,1\leq j\leq r\right\}  \,.
\]
For $R>0$, let
\[
B_{R}:={\operatorname{Ad}}(H\cap K)\left\{  \sum_{j=1}^{r}t_{j}X_{-j}\biggm
|-R<t_{j}<R\right\}  \,.
\]
Then $B_{R}$ is open in $\bar{\frak{n}}$. Let $\beta\colon\mathbf{K}%
_{0}\rightarrow\left(  \mathbf{K}_{0}/\mathbf{N}\right)  \sptilde
=\mathbf{K}$ be the canonical map. Then $\beta$ is a contraction ($\left\|
\beta(f)\right\|  _{J}^{2}=%
\ip{f}{Jf}%
=%
\ip{f}{f}%
_{J}
\leq\left\|  f\right\|  ^{2}$). For $U\subset\Omega$ open, let
\[
\mathcal{C}_{c}^{\infty}(U,\nu):=\left\{  f\in\mathcal{C}_{c}^{\infty}%
(\Omega,\nu)\mid\text{supp}(f)\subset U\right\}
\]
and $\mathbf{K}(U):=\beta(\mathcal{C}_{c}^{\infty}(U))$.

\begin{theorem}
\label{Dense}Let $U\subset\Omega$ be open. Then $\mathbf{K}(U)$ is dense in
$\mathbf{K}$.
\end{theorem}

\begin{proof}
Let $x\in U$. Then we can choose $h\in H$ such that $h\cdot x=0$. As
$\mathcal{C}_{c}^{\infty}(U,\nu)=h\cdot C_{c}^{\infty}(h\cdot U,\nu)$ and $H$
acts unitarily, it follows that we can assume that $0\in U$. Let $R>0$ be such
that $B_{R}\subset U$. Then $\mathcal{C}_{c}^{\infty}(B_{R},\nu)\subset
\mathcal{C}_{c}^{\infty}(U,\nu)$. Hence we can assume that $U=B_{R}$. Let
$g\in\mathcal{C}_{c}^{\infty}(U,\nu)^{\perp}$ and let $f\in\mathcal{C}%
_{c}^{\infty}(\Omega,\nu)$. We want to show that $%
\ip{g}{f}%
_{J}=0$. Choose $0<L<1$ such that $\operatorname*{supp}(f)\subset B_{L}$. For
$t\in\mathbb{R}$ and $a_{t}=\exp(2tX^{0})$ we have $a_{t}\cdot B_{L}%
=B_{e^{-2t}L}$. Thus $\operatorname*{supp}([\pi(\nu)(a_{t})f]_{\Omega})\subset
B_{e^{-2t}L}$. Choose $0<s_{0}$ such that $e^{-2t}L<R$ for every $t>s_{0}$.
Then $\pi(\nu)(a_{t})(f)\in\mathcal{C}_{c}^{\infty}(U,\nu)$ for every
$t>s_{0}$. It follows that for $t>s_{0}$:
\begin{align*}
0  &  =%
\ip{g}{\pi(\nu)(a_{t})f}%
_{J}\\
&  =\int_{\Omega}\int_{\Omega}\overline{g(x)}\left[  \pi(\nu)(a_{t})f\right]
(y)
Q_{\sigma}(x,y)
\,dx\,dy\\
&  =e^{(\lambda+1)t}\int_{\Omega}\int_{\Omega}\overline{g(x)}f(e^{2t}%
y)
Q_{\sigma}(x,y)
\,dx\,dy\\
&  =e^{(\lambda-1)t}\int_{\Omega}\int_{\Omega}\overline{g(x)}f(y)
Q_{\sigma}(x,e^{-2t}y)
\,dx\,dy\,.
\end{align*}
By Lemma \ref{S:DF1} we know that $z\mapsto Q(zX,Y)$ is holomorphic on
$D=\left\{  z\mid\left|  z\right|  <1\right\}  $. As $g$ and $f$ both have
compact support it follows
by (\ref{eqApr.38})
that
\[
F(z):=\int_{\Omega}\int_{\Omega}\overline{g(x)}f(y)Q(x,zy)\,dx\,dy
\]
is holomorphic on $D$. But $F(z)=0$ for $0<z<e^{-2s_{0}}$. Thus $F(z)=0$ for
every $z$. In particular
\[%
\ip{g}{\pi(\nu)(a_{t})f}%
_{J}=0
\]
for every $t>0$. By continuity $%
\ip{g}{f}%
_{J}=0$. Thus $g=0$.
\end{proof}

Let us recall some basic facts from \cite{Jor86}. Let $\rho$ be a local
homomorphism of a neighborhood $U$ of $e$ in $G$ into the space of linear
operators on the Hilbert space $\mathbf{H}$ such that $\rho(g)$ is densely
defined for $g\in U$. Furthermore $\rho|_{(U\cap H)}$ extends to a strongly
continuous representation of $H$ in $\mathbf{H}$. $\rho$ is called a
\emph{local representation} if there exists a dense subspace $\mathbf{D}%
\subset\mathbf{H}$ such that the following hold: \begin{list}{}{\setlength
{\leftmargin}{\customleftmargin}
\setlength{\itemsep}{0.5\customskipamount}
\setlength{\parsep}{0.5\customskipamount}
\setlength{\topsep}{\customskipamount}}
\item[\hss\llap{\rm LR1)}]  $\forall g\in U$, ${\bf D} \subset{\bf D} (\rho
(g))$, where ${\bf D} (\rho(g))$ is the domain of definition for $\rho(g)$.
\item[\hss\llap{\rm LR2)}]  If $g_1,g_2,g_1g_2\in U$ and $u\in{\bf D}$ then
$\rho(g_2)u \in{\bf D} (\rho(g_1))$ and
\[
\rho(g_1)[\rho(g_2)u] = \rho(g_1g_2)u\, .
\]
\item[\hss\llap{\rm LR3)}]  Let $Y\in{\frak h}$ such that $\exp tY\in U$
for $0\le t\le1$. Then for every $u\in{\bf D}$
\[
\lim_{t\to0} \rho(\exp tY) u = u\, .
\]
\item[\hss\llap{\rm LR4)}]  $\rho(Y){\bf D}\subset{\bf D}$ for every $Y\in
{\frak h}$.
\item[\hss\llap{\rm LR5)}]  $\forall u\in{\bf D}\, \exists V_u$ an open
$1$-neighborhod
in $H$ such that $UV_u\subset U^2$ and $\rho(h)u \in{\bf D}$
for every $h\in V_u$.
\item[\hss\llap{\rm LR6)}]  For every $Y\in{\frak q}$ and every $u\in{\bf D}$
the function
\[
h\longmapsto\rho(\exp(\operatorname{Ad}(h)Y))u
\]
is locally integrable on $\{h\in H\mid\exp(\operatorname{Ad}(h)Y)\in U\}$.
\end{list}
\cite{Jor86} now states that every local representation extends to a unitary
representation of $G^{c}$. We now want to use Theorem \ref{Dense} to construct
a local representation of $G$. For that let $0<R<1$ and let $\mathbf{D}%
=\mathbf{K}(B_{R}(0))$. Let $V$ be a symmetric open neighborhood of $1\in G$
such that $V\cdot B_{R}(0)\subset\Omega$. Let $U_{1}$ be a convex symmetric
neighborhood of $0$ in $\frak{g}$ such that with $U:=\exp U_{1}$ we have
$U^{2}\subset V$. If $g\in U$ then obviously (LR1)--(LR3) are satisfied. (LR4)
is satisfied as differentiation does not increase support. (LR6) is also clear
as $u=\beta(f)$ with $f\in\mathcal{C}_{c}^{\infty}(U)$ and hence $\left\|
\rho_{c}(\exp\operatorname{Ad}(h)Y)u\right\|  $ is continuous as a function of
$h$.

(LR5) Let $u=\beta(f)\in\mathbf{K}(B_{R}(0))$. Let $L=\operatorname*{supp}%
(f)\subset B_{R}(0)$. Let $V_{u}$ be such that $V_{u}^{-1}=V_{u}$,
$V_{u}L\subset B_{R}(0)$, and $V_{u}\subset U$. Then $UV_{u}\subset U^{2}$ and
$\tilde{\pi}(\nu)(h)u=\beta(\pi(\nu)(h)f)$ is defined and in $D$. This now
implies that $\tilde{\pi}$ restricted to $U$ is a local representation. Hence
the existence of $\tilde{\pi}^{c}$ follows from \cite{Jor86}. We notice that
this construction of $\tilde{\pi}^{c}$ does not use the full semigroup $S$ but
only $H$ and $\exp\mathbb{R}_{+}X^{o}$.
%
%
%
%
%

\begin{remark}
\label{RemEnright}Here we have used reflection positivity to go from the
generalized principal series to a highest weight module. Similar, but purely
algebraic construction for the special case $G_{\mathbb{C}}$ and
$G_{\mathbb{C}}^{c}=G\times G$ is due to T. J. Enright, \cite{E83}.
\end{remark}

\section{\label{S:Bargmann}The Segal-Bargmann Transform}

We have seen that the reflection positivity and the Osterwalder-Schrader
duality both correspond to transforming functions on the real form
$\Omega\subset\Omega_{\mathbb{C}}$ to a reproducing Hilbert space of
holomorphic functions on $\Omega_{\mathbb{C}}$. A classical example of a
similar situation is given by the Segal-Bargmann transform of the
Schr\"{o}dinger representation realized on $\mathbf{L}^{2}(\mathbb{R}^{n})$ to
the Fock model realized in the space of holomorphic function on $\mathbb{C}%
^{n}$ with finite norm%
\[
\left\|  F\right\|  _{\mathbf{F}}^{2}=\int_{\mathbb{C}^{n}}\left|
F(z)\right|  ^{2}\,e^{-\left|  z\right|  ^{2}}\,d\mu(z)
\]
where $d\mu$ is the $\pi^{-n}$ times the Lebesgue measure on $\mathbb{C}^{n}$.
We normalize the Lebesgue measure on $\mathbb{R}^{n}$ by $(2\pi)^{-n/2}$ times
the usual Lebesgue measure. We denote the resulting measure by $dx$. There are
no further constants in the Fourier inversion formula in this normalization.

Our first observation is that $\mathbb{C}^{n}$ is a ``complexification''\ of
$\mathbb{R}^{n}$. Therefore the ``restriction''\ map%
\[
R\colon\mathbf{F}(\mathbb{C}^{n})\longrightarrow\mathbf{L}^{2}(\mathbb{R}%
^{n})\,,\qquad RF(x):=e^{-x^{2}/2}F(x)
\]
is injective. It can be shown that $R$ is continuous with dense image. Hence
$R^{\ast}\colon\mathbf{L}^{2}(\mathbb{R}^{n})\rightarrow\mathbf{F}%
(\mathbb{C}^{n})$ is well defined and continuous.

The next observation is, that $\mathbf{F}(\mathbb{C}^{n})$ is a reproducing
Hilbert space with reproducing kernel $K(w,z)=K_{z}(w)=e^{w\overline{z}}$.
Therefore%
\begin{align}
R^{\ast}f(z)  &  =%
\ip{K_z}{R^*f}%
\nonumber\\
&  =%
\ip{RK_z}{f}%
\nonumber\\
&  =\int f(x)e^{-x^{2}/2}e^{zx}\,dx\nonumber\\
&  =e^{z^{2}/2}\int f(x)e^{-(z-x)^{2}/2}\,dx\nonumber\\
&  =e^{z^{2}/2}H\ast f(z) \label{E:R*classical}%
\end{align}
where $H(z)=e^{-z^{2}/2}$.

The third observation is that this is just a special case of the \emph{heat
semigroup}%
\[
H_{t}(f)(y)=H_{t}\ast f(y)=t^{-n/2}\int f(x)e^{-(y-x)^{2}/2t}\,dx\,.
\]
As the name indicate the family $\left\{  H_{t}\right\}  _{t>0}$ is a
semigroup, that is
\begin{equation}
H_{t}\ast H_{s}=H_{t+s}. \label{E:Heatsemigr}%
\end{equation}
This gives us the following form for $RR^{\ast}$ and $\left|  R^{\ast}\right|
=\sqrt{RR^{\ast}}$:%
\begin{align}
RR^{\ast}h(x)  &  =e^{-x^{2}/2}R^{\ast}h(x)\nonumber\\
&  =H_{1}\ast h(x)\,. \label{E:RR*classical}%
\end{align}
Thus%
\begin{align}
\sqrt{RR^{\ast}}(h)(x)  &  =H_{1/2}\ast h(x)\label{E:SqrRR*classical}\\
&  =2^{n/2}\int h(y)e^{-(x-y)^{2}}\,dy\,.\nonumber
\end{align}
{}From this we derive the following expression for the unitary part of the
polar decomposition of $R^{\ast}=B\sqrt{RR^{\ast}}$ of $R^{\ast}$:%
\begin{align*}
Bh(z)  &  =R^{\ast}\left|  R^{\ast}\right|  ^{-1}h(z)\\
&  =e^{z^{2}/2}H_{1/2}\ast h(z)\\
&  =2^{n/2}\int h(x)e^{-x^{2}+2zx-z^{2}/2}\,dy\\
&  =2^{n/2}e^{-z^{2}/2}\int h(x)e^{-x^{2}+2xz}\,dx
\end{align*}
which is, up to a scaling factor, the usual Bargmann transform, see
\cite{GF89}, p. 40. In particular it follows directly from our construction
that the Bargmann transform is a unitary isomorphism.

This example shows that we can recover the Bargmann transform from the
\emph{restriction principle}, see \cite{OO96}, that is we have

\begin{itemize}
\item  Manifolds $M\subset M_{\mathbb{C}}$ where $M_{\mathbb{C}}$ is a
``complexification'' of $M$;

\item  Groups $H\subset G$ such that $H$ acts on $M$ and $G$ acts on
$M_{\mathbb{C}}$;

\item  Measures $\mu$ and $\lambda$ on $M$ respectively $M_{\mathbb{C}}$ with
the measure $\mu$ being $H$-invariant;

\item  A reproducing Hilbert space $\mathbf{F}(M_{\mathbb{C}})\subset
\mathbf{L}^{2}(M_{\mathbb{C}},\lambda)\cap\mathcal{O}(M_{\mathbb{C}})$ with a
representation $\pi_{G}$ of $G$ given by%
\[
\lbrack\pi_{G}(g)F](z)=m(g^{-1},z)^{-1}F(g^{-1}\cdot z)
\]
where $m$ is a ``multiplier''.

\item  A function $\chi$ such that the ``restriction''\ map $R(F)(x):=\chi
(x)F(x)$, $x\in M$, from $\mathbf{F}(M_{\mathbb{C}})\rightarrow\mathbf{L}%
^{2}(M,\mu)$ is a closed (or continuous) $H$-intertwining operator with dense image.
\end{itemize}

Denote the reproducing kernel of $\mathbf{F}(M_{\mathbb{C}})$ by
$K(z,w)=K_{w}(z)$. Then $K_{w}\in\mathbf{F}(M_{\mathbb{C}})$ and
$K(z,w)=\overline{K(w,z)}$. The map $R^{\ast}\colon L^{2}(M,\mu)\rightarrow
\mathbf{F}(M_{\mathbb{C}})$ has the form%
\begin{align}
R^{\ast}h(w)  &  =%
\ip{K_w}{R^*h}%
\nonumber\\
&  =%
\ip{RK_w}{h}%
\nonumber\\
&  =\int_{M}h(m)\overline{\chi(m)}K(w,m)\,d\mu(m)\,. \label{E:R*}%
\end{align}
Hence%
\begin{equation}
RR^{\ast}h(x)=\int_{M}h(m)\chi(x)\overline{\chi(m)}\,K(x,m)\,d\mu(m)\,.
\label{E:RestrPrinc}%
\end{equation}
Write $R^{\ast}=B\left|  R^{\ast}\right|  $ for the polar decomposition of
$R^{\ast}$.

\begin{definition}
The unitary isomorphism $B\colon\mathbf{L}^{2}(M,\mu)\rightarrow
\mathbf{F}(M_{\mathbb{C}})$ is called the (generalized) Bargmann transform.
\end{definition}

The natural setting that we are looking at now is $H/H\cap K\simeq
\Omega\subset G^{c}/K^{c}=\Omega_{\mathbb{C}}$ and one of the highest weight
modules as a generalization of the Fock spaces with the representation
$[\rho_{\nu}(g)F](Z)=J_{\lambda}(g^{-1},Z)^{-1}F(g^{-1}\cdot Z)$, see
(\ref{E:intrep}). Here $\lambda$ corresponds to $\nu-\rho$ in our previous
sections. Define%
\begin{equation}
RF(x)=J_{\lambda}(x,0)^{-1}F(x\cdot0),\quad x\in H\mathbb{\,}.
\label{E:DefRbd}%
\end{equation}
Using the multiplier relation $J_{\lambda}(h^{-1}x,0)=J_{\lambda}%
(h^{-1},x\cdot0)J_{\lambda}(x,0)$ which follows from (\ref{E:Multipl}) \ we
get:%
\begin{align*}
R\left[  \rho_{\lambda}(h)F\right]  (x)  &  =J_{\lambda}(x,0)^{-1}\left[
\rho_{\lambda}(h)F\right]  (x\cdot0)\\
&  =J_{\lambda}(x,0)^{-1}J_{\lambda}(h^{-1},x\cdot0)^{-1}F(h^{-1}\cdot
(x\cdot0))\\
&  =J_{\lambda}(h^{-1}x,0)^{-1}F((h^{-1}x)\cdot0)\\
&  =[L(h)RF](x)
\end{align*}
where $L$ stands for the left regular representation of $H$ on $H/H\cap K$.
Hence $R$ is an intertwining operator. Let $\frak{a}_{p}$, $\Delta$ and
$\Delta^{+}$ be as in the end of Section \ref{S:Sssp}. By \cite{'OO88a,'OO88b}
the following is known:

\begin{theorem}
Suppose that $%
\left\langle \lambda+\rho,H_{\alpha}\right\rangle
<0$ for all $\alpha\in\Delta_{+}$. Then $\rho_{\lambda}$ is equivalent to a
discrete summand in $\mathbf{L}^{2}(G^{c}/H)$, that is there exists an
injective $G^{c}$-map $T\colon\mathbf{H}(\rho_{\lambda})\rightarrow
\mathbf{L}^{2}(G^{c}/H)$.
\end{theorem}

We call the resulting discrete part of $\mathbf{L}^{2}(G^{c}/H)$ the
\emph{holomorphic discrete series} of $G^{c}/H$. We have (see \cite{OO96}):

\begin{theorem}
Assume that $\rho_{\lambda}$ is a holomorphic discrete series representation
of $G^{c}/H$ then the following holds:

\begin{enumerate}
\item  The restriction map is injective, closed and with dense image.

\item  The generalized Bargmann transform $B\colon\mathbf{L}^{2}(H/H\cap
K)\rightarrow\mathbf{H}(\rho_{\lambda})$ is a $H$-isomorphism.
\end{enumerate}
\end{theorem}

Denote as before the reproducing kernel of $\mathbf{H}(\rho_{\lambda})$ by
$Q(W,Z)$. Then by (\ref{E:R*}) and Theorem \ref{S:DF1}:\
\begin{align}
R^{\ast}f(Z)  &  =\int_{H/H\cap K}f(h\cdot0)\cdot\overline{J_{\lambda
}(h,0)^{-1}Q(h\cdot0,Z)}\,dh\nonumber\\
&  =\int_{H/H\cap K}f(h\cdot0)J_{\lambda}(h^{-1},Z)\,dh\,. \label{E:R*general}%
\end{align}
If $Z=x\cdot0$, $x\in H$, then $J_{\lambda}(h^{-1},Z)=J_{\lambda}%
(h^{-1}x,0)J(x,0)^{-1}$. Let $\Psi_{\lambda}(h):=J_{\lambda}(h,0)$ and view
$f(h\cdot0)$ as a right $H\cap K$-invariant function on $H$. Then $RR^{\ast}$
is the convolution operator%
\[
RR^{\ast}f(x)=\int_{H}f(h)\Psi_{\lambda}(h^{-1}x)\,dh\,.
\]
The problem that we face here is, that $\lambda\mapsto\Psi_{\lambda}\ast$ is
\emph{not} a semigroup of operators. Hence it becomes harder to find the
square root of $RR^{\ast}$, but it can still be shown, that it is a
convolution operator. The final remark in this section is the following
connection between our map $U$ that comes from the reflection positivity and
$R^{\ast}$.

\begin{theorem}
Let $f\in\mathcal{C}_{c}^{\infty}(\Omega,\nu)$. Define a function $F$ on
$H/H\cap K$ by $F(h)=a_{\bar{N}}(h)^{-\nu-\rho}f(h\cdot_{\operatorname*{opp}%
}0)$. Then%
\[
Uf(Z)=R^{\ast}(F)(Z)\,.
\]
\end{theorem}

\begin{proof}
This follows from (\ref{E:R*general}) and (\ref{E:U})
\end{proof}

\section{\label{Diagonal}The Heisenberg Group}

\setcounter{equation}{0}

A special case of the setup in Definition \ref{ReflectionSymmetric} above
arises as follows: Let the group $G$, and $\tau\in\operatorname{Aut}_{2}(G)$
be as described there. Let $\mathbf{H}_{\pm}$ be two given complex Hilbert
spaces, and $\pi_{\pm}\in\operatorname{Rep}(G,\mathbf{H}_{\pm})$ a pair of
unitary representations. Suppose $T\colon\mathbf{H}_{-}\rightarrow
\mathbf{H}_{+}$ is a unitary operator such that $T\pi_{-}=\left(  \pi_{+}%
\circ\tau\right)  T$, or equivalently,
\begin{equation}
T\pi_{-}(g)f_{-}=\pi_{+}\left(  \tau(g)\right)  Tf_{-} \label{TPiMinus}%
\end{equation}
for all $g\in G$, and all $f_{-}\in\mathbf{H}_{-}$. Form the direct sum
$\mathbf{H}:=\mathbf{H}_{+}\oplus\mathbf{H}_{-}$ with inner product
\begin{equation}%
\ip{f_{+}\oplus f_{-}}{f^{\prime}_{+}\oplus f^{\prime
}_{-}}%
:=%
\ip{f_{+}}{f^{\prime}_{+}}%
_{+}+%
\ip{f_{-}}{f^{\prime}_{-}}%
_{-} \label{IPPlusMinus}%
\end{equation}
where the $\pm$ subscripts are put in to refer to the respective Hilbert
spaces $\mathbf{H}_{\pm}$, and we may form $\pi:=\pi_{+}\oplus\pi_{-}$ as a
unitary representation on $\mathbf{H}=\mathbf{H}_{+}\oplus\mathbf{H}_{-}$ by
\[
\pi(g)\left(  f_{+}\oplus f_{-}\right)  =\pi_{+}(g)f_{+}\oplus\pi_{-}%
(g)f_{-}\,,\quad g\in G,\;f_{\pm}\in\mathbf{H}_{\pm}\,.
\]
Setting
\begin{equation}
J:=\left(
\begin{matrix}
0 & T\\
T^{\ast} & 0
\end{matrix}
\right)  \,, \label{JMatrix}%
\end{equation}
that is $J\left(  f_{+}\oplus f_{-}\right)  =\left(  Tf_{-}\right)
\oplus\left(  T^{\ast}f_{+}\right)  $, it is then clear that properties
(1)--(2) from Definition \ref{ReflectionSymmetric} will be satisfied for the
pair $(J,\pi)$. Formula (\ref{TPiMinus}) may be recovered by writing out the
relation
\begin{equation}
J\pi=\left(  \pi\circ\tau\right)  J \label{JPiRelation}%
\end{equation}
in matrix form, specifically
\[
\left(
\begin{matrix}
0 & T\\
T^{\ast} & 0
\end{matrix}
\right)  \left(
\begin{matrix}
\pi_{+}(g) & 0\\
0 & \pi_{-}(g)
\end{matrix}
\right)  =\left(
\begin{matrix}
\pi_{+}(\tau(g)) & 0\\
0 & \pi_{-}(\tau(g))
\end{matrix}
\right)  \left(
\begin{matrix}
0 & T\\
T^{\ast} & 0
\end{matrix}
\right)  \,.
\]
If, conversely, (\ref{JPiRelation}) is assumed for some unitary period-$2$
operator $J$ on $\mathbf{H}=\mathbf{H}_{+}\oplus\mathbf{H}_{-}$, and, if the
two representations $\pi_{+}$ and $\pi_{-}$ are \textit{disjoint,} in the
sense that no irreducible in one occurs in the other (or, equivalently, there
is no nonzero intertwiner between them), then, in fact, (\ref{TPiMinus}) will
follow from (\ref{JPiRelation}). The diagonal terms in (\ref{JMatrix}) will be
zero if (\ref{JPiRelation}) holds. This last implication is an application of
Schur's lemma.

\begin{lemma}
\label{Operator}Let $0\neq\mathbf{K}_{0}$ be a closed linear subspace of
$\mathbf{H}=\mathbf{H}_{+}\oplus\mathbf{H}_{-}$ satisfying the positivity
condition \textup{(PR3)} in Definition \textup{\ref{Hyperbolic},} that is
\begin{equation}%
\ip{v}{Jv}%
\geq0\,,\quad\forall v\in\mathbf{K}_{0} \label{vJvPositivity}%
\end{equation}
where
\begin{equation}
J=\left(
\begin{matrix}
0 & T\\
T^{\ast} & 0
\end{matrix}
\right)  \label{JMatrixBis}%
\end{equation}
is given from a fixed unitary isomorphism $T\colon\mathbf{H}_{-}%
\rightarrow\mathbf{H}_{+}$ as in \textup{(\ref{TPiMinus}).} For $v=f_{+}\oplus
f_{-}\in\mathbf{H}=\mathbf{H}_{+}\oplus\mathbf{H}_{-}$, set $P_{+}v:=f_{+}$.
The closure of the subspace $P_{+}\mathbf{K}_{0}$ in $\mathbf{H}_{+}$ will be
denoted $\overline{P_{+}\mathbf{K}_{0}}$. Then the subspace
\[
\mathbf{G}=\left\{  \left.  \left(
\begin{matrix}
f_{+}\\
f_{-}%
\end{matrix}
\right)  \in\mathbf{K}_{0}\,\right|  \,f_{-}\in T^{\ast}\left(  \overline
{P_{+}\mathbf{K}_{0}}\right)  \right\}
\]
is the graph of a closed linear operator $M$ with domain
\begin{equation}
\mathbf{D}=\left\{  f_{+}\in\mathbf{H}_{+}\,\left|  \,\exists f_{-}\in
T^{\ast}\left(  \overline{P_{+}\mathbf{K}_{0}}\right)  \mathop{{\rm s.t.}%
}\left(
\begin{matrix}
f_{+}\\
f_{-}%
\end{matrix}
\right)  \in\mathbf{K}_{0}\right.  \right\}  \,; \label{MDomain}%
\end{equation}
and, moreover, the product operator $L:=TM$ is dissipative on this domain,
that is
\begin{equation}%
\ip{Lf_{+}}{f_{+}}%
_{+}+%
\ip{f_{+}}{Lf_{+}}%
_{+}\geq0 \label{LDissipative}%
\end{equation}
holds for all $f_{+}\in\mathbf{D}$.
\end{lemma}

\begin{proof}
The details will only be sketched here, but the reader is referred to
\cite{Sto51} and \cite{Jor80} for definitions and background literature. An
important argument in the proof is the verification that, if a column vector
of the form $\left(
\begin{matrix}
0\\
f_{-}%
\end{matrix}
\right)  $ is in $\mathbf{G}$, then $f_{-}$ must necessarily be zero in
$\mathbf{H}_{-}$. But using positivity, we have
\begin{equation}
\left|
\ip{u}{Jv}%
\right|  ^{2}\leq%
\ip{u}{Ju}%
\ip{v}{Jv}%
\,,\quad\forall u,v\in\mathbf{K}_{0}\,. \label{uJvBound}%
\end{equation}
Using this on the vectors $u=\left(
\begin{matrix}
0\\
f_{-}%
\end{matrix}
\right)  $ and $v=\left(
\begin{matrix}
k_{+}\\
k_{-}%
\end{matrix}
\right)  \in\mathbf{K}_{0}$, we get
\[%
\ip{\left(
\begin{matrix}
0  \\   f_{-}
\end{matrix}
\right)  }{\left(
\begin{matrix}
Tk_{-}  \\   T^{*}k_{+}
\end{matrix}
\right)  }%
=%
\ip{f_{-}}{T^{*}k_{+}}%
=0\,,\quad\forall k_{+}=P_{+}v\,.
\]
But, since $f_{-}$ is also in $T^{\ast}\left(  \overline{P_{+}\mathbf{K}_{0}%
}\right)  $, we conclude that $f_{-}=0$, proving that $\mathbf{G}$ is the
graph of an operator $M$ as specified. The dissipativity of the operator
$L=TM$ is just a restatement of \textup{(PR3)}.
\end{proof}

The above result involves only the operator-theoretic information implied by
the data in Definition \ref{Hyperbolic}, and, in the next lemma, we introduce
the representations:

\begin{lemma}
\label{Normal}Let the representations $\pi_{\pm}$ and the intertwiner $T$ be
given as specified before. Let $H=G^{\tau}$; and suppose we have a cone
$C\subset\frak{q}$ as specified in \textup{(PR1), (PR2}$^{\prime}$\textup{)}
and \textup{(PR3}$^{\prime}$\textup{).} Assume further that\begin{list}%
{}{\setlength{\leftmargin}{\customleftmargin}
\setlength{\itemsep}{0.5\customskipamount}
\setlength{\parsep}{0.5\customskipamount}
\setlength{\topsep}{\customskipamount}}
\item[\hss\llap{\rm PR4)}]  $\mathbf{D}$ is dense in $\mathbf{H}_{+}$;
\item[\hss\llap{\rm PR5)}]  The commutant of $\pi(H)$ is abelian.
\end{list}
Then $L=TM$ is normal.
\end{lemma}

\begin{proof}
Since $T$ is a unitary isomorphism $\mathbf{H}_{-}\rightarrow\mathbf{H}_{+}$
we may make an identification and reduce the proof to the case where
$\mathbf{H}_{+}=\mathbf{H}_{-}$ and $T$ is the identity operator. We then
have
\[
\pi_{-}=T^{-1}\left(  \pi_{+}\circ\tau\right)  T=\pi_{+}\circ\tau\,;
\]
and if $h\in H$, then
\[
\pi_{-}(h)=\pi_{+}\left(  \tau(h)\right)  =\pi_{+}(h)\,;
\]
while, if $\tau(g)=g^{-1}$, then
\[
\pi_{-}(g)=\pi_{+}\left(  \tau(g)\right)  =\pi_{+}\left(  g^{-1}\right)  \,.
\]
Using only the $H$ part from (PR2${}^{\prime}$), we conclude that
$\mathbf{K}_{0}$ is invariant under $\pi_{+}\oplus\pi_{+}(H)$. If the
projection $P_{\mathbf{K}_{0}}$ of $\mathbf{H}_{+}\oplus\mathbf{H}_{+}$ onto
$\mathbf{K}_{0}$ is written as an operator matrix $\left(
\begin{matrix}
P_{11} & P_{12}\\
P_{21} & P_{22}%
\end{matrix}
\right)  $ with entries representing operators in $\mathbf{H}_{+}$, and
satisfying
\begin{align*}
P_{11}^{\ast}  &  =P_{11}\,,\\
P_{22}^{\ast}  &  =P_{22}\,,\\
P_{12}^{\ast}  &  =P_{21}\,,\\
P_{ij}  &  =P_{i1}P_{1j}+P_{i2}P_{2j}\,,
\end{align*}
then it follows that
\begin{equation}
P_{ij}\pi_{+}(h)=\pi_{+}(h)P_{ij}\quad\forall i,j=1,2,\;\forall h\in H\,,
\label{PijPiPlus}%
\end{equation}
which puts each of the four operators $P_{ij}$ in the commutant $\pi
_{+}(H)^{\prime}$ from (PR5). Using (PR4), we then conclude that $L$ is a
dissipative operator with $\mathbf{D}$ as dense domain, and that
$\mathbf{K}_{0}$ is the graph of this operator. Using (PR5), and a theorem of
Stone \cite{Sto51}, we finally conclude that $L$ is a normal operator, that is
it can be represented as a multiplication operator with dense domain
$\mathbf{D}$ in $\mathbf{H}_{+}$.
\end{proof}

We shall consider two cases below (the Heisenberg group, and the
$(ax+b)$-group) when conditions (PR4)--(PR5) can be verified from the context
of the representations. Suppose $G$ has two abelian subgroups $H$, $N$, and
the second $N$ also a normal subgroup, such that $G=HN$ is a product
representation in the sense of Mackey \cite{Mac}. Define $\tau\in
\operatorname{Aut}_{2}(G)$ by setting
\begin{equation}
\tau(h)=h\,,\quad\forall h\in H\,,\text{ and }\tau(n)=n^{-1}\,,\quad\forall
n\in N\,. \label{TauDef}%
\end{equation}

The Heisenberg group is a copy of $\mathbb{R}^{3}$ represented as matrices
$\vphantom{\left( \begin{matrix} 1 & a & c  \\  0 & 1 & b  \\  0 & 0 &
1\end{matrix}\right) _X}\left(
\begin{matrix}
1 & a & c\\
0 & 1 & b\\
0 & 0 & 1
\end{matrix}
\right)  $, or equivalently vectors $(a,b,c)\in\mathbb{R}^{3}$. Setting
$H=\left\{  (a,0,0)\mid a\in\mathbb{R}\right\}  $ and
\begin{equation}
N=\left\{  (0,b,c)\mid b,c\in\mathbb{R}\right\}  \,, \label{NHeisenberg}%
\end{equation}
we arrive at one example.

The $(ax+b)$-group is a copy of $\mathbb{R}^{2}$ represented as matrices
$\left(
\begin{matrix}
a & b\\
0 & 1
\end{matrix}
\right)  $, $a=e^{s}$, $b\in\mathbb{R}$, $s\in\mathbb{R}$. Here we may take
$H=\left\{  \left.  \left(
\begin{matrix}
a & 0\\
0 & 1
\end{matrix}
\right)  \,\right|  \,a\in\mathbb{R}_{+}\right\}  $ and
\begin{equation}
N=\left\{  \left.  \left(
\begin{matrix}
1 & b\\
0 & 1
\end{matrix}
\right)  \,\right|  \,b\in\mathbb{R}\right\}  \,, \label{Naxb}%
\end{equation}
and we have a second example of the Mackey factorization. Generally, if $G=HN$
is specified as described, we use the representations of $G$ which are induced
from one-dimensional representations of $N$. If $G$ is the Heisenberg group,
or the $(ax+b)$-group, we get all the infinite-dimensional irreducible
representations of $G$ by this induction (up to unitary equivalence, of
course). For the Heisenberg group, the representations are indexed by
$\hslash\in\mathbb{R}\setminus\{0\}$, $\hslash$ denoting Planck's constant.
The representation $\pi_{\hslash}$ may be given in $\mathbf{H}=\mathbf{L}%
^{2}(\mathbb{R})$ by
\begin{equation}
\pi_{\hslash}(a,b,c)f(x)=e^{i\hslash(c+bx)}f(x+a)\,,\quad\forall
f\in\mathbf{L}^{2}(\mathbb{R}),\;(a,b,c)\in G\,. \label{PiHBar}%
\end{equation}
The Stone-von Neumann uniqueness theorem asserts that every unitary
representation $\pi$ of $G$ satisfying
\[
\pi(0,0,c)=e^{i\hslash c}I_{\mathbf{H}(\pi)}\quad(\hslash\neq0)
\]
is unitarily equivalent to a direct sum of copies of the representation
$\pi_{\hslash}$ in (\ref{PiHBar}).

The $(ax+b)$-group (in the form $\left\{  \left.  \left(
\begin{matrix}
e^{s} & b\\
0 & 1
\end{matrix}
\right)  \,\right|  \,s,b\in\mathbb{R}\right\}  $) has only two inequivalent
unitary irreducible representations, and they may also be given in the same
Hilbert space $\mathbf{L}^{2}(\mathbb{R})$ by
\begin{equation}
\pi_{\pm}\left(
\begin{matrix}
e^{s} & b\\
0 & 1
\end{matrix}
\right)  f(x)=e^{\pm ie^{x}b}f(x+s)\,,\quad\forall f\in\mathbf{L}%
^{2}(\mathbb{R})\,. \label{PiPlusMinusaxb}%
\end{equation}
There are many references for these standard facts from representation theory;
see, e.g., \cite{Jor88}.

\begin{lemma}
\label{Abel}Let the group $G$ have the form $G=HN$ for locally compact abelian
subgroups $H,N$, with $N$ normal, and $H\cap N=\{e\}$. Let $\chi$ be a
one-dimensional unitary representation of $N$, and let $\pi
=\operatorname*{ind}{}_{N}^{G}(\chi)$ be the corresponding induced
representation. Then the commutant of $\left\{  \pi(H)\mid h\in H\right\}  $
is an abelian von Neumann algebra: in other words, condition \textup{(PR5)} in
Lemma \textup{\ref{Normal}} is satisfied.
\end{lemma}

\begin{proof}
See, e.g., \cite{Jor88}.
\end{proof}

In the rest of the present section, we will treat the case of the
\emph{Heisenberg group,} and the $(ax+b)$\emph{-group} will be the subject of
the next section.

For both groups we get pairs of unitary representations $\pi_{\pm}$ arising
from some $\tau\in\operatorname{Aut}_{2}(G)$ and described as in
(\ref{JPiRelation}) above. But when the two representations $\pi_{+}$ and
$\pi_{-}=\pi_{+}\circ\tau$ are irreducible and disjoint, we will show that
there are no spaces $\mathbf{K}_{0}$ satisfying (PR1), (PR2$^{\prime}$), and
(PR3) such that $\mathbf{K}=\left(  \mathbf{K}_{0}/\mathbf{N}\right)
\sptilde$ is nontrivial. Here (PR2) is replaced by \begin{list}{}%
{\setlength{\leftmargin}{\customleftmargin}
\setlength{\itemsep}{0.5\customskipamount}
\setlength{\parsep}{0.5\customskipamount}
\setlength{\topsep}{\customskipamount}}
\item[\hss\llap{\rm PR2$^{\prime}$)}]  $C$ is a nontrivial cone in $\frak
{q} $.
\end{list}
Since for both groups, and common to all the representations, we noted that
the Hilbert space $\mathbf{H}_{+}$ may be taken as $\mathbf{L}^{2}%
(\mathbb{R})$, we can have $J$ from (\ref{JMatrixBis}) represented in the form
$J=\left(
\begin{matrix}
0 & I\\
I & 0
\end{matrix}
\right)  $. Then the $J$-inner product on $\mathbf{H}_{+}\oplus\mathbf{H}%
_{-}=\mathbf{L}^{2}(\mathbb{R})\oplus\mathbf{L}^{2}(\mathbb{R})\simeq
\mathbf{L}^{2}(\mathbb{R},\mathbb{C}^{2})$ may be brought into the form
\begin{equation}%
\ip{\left( \begin{matrix} f_{+} \\ f_{-}\end{matrix}\right) }{\left
( \begin{matrix} f_{+} \\ f_{-} \end{matrix}\right) }%
_{J}=2\operatorname{Re}%
\ip{f_{+}}{f_{-}}%
=2\int_{-\infty}^{\infty}\operatorname{Re}\left(  \overline{f_{+}(x)}%
f_{-}(x)\right)  \,dx\,. \label{JInnerProduct}%
\end{equation}

For the two examples, we introduce
\[
N_{+}=\left\{  (0,b,c)\mid b,c\in\mathbb{R}_{+}\right\}
\]
where $N$ is defined in (\ref{NHeisenberg}), but $N_{+}$ is not $H$-invariant.
Alternatively, set
\[
N_{+}=\left\{  \left.  \left(
\begin{matrix}
1 & b\\
0 & 1
\end{matrix}
\right)  \,\right|  \,b\in\mathbb{R}_{+}\right\}
\]
for the alternative case where $N$ is defined from (\ref{Naxb}), and note that
this $N_{+}$ is $H$-invariant. In fact there are the following $4$ invariant
cones in $\frak{q}$:
\begin{align*}
C_{1}^{+}  &  =\left\{  (0,0,t)\mid t\geq0\right\} \\
C_{1}^{-}  &  =\left\{  (0,0,t)\mid t\leq0\right\} \\
C_{2}^{+}  &  =\left\{  (0,x,y)\mid x\in\mathbb{R},\;y\geq0\right\} \\
C_{2}^{-}  &  =\left\{  (0,x,y)\mid x\in\mathbb{R},\;y\leq0\right\}
\end{align*}
Let $\pi$ denote one of the representations of $G=HN$ from the discussion
above (see formulas (\ref{PiHBar}) and (\ref{PiPlusMinusaxb})) and let
$\mathcal{D}$ be a closed subspace of $\mathbf{H}=\mathbf{L}^{2}(\mathbb{R})$
which is assumed invariant under $\pi(HN_{+})$. Then it follows that the two
spaces
\begin{align}
\mathcal{D}_{\infty}  &  :=\bigvee\left\{  \pi(n)\mathcal{D}\mid n\in
N\right\} \label{DInfinity}\\
\mathcal{D}_{-\infty}  &  :=\bigwedge\left\{  \pi(n)\mathcal{D}\mid n\in
N\right\}  \label{DMinusInfinity}%
\end{align}
are invariant under $\pi(G)$, where the symbols $\bigvee$ and $\bigwedge$ are
used for the usual lattice operations on closed subspaces in $\mathbf{H}$. We
leave the easy verification to the reader, but the issue is resumed in the
next section. If $P_{\infty}$, resp., $P_{-\infty}$, denotes the projection of
$\mathbf{H}$ onto $\mathcal{D}_{\infty}$, resp., $\mathcal{D}_{-\infty}$, then
we assert that both projections $P_{\pm\infty}$ are in the commutant of
$\pi(G)$. So, if $\pi$ is irreducible, then each $P_{\infty}$, or $P_{-\infty
}$, must be $0$ or $I$. Since $\mathcal{D}_{-\infty}\subset\mathcal{D}%
\subset\mathcal{D}_{\infty}$ from the assumption, it follows that $P_{\infty
}=I$ if $\mathcal{D}\neq\{0\}$.

\begin{lemma}
\label{Heisenberg}Let $G$ be the Heisenberg group, and let the notation be as
described above. Let $\pi_{+}$ be one of the representations $\pi_{\hslash}$
and let $\pi_{-}$ be the corresponding $\pi_{-\hslash}$ representation. Let
$0\neq\mathbf{K}_{0}\subset\mathbf{L}^{2}(\mathbb{R})\oplus\mathbf{L}%
^{2}(\mathbb{R})$ be a closed subspace which is invariant under $\left(
\pi_{+}\oplus\pi_{-}\right)  \left(  HN_{+}\right)  $. Then it follows that
there are only the following possibilities for $\overline{P_{+}\mathbf{K}_{0}%
}$: $\{0\}$, $\mathbf{L}^{2}(\mathbb{R})$, or $A\mathcal{H}_{+}$ where
$\mathcal{H}_{+}$ denotes the Hardy space in $\mathbf{L}^{2}(\mathbb{R})$
consisting of functions $f$ with Fourier transform $\hat{f}$ supported in the
half-line $\left[  0,\infty\right)  $, and where $A\in\mathbf{L}^{\infty
}(\mathbb{R})$ is such that $\left|  A(x)\right|  =1\mathop{{\rm a.e.}}%
x\in\mathbb{R}$. For the space $\overline{P_{-}\mathbf{K}_{0}}$, there are the
possibilities: $\{0\}$, $\mathbf{L}^{2}(\mathbb{R})$, and $A\mathcal{H}_{-}$,
where $A$ is a \textup{(}possibly different\/\textup{)} unitary $\mathbf{L}%
^{\infty}$-function, and $\mathcal{H}_{-}$ denotes the negative Hardy space.
\end{lemma}

\begin{proof}
Immediate from the discussion, and the Beurling-Lax theorem classifying the
closed subspaces in $\mathbf{L}^{2}(\mathbb{R})$ which are invariant under the
multiplication operators, $f(x)\mapsto e^{iax}f(x)$, $a\in\mathbb{R}_{+}$. We
refer to \cite{LaPh}, or \cite{Hel64}, for a review of the Beurling-Lax theorem.
\end{proof}

\begin{corollary}
\label{Positive}Let $\pi_{\pm}$ be the representations of the Heisenberg
group, and suppose that the subspace $\mathbf{K}_{0}$ from Lemma
\textup{\ref{Heisenberg}} is chosen such that \textup{(PR1)--(PR3)} in
Definition \textup{\ref{Hyperbolic}} hold. Then $\left(  \mathbf{K}%
_{0}/\mathbf{N}\right)  \sptilde=\{0\}$.
\end{corollary}

\begin{proof}
Suppose there are unitary functions $A_{\pm}\in\mathbf{L}^{\infty}%
(\mathbb{R})$ such that $\overline{P_{\pm}\mathbf{K}_{0}}=A_{\pm}%
\mathcal{H}_{\pm}$. Then this would violate the Schwarz-estimate
(\ref{uJvBound}), and therefore condition (PR3). Using irreducibility of
$\pi_{+}=\pi_{\hslash}$ and of $\pi_{-}=\pi_{+}\circ\tau=\pi_{-\hslash}$, we
may reduce to considering the cases when one of the spaces $\overline{P_{\pm
}\mathbf{K}_{0}}$ is $\mathbf{L}^{2}(\mathbb{R})$. By Lemma \ref{Normal}, we
are then back to the case when $\mathbf{K}_{0}$ or $\mathbf{K}_{0}^{-1}$ is
the graph of a densely defined normal and dissipative operator $L$, or
$L^{-1}$, respectively. We will consider $L$ only. The other case goes the
same way. Since
\begin{equation}
\left(  \pi_{+}\oplus\pi_{-}\right)  \left(  0,b,0\right)  \left(  f_{+}\oplus
f_{-}\right)  (x)=e^{i\hslash bx}f_{+}(x)\oplus e^{-i\hslash bx}f_{-}(x)
\label{Pibf}%
\end{equation}
it follows that $L$ must anti-commute with the multiplication operator $ix$ on
$\mathbf{L}^{2}(\mathbb{R})$. For deriving this, we used assumption (PR3) at
this point. We also showed in Lemma \ref{Normal} that $L$ must act as a
multiplication operator on the Fourier-transform side. But the
anti-commutativity is inconsistent with a known structure theorem in
\cite{Ped90}, specifically Corollary 3.3 in that paper. Hence there are
unitary functions $A_{\pm}$ in $\mathbf{L}^{\infty}(\mathbb{R})$ such that
$\overline{P_{\pm}\mathbf{K}_{0}}=A_{\pm}\mathcal{H}_{\pm}$. But this
possibility is inconsistent with positivity in the form $\operatorname{Re}%
\ip{f_{+}}{f_{-}}%
\geq0\,,\quad\forall(f_{+},f_{-})\in\mathbf{K}_{0}$ (see (\ref{JInnerProduct}%
)) if $\left(  \mathbf{K}_{0}/\mathbf{N}\right)  \sptilde\neq\{0\}$. To see
this, note that $\mathbf{K}_{0}$ is invariant under the unitary operators
(\ref{Pibf}) for $b\in\mathbb{R}_{+}$. The argument from Lemma
\ref{Heisenberg}, now applied to $\pi_{+}\oplus\pi_{-}$, shows that the two
subspaces
\[
\mathbf{K}_{0}^{\infty}:=\bigvee_{b\in\mathbb{R}}\left(  \pi_{+}\oplus\pi
_{-}\right)  \left(  0,b,0\right)  \mathbf{K}_{0}%
\]
and
\[
\mathbf{K}_{0}^{-\infty}:=\bigwedge_{b\in\mathbb{R}}\left(  \pi_{+}\oplus
\pi_{-}\right)  \left(  0,b,0\right)  \mathbf{K}_{0}%
\]
are both invariant under the whole group $\left(  \pi_{+}\oplus\pi_{-}\right)
(G)$. But the commutant of this is $2$-dimensional: the only projections in
the commutant are represented as one of the following,
\[
\left(
\begin{matrix}
0 & 0\\
0 & 0
\end{matrix}
\right)  ,\;\left(
\begin{matrix}
I & 0\\
0 & 0
\end{matrix}
\right)  ,\;\left(
\begin{matrix}
0 & 0\\
0 & I
\end{matrix}
\right)  ,\text{\quad or\quad}\left(
\begin{matrix}
I & 0\\
0 & I
\end{matrix}
\right)  \,,
\]
relative to the decomposition $\mathbf{L}^{2}(\mathbb{R})\oplus\mathbf{L}%
^{2}(\mathbb{R})$ of $\pi_{+}\oplus\pi_{-}$. The above analysis of the
anti-commutator rules out the cases $\left(
\begin{matrix}
I & 0\\
0 & 0
\end{matrix}
\right)  $ and $\left(
\begin{matrix}
0 & 0\\
0 & I
\end{matrix}
\right)  $, and if $\left(  \mathbf{K}_{0}/\mathbf{N}\right)  \sptilde
\neq\{0\}$, we are left with the cases $\mathbf{K}_{0}^{\infty}=\{0\}$ and
$\mathbf{K}_{0}^{\infty}=\mathbf{L}^{2}(\mathbb{R})\oplus\mathbf{L}%
^{2}(\mathbb{R})$. Recall, generally $\mathbf{K}_{0}^{-\infty}\subset
\mathbf{K}_{0}\subset\mathbf{K}_{0}^{\infty}$, as a starting point for the
analysis. A final application of the Beurling-Lax theorem (as in \cite{LaPh};
see also \cite{DyMc70}) to (\ref{Pibf}) then shows that there must be a pair
of unitary functions $A_{\pm}$ in $\mathbf{L}^{\infty}(\mathbb{R})$ such that
\begin{equation}
\mathbf{K}_{0}=A_{+}\mathcal{H}_{+}\oplus A_{-}\mathcal{H}_{-} \label{KAAHH}%
\end{equation}
where $\mathcal{H}_{\pm}$ are the two Hardy spaces given by having $\hat{f}$
supported in $\left[  0,\infty\right)  $, respectively, $\left(
-\infty,0\right]  $. The argument is now completed by noting that
(\ref{KAAHH}) is inconsistent with the positivity of $\mathbf{K}_{0}$ in
(\ref{vJvPositivity}); that is, we clearly do not have $%
\ip{\left(
\begin{matrix}
A_{+}h_{+}  \\   A_{-}h_{-}
\end{matrix}
\right)  }{J\left(
\begin{matrix}
A_{+}h_{+}  \\   A_{-}h_{-}
\end{matrix}
\right)  }%
=2\operatorname{Re}%
\ip{A_{+}h_{+}}{A_{-}h_{-}}%
$ semidefinite, for all $h_{+}\in\mathcal{H}_{+}$ and all $h_{-}\in
\mathcal{H}_{-}$. This concludes the proof of the Corollary.
\end{proof}

\begin{remark}
\label{UncorrelatedClosedSubspace}At the end of the above proof of Corollary
\ref{Positive}, we arrived at the conclusion (\ref{KAAHH}) for the subspace
$\mathbf{K}_{0}$ under consideration. Motivated by this, we define a closed
subspace $\mathbf{K}_{0}$ in a direct sum Hilbert space $\mathbf{H}_{+}%
\oplus\mathbf{H}_{-}$ to be \emph{uncorrelated} if there are closed subspaces
$\mathbf{D}_{\pm}\subset\mathbf{H}_{\pm}$ in the respective summands such
that
\begin{equation}
\mathbf{K}_{0}=\mathbf{D}_{+}\oplus\mathbf{D}_{-} \label{KDD}%
\end{equation}
Contained in the corollary is then the assertion that every
semigroup-invariant $\mathbf{K}_{0}$ in $\mathbf{L}^{2}(\mathbb{R}%
)\oplus\mathbf{L}^{2}(\mathbb{R})$ is uncorrelated, where the semigroup here
is the subsemigroup $S$ in the Heisenberg group $G$ given by
\begin{equation}
S=\left\{  (a,b,c)\mid b\in\mathbb{R}_{+},\;a,c\in\mathbb{R}\right\}  \,,
\label{SinG}%
\end{equation}
and the parameterization is the one from (\ref{NHeisenberg}). We also had the
representation $\pi$ in the form $\pi_{+}\oplus\pi_{-}$ where the respective
summand representations $\pi_{\pm}$ of $G$ are given by (\ref{PiHBar})
relative to a pair $(\hslash,-\hslash)$, $\hslash\in\mathbb{R}\setminus\{0\}$
some fixed value of Planck's constant. In particular, it is assumed in
Corollary \ref{Positive} that each representation $\pi_{\pm}$ is
\emph{irreducible.} But for proving that some given semigroup-invariant
$\mathbf{K}_{0}$ must be uncorrelated, this last condition can be relaxed
considerably; and this turns out to be relevant for applications to
Lax-Phillips scattering theory for the wave equation with obstacle scattering
\cite{LaPh}. In that context, the spaces $\mathbf{D}_{\pm}$ will be outgoing,
respectively, incoming subspaces; and the wave equation translates backwards,
respectively forwards, according to the unitary one-parameter groups $\pi
_{-}(0,b,0)$, respectively, $\pi_{+}(0,b,0)$, with $b\in\mathbb{R}$
representing the time-variable $t$ for the unitary time-evolution
one-parameter group which solves the wave equation under consideration. The
unitary-equivalence identity (\ref{JPiRelation}) stated before Lemma
\ref{Operator} then implies equivalence of the wave-dynamics before, and
after, the obstacle scattering.
\end{remark}

Before stating our next result, we call attention to the $(2n+1)$-dimensional
Heisenberg group $G_{n}$ in the form $\mathbb{R}^{2n+1}=\mathbb{R}^{n}%
\times\mathbb{R}^{n}\times\mathbb{R}$, in parameter form: $a,b\in
\mathbb{R}^{n}$, $c\in\mathbb{R}$, and product rule
\[
(a,b,c)\cdot(a^{\prime},b^{\prime},c^{\prime})=(a+a^{\prime},b+b^{\prime
},c+c^{\prime}+a\cdot b^{\prime})
\]
where $a+a^{\prime}=(a_{1}+a_{1}^{\prime},\dots,a_{n}+a_{n}^{\prime})$ and
$a\cdot b^{\prime}=\sum_{j=1}^{n}a_{j}b_{j}^{\prime}$. For every (fixed)
$b\in\mathbb{R}^{n}\setminus\{0\}$, we then have a subsemigroup
\begin{equation}
S(b)=\left\{  (a,\beta b,c)\mid\beta\in\mathbb{R}_{+},\;a\in\mathbb{R}%
^{n},\;c\in\mathbb{R}\right\}  \,; \label{SubsemigroupS}%
\end{equation}
and we show in the next result that it is enough to have invariance under such
a semigroup in $G_{n}$, just for a single direction, defined from some fixed
$b\in\mathbb{R}^{n}\setminus\{0\}$.

\begin{theorem}
\label{KNoughtUncorrelated}Let $\pi_{\pm}$ be unitary representations of the
Heisenberg group $G$ on respective Hilbert spaces $\mathbf{H}_{\pm}$, and let
$T\colon\mathbf{H}_{-}\rightarrow\mathbf{H}_{+}$ be a unitary isomorphism
which intertwines $\pi_{-}$ and $\pi_{+}\circ\tau$ as in
\textup{(\ref{TPiMinus})} where
\begin{equation}
\tau(a,b,c)=(a,-b,-c)\,,\quad\forall(a,b,c)\in G\simeq\mathbb{R}^{2n+1}\,.
\label{TauInHeisenberg}%
\end{equation}
Suppose there is $\hslash\in\mathbb{R}\setminus\{0\}$ such that
\begin{equation}
\pi_{+}(0,0,c)=e^{i\hslash c}I_{\mathbf{H}_{+}}\,. \label{PiPlusInHeisenberg}%
\end{equation}
If $\mathbf{K}_{0}\subset\mathbf{H}_{+}\oplus\mathbf{H}_{-}$ is a closed
subspace which is invariant under
\[
\left\{  (\pi_{+}\oplus\pi_{-})(a,\beta b,c)\mid a\in\mathbb{R}^{n},\;\beta
\in\mathbb{R}_{+},\;c\in\mathbb{R}\right\}
\]
from \textup{(\ref{SubsemigroupS}),} $b\in\mathbb{R}^{n}\setminus\{0\}$, then
we conclude that $\mathbf{K}_{0}$ must automatically be uncorrelated.
\end{theorem}

\begin{proof}
The group-law in the Heisenberg group yields the following commutator rule:
\[
(a,0,0)(0,b,0)(-a,0,0)=(0,b,a\cdot b)
\]
for all $a,b\in\mathbb{R}^{n}$. We now apply $\pi=\pi_{+}\oplus\pi_{-}$ to
this, and evaluate on a general vector $f_{+}\oplus f_{-}\in\mathbf{K}%
_{0}\subset\mathbf{H}_{+}\oplus\mathbf{H}_{-}$: abbreviating $\pi(a)$ for
$\pi(a,0,0)$, and $\pi(b)$ for $\pi(0,b,0)$, we get
\[
\pi(a)\pi(\beta b)\pi(-a)(f_{+}\oplus f_{-})=e^{i\hslash\beta a\cdot b}\pi
_{+}(\beta b)f_{+}\oplus e^{-i\hslash\beta a\cdot b}\pi_{-}(\beta b)f_{-}%
\in\mathbf{K}_{0}%
\]
valid for all $a\in\mathbb{R}^{n}$, $\beta\in\mathbb{R}_{+}$. Note, in
(\ref{PiPlusInHeisenberg}), we are assuming that $\pi_{+}$ takes on some
specific value $e^{i\hslash c}$ on the one-dimensional center. Since $\pi_{-}$
is unitarily equivalent to $\pi_{+}\circ\tau$ by assumption (see
(\ref{PiPlusInHeisenberg})), we conclude that
\[
\pi_{-}(0,0,c)=e^{-i\hslash c}I_{\mathbf{H}_{-}}\,,\quad\forall c\in
\mathbb{R}\,.
\]
The argument really only needs that the two representations $\pi_{\pm}$ define
\emph{different} characters on the center. (Clearly $\hslash\neq-\hslash$
since $\hslash\neq0$.) Multiplying through first with $e^{-i\hslash\beta
a\cdot b}$, and integrating the resulting term
\[
\pi_{+}(\beta b)f_{+}\oplus e^{-i2\hslash\beta a\cdot b}\pi_{-}(\beta
b)f_{-}\in\mathbf{K}_{0}%
\]
in the $a$-variable, we get $\pi_{+}(\beta b)f_{+}\oplus0\in\mathbf{K}_{0}$.
The last conclusion is just using that $\mathbf{K}_{0}$ is a closed subspace.
But we can do the same with the term
\[
e^{i2\hslash\beta a\cdot b}\pi_{+}(\beta b)f_{+}\oplus\pi_{-}(\beta b)f_{-}%
\in\mathbf{K}_{0}\,,
\]
and we arrive at $0\oplus\pi_{-}(\beta b)f_{-}\in\mathbf{K}_{0}$. Finally
letting $\beta\rightarrow0_{+}$, and using strong continuity, we get
$f_{+}\oplus0$ and $0\oplus f_{-}$ both in $\mathbf{K}_{0}$. Recalling that
$f_{\pm}$ are general vectors in $P_{\pm}\mathbf{K}_{0}$, we conclude that
$P_{+}\mathbf{K}_{0}\oplus P_{-}\mathbf{K}_{0}\subset\mathbf{K}_{0}$, and
therefore $\overline{P_{+}\mathbf{K}_{0}}\oplus\overline{P_{-}\mathbf{K}_{0}%
}\subset\mathbf{K}_{0}$. Since the converse inclusion is obvious, we arrive at
(\ref{KDD}) with $\mathbf{D}_{\pm}=\overline{P_{\pm}\mathbf{K}_{0}}$.
\end{proof}

The next result shows among other things that there are representations $\pi$
of the Heisenberg group $G_{n}$ (for each $n$) such that the reflected
representation $\pi^{c}$ of $G_{n}^{c}\simeq G_{n}$ (see Theorem
\ref{PiCIrreducible}) acts on a nonzero Hilbert space $\mathbf{H}^{c}=\left(
\mathbf{K}_{0}/\mathbf{N}\right)  \sptilde$. However, because of Lemma
\ref{KernelPiC}, $\pi^{c}\left(  G_{n}^{c}\right)  $ will automatically be an
\emph{abelian} group of operators on $\mathbf{H}^{c}$. To see this, note that
the proof of Theorem \ref{MDomain} shows that $\pi^{c}$ must act as the
identity operator on $\mathbf{H}^{c}$ when restricted to the one-dimensional
center in $G_{n}^{c}\simeq G_{n}$.

It will be convenient for us to read off this result from a more general
context: we shall consider a general Lie group $G$, and we fix a
right-invariant Haar measure on $G$.

A distribution $F$ on the Lie group $G$ will be said to be positive definite
(PD) if
\begin{equation}
\int_{G}\int_{G}F(uv^{-1})\overline{f(u)}f(v)\,du\,dv\geq0 \tag*{(PD)}%
\end{equation}
for all $f\in C_{c}^{\infty}(G)$; and we say that $F$ is PD on some open
subset $\Omega\subset G$ if this holds for all $f\in C_{c}^{\infty}(\Omega)$.
The interpretation of the expression in (PD) is in the sense of distributions.
But presently measurable functions $F$ will serve as the prime examples.

We say that the distribution is reflection-positive (RP) on $\Omega$
((RP$_{\Omega}$) for emphasis) if, for some period-$2$ automorphism $\tau$ of
$G$, we have
\begin{equation}
F\circ\tau=F \label{FTau}%
\end{equation}
and
\begin{equation}
\int_{G}\int_{G}F(\tau(u)v^{-1})\overline{f(u)}f(v)\,du\,dv\geq0 \tag*{(RP$
_\Omega$)}%
\end{equation}
for all $f\in C_{c}^{\infty}(\Omega)$.

We say that $x\in G$ is (RP$_{\Omega}$)-contractive if (RP$_{\Omega}$) holds,
and for all $f\in C_{c}^{\infty}(\Omega)$%
\begin{align*}
0  &  \leq\int_{G}\int_{G}F(\tau(u)v^{-1})\overline{f(ux)}f(vx)\,du\,dv\\
&  \leq\int_{G}\int_{G}F(\tau(u)v^{-1})\overline{f(u)}f(v)\,du\,dv\,.
\end{align*}
Note that, since
\[
\int_{G}\int_{G}F(\tau(u)v^{-1})\overline{f(ux)}f(vx)\,du\,dv=\int_{G}\int
_{G}F(\tau(u)\tau(x)^{-1}xv^{-1})\overline{f(u)}f(v)\,du\,dv\,
\]
it follows that every $x\in H$ is contractive; in fact, isometric. If instead
we have $\tau(x)=x^{-1}$, then contractivity amounts to the estimate
\[
\int_{G}\int_{G}F(\tau(u)x^{2}v^{-1})\overline{f(u)}f(v)\,du\,dv\leq\int
_{G}\int_{G}F(\tau(u)v^{-1})\overline{f(u)}f(v)\,du\,dv\,\,
\]
for all $f\in C_{c}^{\infty}(\Omega)$. Using the Basic Lemma one can also show
that $x$ acts by contractions.

The following result is useful, but an easy consequence of the definitions and
standard techniques for positive definite distributions; see for example
\cite{Jor88,Pra89}.

\begin{theorem}
\label{PiExtends}Let $F$ be a distribution on a Lie group $G$ with a
period-$2$ automorphism $\tau$, and suppose $F$ is $\tau$-invariant,
\textup{(PD)} holds on $G$, and \textup{(RP}$_{\Omega}$\textup{)} holds on
some open, and semigroup-invariant, subset $\Omega$ in $G$. Then define
\[
(\pi(u)f)(v):=f(vu)\,,\quad\forall u,v\in G,\;\forall f\in C_{c}^{\infty
}(G)\,;
\]
and
\[
Jf:=f\circ\tau\,.
\]
Let $\mathbf{H}(F)$ be the Hilbert space obtained from the GNS construction,
applied to \textup{(PD),} with inner product on $C_{c}^{\infty}(G)$ given by
\[%
\ip{f}{g}%
:=\int_{G}\int_{G}F(uv^{-1})\overline{f(u)}g(v)\,du\,dv\,.
\]
Then $\pi$ extends to a unitary representation of $G$ on $\mathbf{H}(F)$, and
$J$ to a unitary operator, such that
\[
J\pi=(\pi\circ\tau)J\,.
\]
If \textup{(RP}$_{\Omega}$\textup{)} further holds, as described, then $\pi$
induces \textup{(}via Theorem \textup{\ref{PiCIrreducible})} a unitary
representation $\pi^{c}$ of $G^{c}$ acting on the new Hilbert space
$\mathbf{H}^{c}$ obtained from completing in the new inner product from
\textup{(RP}$_{\Omega}$\textup{),} and dividing out with the corresponding kernel.
\end{theorem}

The simplest example of a function $F$ on the Heisenberg group $G_{n}$
satisfying (PD), but not (RP$_{\Omega}$), for nontrivial $\Omega$'s, may be
obtained from the Green's function for the sub-Laplacian on $G_{n}$; see
\cite[p. 599]{Ste93} for details.

If complex coordinates are introduced in $G_{n}$, the formula for $F$ takes
the following simple form: let $z\in\mathbb{C}^{n}$, $c\in\mathbb{R}$, and
define
\[
F(z,c)=\frac{1}{\left(  \left|  z\right|  ^{4}+c^{2}\right)  ^{n}}\,.
\]
Then we adapt the product in $G_{n}$ to the modified definition as follows:
\[
(z,c)\cdot(z^{\prime},c^{\prime})=(z+z^{\prime},c+c^{\prime}+\left\langle
z,z^{\prime}\right\rangle )\quad\forall z,z^{\prime}\in\mathbb{C}%
^{n},\;\forall c,c^{\prime}\in\mathbb{R}\,,
\]
where $\left\langle z,z^{\prime}\right\rangle $ is the symplectic form
\[
\left\langle z,z^{\prime}\right\rangle :=2\operatorname{Im}(z\cdot\bar
{z}^{\prime})\,.
\]
The period-$2$ automorphism $\tau$ on $G_{n}$ we take as
\[
\tau(z,c)=(\bar{z},-c)
\]
with $\bar{z}$ denoting complex conjugation $(z_{1},\dots,z_{n})\mapsto
(\bar{z}_{1},\dots,\bar{z}_{n})$.

The simplest example where both (PD) and (RP$_{\Omega}$) hold on the
Heisenberg group $G_{n}$ is the following:

\begin{example}
\label{SubLaplacian}Let $\zeta=(\zeta_{1},\dots\zeta_{n})\in\mathbb{C}^{n}$,
$\xi_{j}=\operatorname{Re}\zeta_{j}$, $\eta_{j}=\operatorname{Im}\zeta_{j}$,
$j=1,\dots,n$. Define
\[
F(z,c)=\int_{\mathbb{R}^{2n}}\frac{e^{i\operatorname{Re}(z\cdot\bar{\zeta})}%
}{\prod_{j=1}^{n}(\left|  \zeta_{j}\right|  ^{2}+1)}\,d\xi_{1}\cdots\,d\xi
_{n}\,d\eta_{1}\cdots\,d\eta_{n}\,.
\]
Let $\Omega:=\{(z,c)\in G_{n}\mid z=(z_{j})_{j=1}^{n},\;\operatorname{Im}%
z_{j}>0\}$. Then (PD) holds on $G_{n}$, and (RP$_{\Omega}$) holds, referring
to this $\Omega$. Since the expression for $F(z,c)$ factors, the problem
reduces to the $(n=1)$ special case. There we have
\[
F(z,c)=\int_{\mathbb{R}^{2}}\frac{e^{i(x\xi+y\eta)}}{\xi^{2}+\eta^{2}+1}%
\,d\xi\,d\eta\,;
\]
and if $f\in C_{c}^{\infty}(\Omega)$ with $\Omega=\left\{  (z,c)\mid
y>0\right\}  $, then
\begin{multline*}
\int_{G_{1}}\int_{G_{1}}F(\tau(u)v^{-1})\overline{f(u)}f(v)\,du\,dv\\
=\int_{\mathbb{R}^{8}}\frac{e^{i(x-x^{\prime})\xi}\,e^{-i(y+y^{\prime})\eta}%
}{\xi^{2}+\eta^{2}+1}\overline{f(x+iy,c)}f(x^{\prime}+iy^{\prime},c^{\prime
})\,d\xi\,d\eta\,dx\,dy\,dc\,dx^{\prime}\,dy^{\prime}\,dc^{\prime}\,.
\end{multline*}
Let $\tilde{f}$ denote the Fourier transform in the $x$-variable, keeping the
last two variables $(y,c)$ separate. Then the integral transforms as follows:
\[
\int_{\mathbb{R}^{5}}\frac{e^{-(y+y^{\prime})\sqrt{1+\xi^{2}}}}{\sqrt
{1+\xi^{2}}}\overline{\tilde{f}(\xi,y,c)}\tilde{f}(\xi,y^{\prime},c^{\prime
})\,d\xi\,dy\,dy^{\prime}\,dc\,dc^{\prime}\,.
\]
Introducing the Laplace transform in the middle variable $y$, we then get
(since $f$ is supported in $y>0$)
\[
\int_{0}^{\infty}e^{-y\sqrt{1+\xi^{2}}}\tilde{f}(\xi,y,c)\,dy=\tilde
{f}_{\lambda}(\xi,\sqrt{1+\xi^{2}},c)\,;
\]
the combined integral reduces further:
\[
\int_{\mathbb{R}}\,\left|  \vphantom{\int\tilde{f}_\lambda(\xi,\sqrt{1+\xi
^2},c)\,dc}\smash{\int_{\mathbb{R}} \tilde{f}_\lambda(\xi,\sqrt{1+\xi
^2},c)\,dc}\right|  ^{2}\frac{d\xi}{1+\xi^{2}}%
\]
which is clearly positive; and we have demonstrated that (RP$_{\Omega}$)
holds. It is immediate that $F$ is $\tau$-invariant (see (\ref{FTau})), and
also that it satisfies (PD) on $G_{n}$.
\end{example}

\section{\label{axb}The $(ax+b)$-Group Revisited}

\setcounter{equation}{0}

We showed that in general we get a unitary representation $\pi^{c}$ of the
group $G^{c}$ from an old one $\pi$ of $G$, provided $\pi$ satisfies the
assumptions of reflection positivity. The construction as we saw uses a
certain cone $C$ and a semigroup $H\exp C$, which are part of the axiom
system. What results is a new class of unitary representations $\pi^{c}$
satisfying a certain spectrum condition (semi-bounded spectrum).

But, for the simplest non-trivial group $G$, this semi-boundedness turns out
\textit{not} to be satisfied in the general case. Nonetheless, we still have a
reflection construction getting us from unitary representations $\pi$ of the
$(ax+b)$-group, such that $\pi\circ\tau\simeq\pi$ (unitary equivalence), to
associated unitary representations $\pi^{c}$ of the same group. The (up to
conjugation) unique non-trivial period-$2$ automorphism $\tau$ of $G$, where
$G$ is the $(ax+b)$-group, is given by
\[
\tau(a,b)=(a,-b)\,.
\]
Recall that the $G$ may be identified with the matrix-group
\[
\left\{  \left.  \left(
\begin{matrix}
a & b\\
0 & 1
\end{matrix}
\right)  \,\right|  \,a>0,b\in\mathbb{R}\right\}
\]
and $(a,b)$ corresponds to the matrix $\left(
\begin{matrix}
a & b\\
0 & 1
\end{matrix}
\right)  $. In this realization the Lie algebra of $G$ has the basis
\[
X=\left(
\begin{matrix}
1 & 0\\
0 & 0
\end{matrix}
\right)  \quad\mbox{and}\quad Y=\left(
\begin{matrix}
0 & 1\\
0 & 0
\end{matrix}
\right)  \,.
\]
We have $\exp(tX)=(e^{t},0)$ and $\exp(sY)=(1,s)$. Hence $\tau(X)=X$ and
$\tau(Y)=-Y$. Thus $\frak{h}=\mathbb{R}X$ and $\frak{q}=\mathbb{R}Y$. We
notice the commutator relation $[X,Y]=Y$. The possible $H$-invariant cones in
$\frak{q}$ are $\pm\{tY\mid t\geq0\}$. It is known from Mackey's theory that
$G$ has two inequivalent, unitary, irreducible, infinite-dimensional
representations $\pi_{\pm}$, and it is immediate that we have the unitary
equivalence (see details below):%
\begin{equation}
\pi_{+}\circ\tau\simeq\pi_{-}\,. \label{E:5.1}%
\end{equation}

Hence, if we set $\pi:=\pi_{+}\oplus\pi_{-}$, then $\pi\circ\tau\simeq\pi$, so
we have the setup for the general theory. We show that $\pi$ may be realized
on $\mathbf{L}^{2}(\mathbb{R})\oplus\mathbf{L}^{2}(\mathbb{R})\simeq
\mathbf{L}^{2}(\mathbb{R},\mathbb{C}^{2})$, and we find and classify the
invariant positive subspaces $\mathbf{K}_{0}\subset\mathbf{L}^{2}%
(\mathbb{R},\mathbb{C}^{2})$. To understand the interesting cases for the
$(ax+b)$-group $G$, we need to relax the invariance condition: We shall
\textit{not} assume invariance of $\mathbf{K}_{0}$ under the semigroup
$\{\pi(1,b)\mid b\geq0\}$, but only under the infinitesimal unbounded
generator $\pi(Y)$. With this, we still get the correspondence $\pi\mapsto
\pi_{\mathbf{K}_{0}}^{c}$ as described above.

We use the above notation. We know from Mackey's theory \cite{Mac} that there
are two inequivalent irreducible infinite-dimensional representations of $G$,
and we shall need them in the following alternative formulations: Let
$\mathcal{L}_{\pm}$ denote the respective Hilbert space $\mathbf{L}%
^{2}(\mathbb{R}_{\pm})$ with the multiplicative invariant measure $d\mu_{\pm
}=dp/|p|$, $p\in\mathbb{R}_{\pm}$. Then the formula%
\begin{equation}
f\longmapsto e^{ipb}f(pa) \label{E:5.4}%
\end{equation}
for functions $f$ on $\mathbb{R}$ restricts to two unitary irreducible
representations, denoted by $\pi_{\pm}$ of $G$ on the respective spaces
$\mathcal{L}_{\pm}$. Let $Q(f)(p):=f(-p)$ denote the canonical mapping from
$\mathcal{L}_{+}$ to $\mathcal{L}_{-}$, or equivalently from $\mathcal{L}_{-}$
to $\mathcal{L}_{+}$. Then we have for $g\in G$ (cf.\ (\ref{E:5.1})):%
\begin{equation}
Q\pi_{+}(g)=\pi_{-}(\tau(g))Q \label{E:5.5}%
\end{equation}

For the representation $\pi:=\pi_{+}\oplus\pi_{-}$ on $\mathbf{H}%
:=\mathcal{L}_{+}\oplus\mathcal{L}_{-}$ we therefore have%
\begin{equation}
J\pi(g)=\pi(\tau(g))J,\quad g\in G\,, \label{E:5.6}%
\end{equation}
where $J$ is the unitary involutive operator on $\mathbf{H}$ given by%
\begin{equation}
J=\left(
\begin{matrix}
0 & Q\\
Q & 0
\end{matrix}
\right)  \,. \label{E:5.7}%
\end{equation}

Instead of the above $p$-realization of $\pi$ we will mainly use the following
$x$-formalism. The map $t\mapsto\pm e^{t}$ defines an isomorphism $L_{\pm
}\colon\mathcal{L}_{\pm}\rightarrow\mathbf{L}^{2}(\mathbb{R})$, where we use
the additive Haar measure $dx$ on $\mathbb{R}$. For $g=(e^{s},b)\in G$ and
$f\in\mathbf{L}^{2}(\mathbb{R})$, set%
\begin{equation}
(\pi_{\pm}(g)f)(x):=e^{\pm ie^{x}b}f(x+s),\quad x\in\mathbb{R}\,.
\label{E:5.8}%
\end{equation}

A simple calculation shows that $L_{\pm}$ intertwines the old and new
construction of $\pi_{\pm}$, excusing our abuse of notation. In this
realization $Q$ becomes simply the identity operator $Q(f)(x)=f(x)$. The
involution $J\colon\mathbf{L}^{2}(\mathbb{R},\mathbb{C}^{2})$ is now simply
given by
\begin{equation}
J(f_{0},f_{1})=(f_{1},f_{0}) \label{eqaxbNew.7}%
\end{equation}
or $J=\left(
\begin{matrix}
0 & 1\\
1 & 0
\end{matrix}
\right)  $.

In this formulation the operator
\begin{equation}
L:=\pi_{\pm}(\Delta_{H}-\Delta_{Q})=\pi_{\pm}(X^{2}-Y^{2}) \label{E:5.9}%
\end{equation}
takes the form
\begin{equation}
L=\left(  \frac{d\,}{dx}\right)  ^{2}\,+\,e^{2x}\,, \label{E:5.10}%
\end{equation}
but it is on $\mathbf{L}^{2}(\mathbb{R})$ and $-\infty<x<\infty$. This
operator is known to have defect indices $(1,0)$ \cite{Jor75,NeSt59}, which
means that it cannot be extended to a selfadjoint operator on $\mathbf{L}%
^{2}(\mathbb{R})$. Using a theorem from \cite{Jor75,ReSi75} we can see this by
comparing the quantum mechanical problem for a particle governed by $-L$ as a
Schr\"{o}dinger operator (that is a strongly repulsive force) with the
corresponding classical one governed (on each energy surface) by
\[
E_{\mathrm{kin}}+E_{\mathrm{pot}}=\left(  \frac{dx}{dt}\right)  ^{2}%
\,-\,e^{2x}=E\,.
\]
The escape time for this particle to $x=\pm\infty$ is
\begin{equation}
t_{\pm}=\int_{\mathrm{finite}}^{\pm\infty}\frac{dx}{\sqrt{E+e^{2x}}\,}\,,
\label{E:5.11}%
\end{equation}
that is $t_{\infty}$ is finite, and $t_{-\infty}=\infty$. We elaborate on this
point below. The nonzero defect vector for the quantum mechanical problem
corresponds to a boundary condition at $x=\infty$ since this is the
singularity which is reached in finite time.

The fact from \cite{Jor75} we use for the defect index assertion is this: The
Schr\"{o}dinger operator $H=-\left(  \frac{d\, }{dx}\right)  ^{2}\, + V(x)$
for a single particle has nonzero defect solutions $f_{\pm}\in\mathbf{L}
^{2}(\mathbb{R} )$ to $H^{*}f_{\pm}=\pm if_{\pm}$ iff there are solutions
$t\mapsto x(t) $ to the corresponding classical problem
\[
E=\left(  \frac{dx(t)}{dt}\right)  ^{2}+V(x(t))
\]
with finite travel-time to $x=+\infty$, respectively, $x=-\infty$. The
respective (possibly infinite) travel-times are
\[
t_{\pm\infty} = \int_{\mathrm{finite}}^{\pm\infty}\frac{dx}{\sqrt{E-V(x)}}\,
.
\]
The correspondence principle states that one finite travel-time to $+\infty$
(say) yields a dimension in the associated defect space, and similarly for the
other travel-time to $-\infty$.

In the $x$-formalism, (\ref{E:5.5}) from above then simplifies to the
following identity for operators on the \textit{same} Hilbert space
$\mathbf{L}^{2}(\mathbb{R})$ (carrying the two inequivalent representations
$\pi_{+}$ and $\pi_{-}$):
\begin{equation}
\pi_{+}(g)=\pi_{-}(\tau(g)),\quad g\in G\,. \label{E:5.12}%
\end{equation}
We realize the representation $\pi=\pi_{+}\oplus\pi_{-}$ in the Hilbert space
$\mathbf{H}=\mathbf{L}^{2}(\mathbb{R})\oplus\mathbf{L}^{2}(\mathbb{R}%
)=\mathbf{L}^{2}(X_{2})$ where $X_{2}={0}\times\mathbb{R}\cup{1}%
\times\mathbb{R}$. We may represent $J$ by an automorphism $\theta\colon
X_{2}\rightarrow X_{2}$ (as illustrated in Proposition \ref{P:3.3}):
\[
\theta(0,x):=(1,x)\quad\mbox{and}\quad\theta(1,y)=(0,y)\,,\quad x,y\in
\mathbb{R}\,,
\]
and
\[
J(f)(\omega)=f(\theta(\omega))\,,\quad\omega\in X_{2}\,.
\]
Notice that the subset
\[
X_{2}^{\theta}=\{\omega\in X_{2}\mid\theta(\omega)=\omega\}
\]
is empty. Define for $f\in\mathbf{L}^{2}(X_{2})$, $f_{k}(x)=f(k,x)$, $k=0,1$,
$x\in\mathbb{R}$. We have for $g=(e^{s},b)\in G$:
\[
\left(  \pi(g)f\right)  _{0}(x)=e^{ibe^{x}}f_{0}(x+s)=(\bar{\pi}_{+}%
(g)f_{0})(x)
\]
and
\[
\left(  \pi(g)f\right)  _{1}(x)=e^{-ibe^{x}}f_{1}(x+s)=(\bar{\pi}_{-}%
(g)f_{1})(x)\,.
\]

\begin{proposition}
\label{Negative}Let $\pi=\pi_{+}\oplus\pi_{-}$ be the representation from
\textup{(\ref{E:5.1})--(\ref{E:5.6})} above of the $(ax+b)$-group $G$. Then
the only choices of reflections $\mathbf{K}_{0}$ as in Remark
\textup{\ref{WhatToAssume}} for the sub-semigroup $S=\left\{  \left(
a,b\right)  \in G\mid b>0\right\}  $ will have $\mathbf{K}=\left(
\mathbf{K}_{0}/\mathbf{N}\right)  \sptilde$ equal to $0$.
\end{proposition}

\begin{lemma}
\label{QField}Let $Q$ be the projection in $\mathbf{L}^{2}(\mathbb{R}%
)\oplus\mathbf{L}^{2}(\mathbb{R})$ onto a translation-invariant $J$-positive
subspace. Then $Q$ is represented by a measurable field of $2\times2$ complex
matrices $\mathbb{R}\ni\xi\mapsto\left(  Q_{ij}(\xi)\right)  _{ij=1}^{2}$ such
that $\left|  Q_{12}(\xi)\right|  ^{2}=Q_{11}(\xi)Q_{22}(\xi)$ a.e. on
$\mathbb{R}$, and $Q_{12}(\xi)+Q_{21}(\xi)\geq0$ a.e.; and conversely.
\end{lemma}

\begin{corollary}
\label{QMatrix} These relations imply the following for the matrix $Q$:

\begin{enumerate}
\item [\hss\llap{\rm1)}]If $Q_{12}(\xi)=0$ then we have the three
possibilities:
\begin{align*}
Q(\xi)  &  =0\,,\\
Q(\xi)  &  =\left(
\begin{matrix}
1 & 0\\
0 & 0
\end{matrix}
\right)  \,,\mbox{ and }\\
Q(\xi)  &  =\left(
\begin{matrix}
0 & 0\\
0 & 1
\end{matrix}
\right)  \,.
\end{align*}
In all those cases, we have $Q(\xi)JQ(\xi)=0$.

\item[\hss\llap{\rm2)}] If $Q_{12}(\xi)\neq0$, then $0<Q_{22}(\xi
)=1-Q_{11}(\xi)<1$. Let $\mu(\xi)=Q_{12}(\xi)/Q_{11}(\xi)$. Then by
$\operatorname*{Tr}(Q(\xi)JQ(\xi))\geq0$ we have $\operatorname{Re}\mu
(\xi)\geq0$ and
\begin{equation}
Q(\xi)=\frac{1}{1+\left|  \mu(\xi)\right|  ^{2}}\left(
\begin{matrix}
1 & \mu(\xi)\\
\overline{\mu(\xi)} & \left|  \mu(\xi)\right|  ^{2}%
\end{matrix}
\right)  \,. \label{QMatrixFormula}%
\end{equation}
With $\lambda=\bar{\mu}$ we get that the image of $Q(\xi)$ is given by
\[
\left\{  \left.  u(\xi)\left(
\begin{matrix}
1\\
\lambda(\xi)
\end{matrix}
\right)  \,\right|  \,u(\xi)\in\mathbb{C}\right\}  \,.
\]
Specifying to our situation, $f=\left(
\begin{matrix}
f_{0}\\
f_{1}%
\end{matrix}
\right)  \in\mathbf{K}_{0}$ if and only if
\begin{equation}
\hat{f}_{1}(\xi)=\lambda(\xi)\hat{f}_{0}(\xi)\,. \label{flamf}%
\end{equation}
\end{enumerate}
\end{corollary}

Since $Q(\xi)$ is a measurable field of projections, the function $\mathbb{R}
\ni\xi\mapsto\lambda(\xi)$ must be measurable, but it may be unbounded. This
also means that $P_{\mathbf{K}_{0}}$ is the projection onto the graph of the
operator $T_{0}\colon f_{0}\mapsto f_{1}$ where $f_{0}$ and $f_{1}$ are
related as in (\ref{flamf}), and the Fourier transform $\hat{\cdot}$ is in the
$\mathbf{L}^{2}$-sense. \smallskip

The following argument deals with the general case, avoiding the separation of
the proof into the two cases (I) and (II): If vectors $v\in\mathbf{K}_{0}$ are
expanded as $v=\left(
\begin{matrix}
h\\
k
\end{matrix}
\right)  $, $h=Q_{11}h+Q_{12}k$, $k=Q_{21}h+Q_{22}k$, we can introduce
$\mathcal{D}=\left\{  h\in\mathbf{L}^{2}(\mathbb{R})\,\left|  \,\exists
k\in\mathbf{L}^{2}(\mathbb{R})\mathop{{\rm s.t.}}\left(
\begin{matrix}
h\\
k
\end{matrix}
\right)  \in\mathbf{N}\right.  \right\}  $. If $b>0$, then:
\begin{align*}
\pi_{+}(b)h  &  =Q_{11}\pi_{+}(b)h+Q_{12}\pi_{+}(-b)k\,,\\
\pi_{+}(-b)k  &  =Q_{21}\pi_{+}(b)h+Q_{22}\pi_{+}(-b)k\,,
\end{align*}
valid for any $\left(
\begin{matrix}
h\\
k
\end{matrix}
\right)  \in\mathbf{K}_{0}$, and $b\in\mathbb{R}_{+}$. So it follows from
Lemma \ref{BasicLemma} again that $\mathcal{D}$ is invariant under $\left\{
\pi_{+}(b)\mid b>0\right\}  $, and also under the whole semigroup $\left\{
\pi_{+}(g)\mid g\in S\right\}  $ where $\pi_{+}$ is now denoting the
corresponding irreducible representation of $G$ on $\mathbf{L}^{2}%
(\mathbb{R})$. Let
\begin{align}
\mathcal{D}_{\infty}  &  :=\bigvee_{b\in\mathbb{R}}\pi_{+}(b)\mathcal{D}%
\,,\label{DInfinitybis}\\
\mathcal{D}_{-\infty}  &  :=\bigwedge_{b\in\mathbb{R}}\pi_{+}(b)\mathcal{D}\,,
\label{DMinusInfinitybis}%
\end{align}
where $\bigvee$ and $\bigwedge$ denote the lattice operations on closed
subspaces in $\mathbf{L}^{2}(\mathbb{R})$, and
\[
\left(  \pi_{+}(b)f\right)  (x)=e^{ibe^{x}}f(x)\,,\quad f\in\mathbf{L}%
^{2}(\mathbb{R}),\;b,x\in\mathbb{R}.
\]
We may now apply the Lax-Phillips argument to the spaces $\mathcal{D}%
_{\pm\infty}$. If $\left(  \mathbf{K}_{0}/\mathbf{N}\right)  \sptilde$ should
be $\neq\{0\}$, then $\mathcal{D}=\{0\}$ by the argument. Since we are
assuming $\left(  \mathbf{K}_{0}/\mathbf{N}\right)  \sptilde
\neq\{0\}$, we get $\mathcal{D}=\{0\}$, and as a consequence the following
operator graph representation for $\mathbf{K}_{0}$: $\left(  \mathbf{K}%
_{0}/\mathbf{N}\right)  \sptilde=\beta\left(  G(L)\right)  $ where $G(L)$ is
the graph of a closed operator $L$ in $\mathbf{L}^{2}(\mathbb{R})$.
Specifically, this means that the linear mapping $\mathbf{K}_{0}/\mathbf{N}%
\ni\left(
\begin{matrix}
h\\
k
\end{matrix}
\right)  +\mathbf{N}\mapsto h$ is well-defined as a linear closed operator.
This in turn means that $\mathbf{K}_{0}$ may be represented as the graph of a
closable operator in $\mathbf{L}^{2}(\mathbb{R})$ as discussed in the first
part of the proof. Hence such a representation could have been assumed at the outset.

\begin{remark}
\label{BorchersCMP92}In a recent paper on local quantum field theory
\cite{Bor92}, Borchers considers in his Theorem II.9 a representation $\pi$ of
the $(ax+b)$-group $G$ on a Hilbert space $\mathbf{H}$ such that there is a
\emph{conjugate linear} $J$ (that is a period-$2$ antiunitary) such that $J\pi
J=\pi\circ\tau$ where $\tau$ is the period-$2$ automorphism of $G$ given by
$\tau(a,b):=(a,-b)$. In Borchers's example, the one-parameter subgroup
$b\mapsto\pi(1,b)$ has semibounded spectrum, and there is a unit-vector
$v_{0}\in\mathbf{H}$ such that $\pi(1,b)v_{0}=v_{0}$, $\forall b\in\mathbb{R}%
$. The vector $v_{0}$ is cyclic and separating for a von Neumann algebra $M$
such that $\pi(1,b)M\pi(1,-b)\subset M$, $\forall b\in\mathbb{R}_{+}$. Let
$a=e^{t}$, $t\in\mathbb{R}$. Then, in Borchers's construction, the other
one-parameter subgroup $t\mapsto\pi(e^{t},0)$ is the modular group
$\Delta^{it}$ associated with the cyclic and separating vector $v_{0}$ (from
Tomita-Takesaki theory \cite[vol. I]{BrRo}). Finally, $J$ is the corresponding
modular conjugation satisfying $JMJ=M^{\prime}$ when $M^{\prime}$ is the
commutant of $M$.
\end{remark}

\begin{acknowledgements}
The authors would like to thank Tom Branson, Paul Muhly, Bent \O rsted and
Steen Pedersen for helpful discussions, and the referee for constructive
suggestions. We would like to thank Brian Treadway for his excellent
typesetting of the final version of the paper in \LaTeX. Both authors were
also supported in part by the National Science Foundation.
\end{acknowledgements}

\end{document}